\documentclass[preprint,3p,sort&compress]{elsarticle}
\journal{{\tt arXiv}}

\usepackage[dvips]{epsfig}
\usepackage{graphicx} 		
\usepackage{pdfpages}
\usepackage{latexsym}
\usepackage{verbatim}
\usepackage{amsmath,amsthm}
\usepackage[utf8]{inputenc}
\usepackage{amssymb}
\usepackage{bbm}
\usepackage{fonts}
\usepackage{enumerate}
\usepackage{placeins}
\usepackage{standalone}
\usepackage{paralist}
\usepackage{subcaption}
\usepackage{booktabs}
\usepackage{array}
\newcolumntype{C}[1]{>{\centering\let\newline\\\arraybackslash\hspace{0pt}}m{#1}}
\usepackage{empheq}
\usepackage{algorithm, algorithmic}
\usepackage{mdframed}  

\usepackage[]{hyperref} 

\usepackage{pgfplots}
\pgfplotsset{compat=newest}       
\usepackage[]{externaltikz}

 \usepackage{relinput}
\newtheorem{theorem}{Theorem}[section]

\newtheorem{remark}[theorem]{Remark}

\newcounter{tikzsubfigcounter}[figure]
\renewcommand{\thetikzsubfigcounter}{\the\numexpr\value{figure}+1\relax\alph{tikzsubfigcounter}}

\newcounter{tikzsubfigcounterinvisible}[figure]
\renewcommand{\thetikzsubfigcounterinvisible}{\the\numexpr\value{figure}+1\relax\alph{tikzsubfigcounterinvisible}}

\newcommand{\refone}[1]{\textcolor{black}{#1}}
\newcommand{\reftwo}[1]{\textcolor{black}{#1}}

\numberwithin{equation}{section}

\title{Parameter identification in uncertain scalar conservation laws discretized with the discontinuous stochastic Galerkin Scheme
}

\author[ls]{Louisa Schlachter}
\address[ls]{Fachbereich Mathematik, TU Kaiserslautern, Erwin-Schr\"odinger-Str., 67663 Kaiserslautern, Germany, {\tt schlacht@mathematik.uni-kl.de}}

\author[ok]{Claudia Totzeck}
\address[ok]{Fachbereich Mathematik, TU Kaiserslautern, Erwin-Schr\"odinger-Str., 67663 Kaiserslautern, Germany, {\tt totzeck@mathematik.uni-kl.de}}

\date{}

\definecolor{greenyellow}   {cmyk}{0.15, 0   , 0.69, 0   }
\definecolor{yellow}        {cmyk}{0   , 0   , 1   , 0   }
\definecolor{goldenrod}     {cmyk}{0   , 0.10, 0.84, 0   }
\definecolor{dandelion}     {cmyk}{0   , 0.29, 0.84, 0   }
\definecolor{apricot}       {cmyk}{0   , 0.32, 0.52, 0   }
\definecolor{peach}         {cmyk}{0   , 0.50, 0.70, 0   }
\definecolor{melon}         {cmyk}{0   , 0.46, 0.50, 0   }
\definecolor{yelloworange}  {cmyk}{0   , 0.42, 1   , 0   }
\definecolor{orange}        {cmyk}{0   , 0.61, 0.87, 0   }
\definecolor{burntorange}   {cmyk}{0   , 0.51, 1   , 0   }
\definecolor{bittersweet}   {cmyk}{0   , 0.75, 1   , 0.24}
\definecolor{redorange}     {cmyk}{0   , 0.77, 0.87, 0   }
\definecolor{mahogany}      {cmyk}{0   , 0.85, 0.87, 0.35}
\definecolor{maroon}        {cmyk}{0   , 0.87, 0.68, 0.32}
\definecolor{brickred}      {cmyk}{0   , 0.89, 0.94, 0.28}
\definecolor{red}           {cmyk}{0   , 1   , 1   , 0   }
\definecolor{orangered}     {cmyk}{0   , 1   , 0.50, 0   }
\definecolor{rubinered}     {cmyk}{0   , 1   , 0.13, 0   }
\definecolor{wildstrawberry}{cmyk}{0   , 0.96, 0.39, 0   }
\definecolor{salmon}        {cmyk}{0   , 0.53, 0.38, 0   }
\definecolor{carnationpink} {cmyk}{0   , 0.63, 0   , 0   }
\definecolor{magenta}       {cmyk}{0   , 1   , 0   , 0   }
\definecolor{violetred}     {cmyk}{0   , 0.81, 0   , 0   }
\definecolor{rhodamine}     {cmyk}{0   , 0.82, 0   , 0   }
\definecolor{mulberry}      {cmyk}{0.34, 0.90, 0   , 0.02}
\definecolor{redviolet}     {cmyk}{0.07, 0.90, 0   , 0.34}
\definecolor{fuchsia}       {cmyk}{0.47, 0.91, 0   , 0.08}
\definecolor{lavender}      {cmyk}{0   , 0.48, 0   , 0   }
\definecolor{thistle}       {cmyk}{0.12, 0.59, 0   , 0   }
\definecolor{orchid}        {cmyk}{0.32, 0.64, 0   , 0   }
\definecolor{darkorchid}    {cmyk}{0.40, 0.80, 0.20, 0   }
\definecolor{purple}        {cmyk}{0.45, 0.86, 0   , 0   }
\definecolor{plum}          {cmyk}{0.50, 1   , 0   , 0   }
\definecolor{violet}        {cmyk}{0.79, 0.88, 0   , 0   }
\definecolor{royalpurple}   {cmyk}{0.75, 0.90, 0   , 0   }
\definecolor{blueviolet}    {cmyk}{0.86, 0.91, 0   , 0.04}
\definecolor{periwinkle}    {cmyk}{0.57, 0.55, 0   , 0   }
\definecolor{cadetblue}     {cmyk}{0.62, 0.57, 0.23, 0   }
\definecolor{cornflowerblue}{cmyk}{0.65, 0.13, 0   , 0   }
\definecolor{midnightblue}  {cmyk}{0.98, 0.13, 0   , 0.43}
\definecolor{navyblue}      {cmyk}{0.94, 0.54, 0   , 0   }
\definecolor{royalblue}     {cmyk}{1   , 0.50, 0   , 0   }
\definecolor{blue}          {cmyk}{1   , 1   , 0   , 0   }
\definecolor{cerulean}      {cmyk}{0.94, 0.11, 0   , 0   }
\definecolor{cyan}          {cmyk}{1   , 0   , 0   , 0   }
\definecolor{processblue}   {cmyk}{0.96, 0   , 0   , 0   }
\definecolor{skyblue}       {cmyk}{0.62, 0   , 0.12, 0   }
\definecolor{turquoise}     {cmyk}{0.85, 0   , 0.20, 0   }
\definecolor{tealblue}      {cmyk}{0.86, 0   , 0.34, 0.02}
\definecolor{aquamarine}    {cmyk}{0.82, 0   , 0.30, 0   }
\definecolor{bluegreen}     {cmyk}{0.85, 0   , 0.33, 0   }
\definecolor{emerald}       {cmyk}{1   , 0   , 0.50, 0   }
\definecolor{junglegreen}   {cmyk}{0.99, 0   , 0.52, 0   }
\definecolor{seagreen}      {cmyk}{0.69, 0   , 0.50, 0   }
\definecolor{green}         {cmyk}{1   , 0   , 1   , 0   }
\definecolor{forestgreen}   {cmyk}{0.91, 0   , 0.88, 0.12}
\definecolor{pinegreen}     {cmyk}{0.92, 0   , 0.59, 0.25}
\definecolor{limegreen}     {cmyk}{0.50, 0   , 1   , 0   }
\definecolor{yellowgreen}   {cmyk}{0.44, 0   , 0.74, 0   }
\definecolor{springgreen}   {cmyk}{0.26, 0   , 0.76, 0   }
\definecolor{olivegreen}    {cmyk}{0.64, 0   , 0.95, 0.40}
\definecolor{rawsienna}     {cmyk}{0   , 0.72, 1   , 0.45}
\definecolor{sepia}         {cmyk}{0   , 0.83, 1   , 0.70}
\definecolor{brown}         {cmyk}{0   , 0.81, 1   , 0.60}
\definecolor{tan}           {cmyk}{0.14, 0.42, 0.56, 0   }
\definecolor{gray}          {cmyk}{0   , 0   , 0   , 0.50}
\definecolor{black}         {cmyk}{0   , 0   , 0   , 1   }
\definecolor{white}         {cmyk}{0   , 0   , 0   , 0   }
\definecolor{tuklblue}{RGB}{0,95,140}
\definecolor{tuklred}{RGB}{185,40,25} 

\usepgfplotslibrary{groupplots}

\pgfplotsset{
	colormap={jet}{
rgb(0.000000 pt)=(0.000000,0.000000,0.504000);
rgb(1.000000 pt)=(0.000000,0.000000,0.508000);
rgb(2.000000 pt)=(0.000000,0.000000,0.512000);
rgb(3.000000 pt)=(0.000000,0.000000,0.516000);
rgb(4.000000 pt)=(0.000000,0.000000,0.520000);
rgb(5.000000 pt)=(0.000000,0.000000,0.524000);
rgb(6.000000 pt)=(0.000000,0.000000,0.528000);
rgb(7.000000 pt)=(0.000000,0.000000,0.532000);
rgb(8.000000 pt)=(0.000000,0.000000,0.536000);
rgb(9.000000 pt)=(0.000000,0.000000,0.540000);
rgb(10.000000 pt)=(0.000000,0.000000,0.544000);
rgb(11.000000 pt)=(0.000000,0.000000,0.548000);
rgb(12.000000 pt)=(0.000000,0.000000,0.552000);
rgb(13.000000 pt)=(0.000000,0.000000,0.556000);
rgb(14.000000 pt)=(0.000000,0.000000,0.560000);
rgb(15.000000 pt)=(0.000000,0.000000,0.564000);
rgb(16.000000 pt)=(0.000000,0.000000,0.568000);
rgb(17.000000 pt)=(0.000000,0.000000,0.572000);
rgb(18.000000 pt)=(0.000000,0.000000,0.576000);
rgb(19.000000 pt)=(0.000000,0.000000,0.580000);
rgb(20.000000 pt)=(0.000000,0.000000,0.584000);
rgb(21.000000 pt)=(0.000000,0.000000,0.588000);
rgb(22.000000 pt)=(0.000000,0.000000,0.592000);
rgb(23.000000 pt)=(0.000000,0.000000,0.596000);
rgb(24.000000 pt)=(0.000000,0.000000,0.600000);
rgb(25.000000 pt)=(0.000000,0.000000,0.604000);
rgb(26.000000 pt)=(0.000000,0.000000,0.608000);
rgb(27.000000 pt)=(0.000000,0.000000,0.612000);
rgb(28.000000 pt)=(0.000000,0.000000,0.616000);
rgb(29.000000 pt)=(0.000000,0.000000,0.620000);
rgb(30.000000 pt)=(0.000000,0.000000,0.624000);
rgb(31.000000 pt)=(0.000000,0.000000,0.628000);
rgb(32.000000 pt)=(0.000000,0.000000,0.632000);
rgb(33.000000 pt)=(0.000000,0.000000,0.636000);
rgb(34.000000 pt)=(0.000000,0.000000,0.640000);
rgb(35.000000 pt)=(0.000000,0.000000,0.644000);
rgb(36.000000 pt)=(0.000000,0.000000,0.648000);
rgb(37.000000 pt)=(0.000000,0.000000,0.652000);
rgb(38.000000 pt)=(0.000000,0.000000,0.656000);
rgb(39.000000 pt)=(0.000000,0.000000,0.660000);
rgb(40.000000 pt)=(0.000000,0.000000,0.664000);
rgb(41.000000 pt)=(0.000000,0.000000,0.668000);
rgb(42.000000 pt)=(0.000000,0.000000,0.672000);
rgb(43.000000 pt)=(0.000000,0.000000,0.676000);
rgb(44.000000 pt)=(0.000000,0.000000,0.680000);
rgb(45.000000 pt)=(0.000000,0.000000,0.684000);
rgb(46.000000 pt)=(0.000000,0.000000,0.688000);
rgb(47.000000 pt)=(0.000000,0.000000,0.692000);
rgb(48.000000 pt)=(0.000000,0.000000,0.696000);
rgb(49.000000 pt)=(0.000000,0.000000,0.700000);
rgb(50.000000 pt)=(0.000000,0.000000,0.704000);
rgb(51.000000 pt)=(0.000000,0.000000,0.708000);
rgb(52.000000 pt)=(0.000000,0.000000,0.712000);
rgb(53.000000 pt)=(0.000000,0.000000,0.716000);
rgb(54.000000 pt)=(0.000000,0.000000,0.720000);
rgb(55.000000 pt)=(0.000000,0.000000,0.724000);
rgb(56.000000 pt)=(0.000000,0.000000,0.728000);
rgb(57.000000 pt)=(0.000000,0.000000,0.732000);
rgb(58.000000 pt)=(0.000000,0.000000,0.736000);
rgb(59.000000 pt)=(0.000000,0.000000,0.740000);
rgb(60.000000 pt)=(0.000000,0.000000,0.744000);
rgb(61.000000 pt)=(0.000000,0.000000,0.748000);
rgb(62.000000 pt)=(0.000000,0.000000,0.752000);
rgb(63.000000 pt)=(0.000000,0.000000,0.756000);
rgb(64.000000 pt)=(0.000000,0.000000,0.760000);
rgb(65.000000 pt)=(0.000000,0.000000,0.764000);
rgb(66.000000 pt)=(0.000000,0.000000,0.768000);
rgb(67.000000 pt)=(0.000000,0.000000,0.772000);
rgb(68.000000 pt)=(0.000000,0.000000,0.776000);
rgb(69.000000 pt)=(0.000000,0.000000,0.780000);
rgb(70.000000 pt)=(0.000000,0.000000,0.784000);
rgb(71.000000 pt)=(0.000000,0.000000,0.788000);
rgb(72.000000 pt)=(0.000000,0.000000,0.792000);
rgb(73.000000 pt)=(0.000000,0.000000,0.796000);
rgb(74.000000 pt)=(0.000000,0.000000,0.800000);
rgb(75.000000 pt)=(0.000000,0.000000,0.804000);
rgb(76.000000 pt)=(0.000000,0.000000,0.808000);
rgb(77.000000 pt)=(0.000000,0.000000,0.812000);
rgb(78.000000 pt)=(0.000000,0.000000,0.816000);
rgb(79.000000 pt)=(0.000000,0.000000,0.820000);
rgb(80.000000 pt)=(0.000000,0.000000,0.824000);
rgb(81.000000 pt)=(0.000000,0.000000,0.828000);
rgb(82.000000 pt)=(0.000000,0.000000,0.832000);
rgb(83.000000 pt)=(0.000000,0.000000,0.836000);
rgb(84.000000 pt)=(0.000000,0.000000,0.840000);
rgb(85.000000 pt)=(0.000000,0.000000,0.844000);
rgb(86.000000 pt)=(0.000000,0.000000,0.848000);
rgb(87.000000 pt)=(0.000000,0.000000,0.852000);
rgb(88.000000 pt)=(0.000000,0.000000,0.856000);
rgb(89.000000 pt)=(0.000000,0.000000,0.860000);
rgb(90.000000 pt)=(0.000000,0.000000,0.864000);
rgb(91.000000 pt)=(0.000000,0.000000,0.868000);
rgb(92.000000 pt)=(0.000000,0.000000,0.872000);
rgb(93.000000 pt)=(0.000000,0.000000,0.876000);
rgb(94.000000 pt)=(0.000000,0.000000,0.880000);
rgb(95.000000 pt)=(0.000000,0.000000,0.884000);
rgb(96.000000 pt)=(0.000000,0.000000,0.888000);
rgb(97.000000 pt)=(0.000000,0.000000,0.892000);
rgb(98.000000 pt)=(0.000000,0.000000,0.896000);
rgb(99.000000 pt)=(0.000000,0.000000,0.900000);
rgb(100.000000 pt)=(0.000000,0.000000,0.904000);
rgb(101.000000 pt)=(0.000000,0.000000,0.908000);
rgb(102.000000 pt)=(0.000000,0.000000,0.912000);
rgb(103.000000 pt)=(0.000000,0.000000,0.916000);
rgb(104.000000 pt)=(0.000000,0.000000,0.920000);
rgb(105.000000 pt)=(0.000000,0.000000,0.924000);
rgb(106.000000 pt)=(0.000000,0.000000,0.928000);
rgb(107.000000 pt)=(0.000000,0.000000,0.932000);
rgb(108.000000 pt)=(0.000000,0.000000,0.936000);
rgb(109.000000 pt)=(0.000000,0.000000,0.940000);
rgb(110.000000 pt)=(0.000000,0.000000,0.944000);
rgb(111.000000 pt)=(0.000000,0.000000,0.948000);
rgb(112.000000 pt)=(0.000000,0.000000,0.952000);
rgb(113.000000 pt)=(0.000000,0.000000,0.956000);
rgb(114.000000 pt)=(0.000000,0.000000,0.960000);
rgb(115.000000 pt)=(0.000000,0.000000,0.964000);
rgb(116.000000 pt)=(0.000000,0.000000,0.968000);
rgb(117.000000 pt)=(0.000000,0.000000,0.972000);
rgb(118.000000 pt)=(0.000000,0.000000,0.976000);
rgb(119.000000 pt)=(0.000000,0.000000,0.980000);
rgb(120.000000 pt)=(0.000000,0.000000,0.984000);
rgb(121.000000 pt)=(0.000000,0.000000,0.988000);
rgb(122.000000 pt)=(0.000000,0.000000,0.992000);
rgb(123.000000 pt)=(0.000000,0.000000,0.996000);
rgb(124.000000 pt)=(0.000000,0.000000,1.000000);
rgb(125.000000 pt)=(0.000000,0.004000,1.000000);
rgb(126.000000 pt)=(0.000000,0.008000,1.000000);
rgb(127.000000 pt)=(0.000000,0.012000,1.000000);
rgb(128.000000 pt)=(0.000000,0.016000,1.000000);
rgb(129.000000 pt)=(0.000000,0.020000,1.000000);
rgb(130.000000 pt)=(0.000000,0.024000,1.000000);
rgb(131.000000 pt)=(0.000000,0.028000,1.000000);
rgb(132.000000 pt)=(0.000000,0.032000,1.000000);
rgb(133.000000 pt)=(0.000000,0.036000,1.000000);
rgb(134.000000 pt)=(0.000000,0.040000,1.000000);
rgb(135.000000 pt)=(0.000000,0.044000,1.000000);
rgb(136.000000 pt)=(0.000000,0.048000,1.000000);
rgb(137.000000 pt)=(0.000000,0.052000,1.000000);
rgb(138.000000 pt)=(0.000000,0.056000,1.000000);
rgb(139.000000 pt)=(0.000000,0.060000,1.000000);
rgb(140.000000 pt)=(0.000000,0.064000,1.000000);
rgb(141.000000 pt)=(0.000000,0.068000,1.000000);
rgb(142.000000 pt)=(0.000000,0.072000,1.000000);
rgb(143.000000 pt)=(0.000000,0.076000,1.000000);
rgb(144.000000 pt)=(0.000000,0.080000,1.000000);
rgb(145.000000 pt)=(0.000000,0.084000,1.000000);
rgb(146.000000 pt)=(0.000000,0.088000,1.000000);
rgb(147.000000 pt)=(0.000000,0.092000,1.000000);
rgb(148.000000 pt)=(0.000000,0.096000,1.000000);
rgb(149.000000 pt)=(0.000000,0.100000,1.000000);
rgb(150.000000 pt)=(0.000000,0.104000,1.000000);
rgb(151.000000 pt)=(0.000000,0.108000,1.000000);
rgb(152.000000 pt)=(0.000000,0.112000,1.000000);
rgb(153.000000 pt)=(0.000000,0.116000,1.000000);
rgb(154.000000 pt)=(0.000000,0.120000,1.000000);
rgb(155.000000 pt)=(0.000000,0.124000,1.000000);
rgb(156.000000 pt)=(0.000000,0.128000,1.000000);
rgb(157.000000 pt)=(0.000000,0.132000,1.000000);
rgb(158.000000 pt)=(0.000000,0.136000,1.000000);
rgb(159.000000 pt)=(0.000000,0.140000,1.000000);
rgb(160.000000 pt)=(0.000000,0.144000,1.000000);
rgb(161.000000 pt)=(0.000000,0.148000,1.000000);
rgb(162.000000 pt)=(0.000000,0.152000,1.000000);
rgb(163.000000 pt)=(0.000000,0.156000,1.000000);
rgb(164.000000 pt)=(0.000000,0.160000,1.000000);
rgb(165.000000 pt)=(0.000000,0.164000,1.000000);
rgb(166.000000 pt)=(0.000000,0.168000,1.000000);
rgb(167.000000 pt)=(0.000000,0.172000,1.000000);
rgb(168.000000 pt)=(0.000000,0.176000,1.000000);
rgb(169.000000 pt)=(0.000000,0.180000,1.000000);
rgb(170.000000 pt)=(0.000000,0.184000,1.000000);
rgb(171.000000 pt)=(0.000000,0.188000,1.000000);
rgb(172.000000 pt)=(0.000000,0.192000,1.000000);
rgb(173.000000 pt)=(0.000000,0.196000,1.000000);
rgb(174.000000 pt)=(0.000000,0.200000,1.000000);
rgb(175.000000 pt)=(0.000000,0.204000,1.000000);
rgb(176.000000 pt)=(0.000000,0.208000,1.000000);
rgb(177.000000 pt)=(0.000000,0.212000,1.000000);
rgb(178.000000 pt)=(0.000000,0.216000,1.000000);
rgb(179.000000 pt)=(0.000000,0.220000,1.000000);
rgb(180.000000 pt)=(0.000000,0.224000,1.000000);
rgb(181.000000 pt)=(0.000000,0.228000,1.000000);
rgb(182.000000 pt)=(0.000000,0.232000,1.000000);
rgb(183.000000 pt)=(0.000000,0.236000,1.000000);
rgb(184.000000 pt)=(0.000000,0.240000,1.000000);
rgb(185.000000 pt)=(0.000000,0.244000,1.000000);
rgb(186.000000 pt)=(0.000000,0.248000,1.000000);
rgb(187.000000 pt)=(0.000000,0.252000,1.000000);
rgb(188.000000 pt)=(0.000000,0.256000,1.000000);
rgb(189.000000 pt)=(0.000000,0.260000,1.000000);
rgb(190.000000 pt)=(0.000000,0.264000,1.000000);
rgb(191.000000 pt)=(0.000000,0.268000,1.000000);
rgb(192.000000 pt)=(0.000000,0.272000,1.000000);
rgb(193.000000 pt)=(0.000000,0.276000,1.000000);
rgb(194.000000 pt)=(0.000000,0.280000,1.000000);
rgb(195.000000 pt)=(0.000000,0.284000,1.000000);
rgb(196.000000 pt)=(0.000000,0.288000,1.000000);
rgb(197.000000 pt)=(0.000000,0.292000,1.000000);
rgb(198.000000 pt)=(0.000000,0.296000,1.000000);
rgb(199.000000 pt)=(0.000000,0.300000,1.000000);
rgb(200.000000 pt)=(0.000000,0.304000,1.000000);
rgb(201.000000 pt)=(0.000000,0.308000,1.000000);
rgb(202.000000 pt)=(0.000000,0.312000,1.000000);
rgb(203.000000 pt)=(0.000000,0.316000,1.000000);
rgb(204.000000 pt)=(0.000000,0.320000,1.000000);
rgb(205.000000 pt)=(0.000000,0.324000,1.000000);
rgb(206.000000 pt)=(0.000000,0.328000,1.000000);
rgb(207.000000 pt)=(0.000000,0.332000,1.000000);
rgb(208.000000 pt)=(0.000000,0.336000,1.000000);
rgb(209.000000 pt)=(0.000000,0.340000,1.000000);
rgb(210.000000 pt)=(0.000000,0.344000,1.000000);
rgb(211.000000 pt)=(0.000000,0.348000,1.000000);
rgb(212.000000 pt)=(0.000000,0.352000,1.000000);
rgb(213.000000 pt)=(0.000000,0.356000,1.000000);
rgb(214.000000 pt)=(0.000000,0.360000,1.000000);
rgb(215.000000 pt)=(0.000000,0.364000,1.000000);
rgb(216.000000 pt)=(0.000000,0.368000,1.000000);
rgb(217.000000 pt)=(0.000000,0.372000,1.000000);
rgb(218.000000 pt)=(0.000000,0.376000,1.000000);
rgb(219.000000 pt)=(0.000000,0.380000,1.000000);
rgb(220.000000 pt)=(0.000000,0.384000,1.000000);
rgb(221.000000 pt)=(0.000000,0.388000,1.000000);
rgb(222.000000 pt)=(0.000000,0.392000,1.000000);
rgb(223.000000 pt)=(0.000000,0.396000,1.000000);
rgb(224.000000 pt)=(0.000000,0.400000,1.000000);
rgb(225.000000 pt)=(0.000000,0.404000,1.000000);
rgb(226.000000 pt)=(0.000000,0.408000,1.000000);
rgb(227.000000 pt)=(0.000000,0.412000,1.000000);
rgb(228.000000 pt)=(0.000000,0.416000,1.000000);
rgb(229.000000 pt)=(0.000000,0.420000,1.000000);
rgb(230.000000 pt)=(0.000000,0.424000,1.000000);
rgb(231.000000 pt)=(0.000000,0.428000,1.000000);
rgb(232.000000 pt)=(0.000000,0.432000,1.000000);
rgb(233.000000 pt)=(0.000000,0.436000,1.000000);
rgb(234.000000 pt)=(0.000000,0.440000,1.000000);
rgb(235.000000 pt)=(0.000000,0.444000,1.000000);
rgb(236.000000 pt)=(0.000000,0.448000,1.000000);
rgb(237.000000 pt)=(0.000000,0.452000,1.000000);
rgb(238.000000 pt)=(0.000000,0.456000,1.000000);
rgb(239.000000 pt)=(0.000000,0.460000,1.000000);
rgb(240.000000 pt)=(0.000000,0.464000,1.000000);
rgb(241.000000 pt)=(0.000000,0.468000,1.000000);
rgb(242.000000 pt)=(0.000000,0.472000,1.000000);
rgb(243.000000 pt)=(0.000000,0.476000,1.000000);
rgb(244.000000 pt)=(0.000000,0.480000,1.000000);
rgb(245.000000 pt)=(0.000000,0.484000,1.000000);
rgb(246.000000 pt)=(0.000000,0.488000,1.000000);
rgb(247.000000 pt)=(0.000000,0.492000,1.000000);
rgb(248.000000 pt)=(0.000000,0.496000,1.000000);
rgb(249.000000 pt)=(0.000000,0.500000,1.000000);
rgb(250.000000 pt)=(0.000000,0.504000,1.000000);
rgb(251.000000 pt)=(0.000000,0.508000,1.000000);
rgb(252.000000 pt)=(0.000000,0.512000,1.000000);
rgb(253.000000 pt)=(0.000000,0.516000,1.000000);
rgb(254.000000 pt)=(0.000000,0.520000,1.000000);
rgb(255.000000 pt)=(0.000000,0.524000,1.000000);
rgb(256.000000 pt)=(0.000000,0.528000,1.000000);
rgb(257.000000 pt)=(0.000000,0.532000,1.000000);
rgb(258.000000 pt)=(0.000000,0.536000,1.000000);
rgb(259.000000 pt)=(0.000000,0.540000,1.000000);
rgb(260.000000 pt)=(0.000000,0.544000,1.000000);
rgb(261.000000 pt)=(0.000000,0.548000,1.000000);
rgb(262.000000 pt)=(0.000000,0.552000,1.000000);
rgb(263.000000 pt)=(0.000000,0.556000,1.000000);
rgb(264.000000 pt)=(0.000000,0.560000,1.000000);
rgb(265.000000 pt)=(0.000000,0.564000,1.000000);
rgb(266.000000 pt)=(0.000000,0.568000,1.000000);
rgb(267.000000 pt)=(0.000000,0.572000,1.000000);
rgb(268.000000 pt)=(0.000000,0.576000,1.000000);
rgb(269.000000 pt)=(0.000000,0.580000,1.000000);
rgb(270.000000 pt)=(0.000000,0.584000,1.000000);
rgb(271.000000 pt)=(0.000000,0.588000,1.000000);
rgb(272.000000 pt)=(0.000000,0.592000,1.000000);
rgb(273.000000 pt)=(0.000000,0.596000,1.000000);
rgb(274.000000 pt)=(0.000000,0.600000,1.000000);
rgb(275.000000 pt)=(0.000000,0.604000,1.000000);
rgb(276.000000 pt)=(0.000000,0.608000,1.000000);
rgb(277.000000 pt)=(0.000000,0.612000,1.000000);
rgb(278.000000 pt)=(0.000000,0.616000,1.000000);
rgb(279.000000 pt)=(0.000000,0.620000,1.000000);
rgb(280.000000 pt)=(0.000000,0.624000,1.000000);
rgb(281.000000 pt)=(0.000000,0.628000,1.000000);
rgb(282.000000 pt)=(0.000000,0.632000,1.000000);
rgb(283.000000 pt)=(0.000000,0.636000,1.000000);
rgb(284.000000 pt)=(0.000000,0.640000,1.000000);
rgb(285.000000 pt)=(0.000000,0.644000,1.000000);
rgb(286.000000 pt)=(0.000000,0.648000,1.000000);
rgb(287.000000 pt)=(0.000000,0.652000,1.000000);
rgb(288.000000 pt)=(0.000000,0.656000,1.000000);
rgb(289.000000 pt)=(0.000000,0.660000,1.000000);
rgb(290.000000 pt)=(0.000000,0.664000,1.000000);
rgb(291.000000 pt)=(0.000000,0.668000,1.000000);
rgb(292.000000 pt)=(0.000000,0.672000,1.000000);
rgb(293.000000 pt)=(0.000000,0.676000,1.000000);
rgb(294.000000 pt)=(0.000000,0.680000,1.000000);
rgb(295.000000 pt)=(0.000000,0.684000,1.000000);
rgb(296.000000 pt)=(0.000000,0.688000,1.000000);
rgb(297.000000 pt)=(0.000000,0.692000,1.000000);
rgb(298.000000 pt)=(0.000000,0.696000,1.000000);
rgb(299.000000 pt)=(0.000000,0.700000,1.000000);
rgb(300.000000 pt)=(0.000000,0.704000,1.000000);
rgb(301.000000 pt)=(0.000000,0.708000,1.000000);
rgb(302.000000 pt)=(0.000000,0.712000,1.000000);
rgb(303.000000 pt)=(0.000000,0.716000,1.000000);
rgb(304.000000 pt)=(0.000000,0.720000,1.000000);
rgb(305.000000 pt)=(0.000000,0.724000,1.000000);
rgb(306.000000 pt)=(0.000000,0.728000,1.000000);
rgb(307.000000 pt)=(0.000000,0.732000,1.000000);
rgb(308.000000 pt)=(0.000000,0.736000,1.000000);
rgb(309.000000 pt)=(0.000000,0.740000,1.000000);
rgb(310.000000 pt)=(0.000000,0.744000,1.000000);
rgb(311.000000 pt)=(0.000000,0.748000,1.000000);
rgb(312.000000 pt)=(0.000000,0.752000,1.000000);
rgb(313.000000 pt)=(0.000000,0.756000,1.000000);
rgb(314.000000 pt)=(0.000000,0.760000,1.000000);
rgb(315.000000 pt)=(0.000000,0.764000,1.000000);
rgb(316.000000 pt)=(0.000000,0.768000,1.000000);
rgb(317.000000 pt)=(0.000000,0.772000,1.000000);
rgb(318.000000 pt)=(0.000000,0.776000,1.000000);
rgb(319.000000 pt)=(0.000000,0.780000,1.000000);
rgb(320.000000 pt)=(0.000000,0.784000,1.000000);
rgb(321.000000 pt)=(0.000000,0.788000,1.000000);
rgb(322.000000 pt)=(0.000000,0.792000,1.000000);
rgb(323.000000 pt)=(0.000000,0.796000,1.000000);
rgb(324.000000 pt)=(0.000000,0.800000,1.000000);
rgb(325.000000 pt)=(0.000000,0.804000,1.000000);
rgb(326.000000 pt)=(0.000000,0.808000,1.000000);
rgb(327.000000 pt)=(0.000000,0.812000,1.000000);
rgb(328.000000 pt)=(0.000000,0.816000,1.000000);
rgb(329.000000 pt)=(0.000000,0.820000,1.000000);
rgb(330.000000 pt)=(0.000000,0.824000,1.000000);
rgb(331.000000 pt)=(0.000000,0.828000,1.000000);
rgb(332.000000 pt)=(0.000000,0.832000,1.000000);
rgb(333.000000 pt)=(0.000000,0.836000,1.000000);
rgb(334.000000 pt)=(0.000000,0.840000,1.000000);
rgb(335.000000 pt)=(0.000000,0.844000,1.000000);
rgb(336.000000 pt)=(0.000000,0.848000,1.000000);
rgb(337.000000 pt)=(0.000000,0.852000,1.000000);
rgb(338.000000 pt)=(0.000000,0.856000,1.000000);
rgb(339.000000 pt)=(0.000000,0.860000,1.000000);
rgb(340.000000 pt)=(0.000000,0.864000,1.000000);
rgb(341.000000 pt)=(0.000000,0.868000,1.000000);
rgb(342.000000 pt)=(0.000000,0.872000,1.000000);
rgb(343.000000 pt)=(0.000000,0.876000,1.000000);
rgb(344.000000 pt)=(0.000000,0.880000,1.000000);
rgb(345.000000 pt)=(0.000000,0.884000,1.000000);
rgb(346.000000 pt)=(0.000000,0.888000,1.000000);
rgb(347.000000 pt)=(0.000000,0.892000,1.000000);
rgb(348.000000 pt)=(0.000000,0.896000,1.000000);
rgb(349.000000 pt)=(0.000000,0.900000,1.000000);
rgb(350.000000 pt)=(0.000000,0.904000,1.000000);
rgb(351.000000 pt)=(0.000000,0.908000,1.000000);
rgb(352.000000 pt)=(0.000000,0.912000,1.000000);
rgb(353.000000 pt)=(0.000000,0.916000,1.000000);
rgb(354.000000 pt)=(0.000000,0.920000,1.000000);
rgb(355.000000 pt)=(0.000000,0.924000,1.000000);
rgb(356.000000 pt)=(0.000000,0.928000,1.000000);
rgb(357.000000 pt)=(0.000000,0.932000,1.000000);
rgb(358.000000 pt)=(0.000000,0.936000,1.000000);
rgb(359.000000 pt)=(0.000000,0.940000,1.000000);
rgb(360.000000 pt)=(0.000000,0.944000,1.000000);
rgb(361.000000 pt)=(0.000000,0.948000,1.000000);
rgb(362.000000 pt)=(0.000000,0.952000,1.000000);
rgb(363.000000 pt)=(0.000000,0.956000,1.000000);
rgb(364.000000 pt)=(0.000000,0.960000,1.000000);
rgb(365.000000 pt)=(0.000000,0.964000,1.000000);
rgb(366.000000 pt)=(0.000000,0.968000,1.000000);
rgb(367.000000 pt)=(0.000000,0.972000,1.000000);
rgb(368.000000 pt)=(0.000000,0.976000,1.000000);
rgb(369.000000 pt)=(0.000000,0.980000,1.000000);
rgb(370.000000 pt)=(0.000000,0.984000,1.000000);
rgb(371.000000 pt)=(0.000000,0.988000,1.000000);
rgb(372.000000 pt)=(0.000000,0.992000,1.000000);
rgb(373.000000 pt)=(0.000000,0.996000,1.000000);
rgb(374.000000 pt)=(0.000000,1.000000,1.000000);
rgb(375.000000 pt)=(0.004000,1.000000,0.996000);
rgb(376.000000 pt)=(0.008000,1.000000,0.992000);
rgb(377.000000 pt)=(0.012000,1.000000,0.988000);
rgb(378.000000 pt)=(0.016000,1.000000,0.984000);
rgb(379.000000 pt)=(0.020000,1.000000,0.980000);
rgb(380.000000 pt)=(0.024000,1.000000,0.976000);
rgb(381.000000 pt)=(0.028000,1.000000,0.972000);
rgb(382.000000 pt)=(0.032000,1.000000,0.968000);
rgb(383.000000 pt)=(0.036000,1.000000,0.964000);
rgb(384.000000 pt)=(0.040000,1.000000,0.960000);
rgb(385.000000 pt)=(0.044000,1.000000,0.956000);
rgb(386.000000 pt)=(0.048000,1.000000,0.952000);
rgb(387.000000 pt)=(0.052000,1.000000,0.948000);
rgb(388.000000 pt)=(0.056000,1.000000,0.944000);
rgb(389.000000 pt)=(0.060000,1.000000,0.940000);
rgb(390.000000 pt)=(0.064000,1.000000,0.936000);
rgb(391.000000 pt)=(0.068000,1.000000,0.932000);
rgb(392.000000 pt)=(0.072000,1.000000,0.928000);
rgb(393.000000 pt)=(0.076000,1.000000,0.924000);
rgb(394.000000 pt)=(0.080000,1.000000,0.920000);
rgb(395.000000 pt)=(0.084000,1.000000,0.916000);
rgb(396.000000 pt)=(0.088000,1.000000,0.912000);
rgb(397.000000 pt)=(0.092000,1.000000,0.908000);
rgb(398.000000 pt)=(0.096000,1.000000,0.904000);
rgb(399.000000 pt)=(0.100000,1.000000,0.900000);
rgb(400.000000 pt)=(0.104000,1.000000,0.896000);
rgb(401.000000 pt)=(0.108000,1.000000,0.892000);
rgb(402.000000 pt)=(0.112000,1.000000,0.888000);
rgb(403.000000 pt)=(0.116000,1.000000,0.884000);
rgb(404.000000 pt)=(0.120000,1.000000,0.880000);
rgb(405.000000 pt)=(0.124000,1.000000,0.876000);
rgb(406.000000 pt)=(0.128000,1.000000,0.872000);
rgb(407.000000 pt)=(0.132000,1.000000,0.868000);
rgb(408.000000 pt)=(0.136000,1.000000,0.864000);
rgb(409.000000 pt)=(0.140000,1.000000,0.860000);
rgb(410.000000 pt)=(0.144000,1.000000,0.856000);
rgb(411.000000 pt)=(0.148000,1.000000,0.852000);
rgb(412.000000 pt)=(0.152000,1.000000,0.848000);
rgb(413.000000 pt)=(0.156000,1.000000,0.844000);
rgb(414.000000 pt)=(0.160000,1.000000,0.840000);
rgb(415.000000 pt)=(0.164000,1.000000,0.836000);
rgb(416.000000 pt)=(0.168000,1.000000,0.832000);
rgb(417.000000 pt)=(0.172000,1.000000,0.828000);
rgb(418.000000 pt)=(0.176000,1.000000,0.824000);
rgb(419.000000 pt)=(0.180000,1.000000,0.820000);
rgb(420.000000 pt)=(0.184000,1.000000,0.816000);
rgb(421.000000 pt)=(0.188000,1.000000,0.812000);
rgb(422.000000 pt)=(0.192000,1.000000,0.808000);
rgb(423.000000 pt)=(0.196000,1.000000,0.804000);
rgb(424.000000 pt)=(0.200000,1.000000,0.800000);
rgb(425.000000 pt)=(0.204000,1.000000,0.796000);
rgb(426.000000 pt)=(0.208000,1.000000,0.792000);
rgb(427.000000 pt)=(0.212000,1.000000,0.788000);
rgb(428.000000 pt)=(0.216000,1.000000,0.784000);
rgb(429.000000 pt)=(0.220000,1.000000,0.780000);
rgb(430.000000 pt)=(0.224000,1.000000,0.776000);
rgb(431.000000 pt)=(0.228000,1.000000,0.772000);
rgb(432.000000 pt)=(0.232000,1.000000,0.768000);
rgb(433.000000 pt)=(0.236000,1.000000,0.764000);
rgb(434.000000 pt)=(0.240000,1.000000,0.760000);
rgb(435.000000 pt)=(0.244000,1.000000,0.756000);
rgb(436.000000 pt)=(0.248000,1.000000,0.752000);
rgb(437.000000 pt)=(0.252000,1.000000,0.748000);
rgb(438.000000 pt)=(0.256000,1.000000,0.744000);
rgb(439.000000 pt)=(0.260000,1.000000,0.740000);
rgb(440.000000 pt)=(0.264000,1.000000,0.736000);
rgb(441.000000 pt)=(0.268000,1.000000,0.732000);
rgb(442.000000 pt)=(0.272000,1.000000,0.728000);
rgb(443.000000 pt)=(0.276000,1.000000,0.724000);
rgb(444.000000 pt)=(0.280000,1.000000,0.720000);
rgb(445.000000 pt)=(0.284000,1.000000,0.716000);
rgb(446.000000 pt)=(0.288000,1.000000,0.712000);
rgb(447.000000 pt)=(0.292000,1.000000,0.708000);
rgb(448.000000 pt)=(0.296000,1.000000,0.704000);
rgb(449.000000 pt)=(0.300000,1.000000,0.700000);
rgb(450.000000 pt)=(0.304000,1.000000,0.696000);
rgb(451.000000 pt)=(0.308000,1.000000,0.692000);
rgb(452.000000 pt)=(0.312000,1.000000,0.688000);
rgb(453.000000 pt)=(0.316000,1.000000,0.684000);
rgb(454.000000 pt)=(0.320000,1.000000,0.680000);
rgb(455.000000 pt)=(0.324000,1.000000,0.676000);
rgb(456.000000 pt)=(0.328000,1.000000,0.672000);
rgb(457.000000 pt)=(0.332000,1.000000,0.668000);
rgb(458.000000 pt)=(0.336000,1.000000,0.664000);
rgb(459.000000 pt)=(0.340000,1.000000,0.660000);
rgb(460.000000 pt)=(0.344000,1.000000,0.656000);
rgb(461.000000 pt)=(0.348000,1.000000,0.652000);
rgb(462.000000 pt)=(0.352000,1.000000,0.648000);
rgb(463.000000 pt)=(0.356000,1.000000,0.644000);
rgb(464.000000 pt)=(0.360000,1.000000,0.640000);
rgb(465.000000 pt)=(0.364000,1.000000,0.636000);
rgb(466.000000 pt)=(0.368000,1.000000,0.632000);
rgb(467.000000 pt)=(0.372000,1.000000,0.628000);
rgb(468.000000 pt)=(0.376000,1.000000,0.624000);
rgb(469.000000 pt)=(0.380000,1.000000,0.620000);
rgb(470.000000 pt)=(0.384000,1.000000,0.616000);
rgb(471.000000 pt)=(0.388000,1.000000,0.612000);
rgb(472.000000 pt)=(0.392000,1.000000,0.608000);
rgb(473.000000 pt)=(0.396000,1.000000,0.604000);
rgb(474.000000 pt)=(0.400000,1.000000,0.600000);
rgb(475.000000 pt)=(0.404000,1.000000,0.596000);
rgb(476.000000 pt)=(0.408000,1.000000,0.592000);
rgb(477.000000 pt)=(0.412000,1.000000,0.588000);
rgb(478.000000 pt)=(0.416000,1.000000,0.584000);
rgb(479.000000 pt)=(0.420000,1.000000,0.580000);
rgb(480.000000 pt)=(0.424000,1.000000,0.576000);
rgb(481.000000 pt)=(0.428000,1.000000,0.572000);
rgb(482.000000 pt)=(0.432000,1.000000,0.568000);
rgb(483.000000 pt)=(0.436000,1.000000,0.564000);
rgb(484.000000 pt)=(0.440000,1.000000,0.560000);
rgb(485.000000 pt)=(0.444000,1.000000,0.556000);
rgb(486.000000 pt)=(0.448000,1.000000,0.552000);
rgb(487.000000 pt)=(0.452000,1.000000,0.548000);
rgb(488.000000 pt)=(0.456000,1.000000,0.544000);
rgb(489.000000 pt)=(0.460000,1.000000,0.540000);
rgb(490.000000 pt)=(0.464000,1.000000,0.536000);
rgb(491.000000 pt)=(0.468000,1.000000,0.532000);
rgb(492.000000 pt)=(0.472000,1.000000,0.528000);
rgb(493.000000 pt)=(0.476000,1.000000,0.524000);
rgb(494.000000 pt)=(0.480000,1.000000,0.520000);
rgb(495.000000 pt)=(0.484000,1.000000,0.516000);
rgb(496.000000 pt)=(0.488000,1.000000,0.512000);
rgb(497.000000 pt)=(0.492000,1.000000,0.508000);
rgb(498.000000 pt)=(0.496000,1.000000,0.504000);
rgb(499.000000 pt)=(0.500000,1.000000,0.500000);
rgb(500.000000 pt)=(0.504000,1.000000,0.496000);
rgb(501.000000 pt)=(0.508000,1.000000,0.492000);
rgb(502.000000 pt)=(0.512000,1.000000,0.488000);
rgb(503.000000 pt)=(0.516000,1.000000,0.484000);
rgb(504.000000 pt)=(0.520000,1.000000,0.480000);
rgb(505.000000 pt)=(0.524000,1.000000,0.476000);
rgb(506.000000 pt)=(0.528000,1.000000,0.472000);
rgb(507.000000 pt)=(0.532000,1.000000,0.468000);
rgb(508.000000 pt)=(0.536000,1.000000,0.464000);
rgb(509.000000 pt)=(0.540000,1.000000,0.460000);
rgb(510.000000 pt)=(0.544000,1.000000,0.456000);
rgb(511.000000 pt)=(0.548000,1.000000,0.452000);
rgb(512.000000 pt)=(0.552000,1.000000,0.448000);
rgb(513.000000 pt)=(0.556000,1.000000,0.444000);
rgb(514.000000 pt)=(0.560000,1.000000,0.440000);
rgb(515.000000 pt)=(0.564000,1.000000,0.436000);
rgb(516.000000 pt)=(0.568000,1.000000,0.432000);
rgb(517.000000 pt)=(0.572000,1.000000,0.428000);
rgb(518.000000 pt)=(0.576000,1.000000,0.424000);
rgb(519.000000 pt)=(0.580000,1.000000,0.420000);
rgb(520.000000 pt)=(0.584000,1.000000,0.416000);
rgb(521.000000 pt)=(0.588000,1.000000,0.412000);
rgb(522.000000 pt)=(0.592000,1.000000,0.408000);
rgb(523.000000 pt)=(0.596000,1.000000,0.404000);
rgb(524.000000 pt)=(0.600000,1.000000,0.400000);
rgb(525.000000 pt)=(0.604000,1.000000,0.396000);
rgb(526.000000 pt)=(0.608000,1.000000,0.392000);
rgb(527.000000 pt)=(0.612000,1.000000,0.388000);
rgb(528.000000 pt)=(0.616000,1.000000,0.384000);
rgb(529.000000 pt)=(0.620000,1.000000,0.380000);
rgb(530.000000 pt)=(0.624000,1.000000,0.376000);
rgb(531.000000 pt)=(0.628000,1.000000,0.372000);
rgb(532.000000 pt)=(0.632000,1.000000,0.368000);
rgb(533.000000 pt)=(0.636000,1.000000,0.364000);
rgb(534.000000 pt)=(0.640000,1.000000,0.360000);
rgb(535.000000 pt)=(0.644000,1.000000,0.356000);
rgb(536.000000 pt)=(0.648000,1.000000,0.352000);
rgb(537.000000 pt)=(0.652000,1.000000,0.348000);
rgb(538.000000 pt)=(0.656000,1.000000,0.344000);
rgb(539.000000 pt)=(0.660000,1.000000,0.340000);
rgb(540.000000 pt)=(0.664000,1.000000,0.336000);
rgb(541.000000 pt)=(0.668000,1.000000,0.332000);
rgb(542.000000 pt)=(0.672000,1.000000,0.328000);
rgb(543.000000 pt)=(0.676000,1.000000,0.324000);
rgb(544.000000 pt)=(0.680000,1.000000,0.320000);
rgb(545.000000 pt)=(0.684000,1.000000,0.316000);
rgb(546.000000 pt)=(0.688000,1.000000,0.312000);
rgb(547.000000 pt)=(0.692000,1.000000,0.308000);
rgb(548.000000 pt)=(0.696000,1.000000,0.304000);
rgb(549.000000 pt)=(0.700000,1.000000,0.300000);
rgb(550.000000 pt)=(0.704000,1.000000,0.296000);
rgb(551.000000 pt)=(0.708000,1.000000,0.292000);
rgb(552.000000 pt)=(0.712000,1.000000,0.288000);
rgb(553.000000 pt)=(0.716000,1.000000,0.284000);
rgb(554.000000 pt)=(0.720000,1.000000,0.280000);
rgb(555.000000 pt)=(0.724000,1.000000,0.276000);
rgb(556.000000 pt)=(0.728000,1.000000,0.272000);
rgb(557.000000 pt)=(0.732000,1.000000,0.268000);
rgb(558.000000 pt)=(0.736000,1.000000,0.264000);
rgb(559.000000 pt)=(0.740000,1.000000,0.260000);
rgb(560.000000 pt)=(0.744000,1.000000,0.256000);
rgb(561.000000 pt)=(0.748000,1.000000,0.252000);
rgb(562.000000 pt)=(0.752000,1.000000,0.248000);
rgb(563.000000 pt)=(0.756000,1.000000,0.244000);
rgb(564.000000 pt)=(0.760000,1.000000,0.240000);
rgb(565.000000 pt)=(0.764000,1.000000,0.236000);
rgb(566.000000 pt)=(0.768000,1.000000,0.232000);
rgb(567.000000 pt)=(0.772000,1.000000,0.228000);
rgb(568.000000 pt)=(0.776000,1.000000,0.224000);
rgb(569.000000 pt)=(0.780000,1.000000,0.220000);
rgb(570.000000 pt)=(0.784000,1.000000,0.216000);
rgb(571.000000 pt)=(0.788000,1.000000,0.212000);
rgb(572.000000 pt)=(0.792000,1.000000,0.208000);
rgb(573.000000 pt)=(0.796000,1.000000,0.204000);
rgb(574.000000 pt)=(0.800000,1.000000,0.200000);
rgb(575.000000 pt)=(0.804000,1.000000,0.196000);
rgb(576.000000 pt)=(0.808000,1.000000,0.192000);
rgb(577.000000 pt)=(0.812000,1.000000,0.188000);
rgb(578.000000 pt)=(0.816000,1.000000,0.184000);
rgb(579.000000 pt)=(0.820000,1.000000,0.180000);
rgb(580.000000 pt)=(0.824000,1.000000,0.176000);
rgb(581.000000 pt)=(0.828000,1.000000,0.172000);
rgb(582.000000 pt)=(0.832000,1.000000,0.168000);
rgb(583.000000 pt)=(0.836000,1.000000,0.164000);
rgb(584.000000 pt)=(0.840000,1.000000,0.160000);
rgb(585.000000 pt)=(0.844000,1.000000,0.156000);
rgb(586.000000 pt)=(0.848000,1.000000,0.152000);
rgb(587.000000 pt)=(0.852000,1.000000,0.148000);
rgb(588.000000 pt)=(0.856000,1.000000,0.144000);
rgb(589.000000 pt)=(0.860000,1.000000,0.140000);
rgb(590.000000 pt)=(0.864000,1.000000,0.136000);
rgb(591.000000 pt)=(0.868000,1.000000,0.132000);
rgb(592.000000 pt)=(0.872000,1.000000,0.128000);
rgb(593.000000 pt)=(0.876000,1.000000,0.124000);
rgb(594.000000 pt)=(0.880000,1.000000,0.120000);
rgb(595.000000 pt)=(0.884000,1.000000,0.116000);
rgb(596.000000 pt)=(0.888000,1.000000,0.112000);
rgb(597.000000 pt)=(0.892000,1.000000,0.108000);
rgb(598.000000 pt)=(0.896000,1.000000,0.104000);
rgb(599.000000 pt)=(0.900000,1.000000,0.100000);
rgb(600.000000 pt)=(0.904000,1.000000,0.096000);
rgb(601.000000 pt)=(0.908000,1.000000,0.092000);
rgb(602.000000 pt)=(0.912000,1.000000,0.088000);
rgb(603.000000 pt)=(0.916000,1.000000,0.084000);
rgb(604.000000 pt)=(0.920000,1.000000,0.080000);
rgb(605.000000 pt)=(0.924000,1.000000,0.076000);
rgb(606.000000 pt)=(0.928000,1.000000,0.072000);
rgb(607.000000 pt)=(0.932000,1.000000,0.068000);
rgb(608.000000 pt)=(0.936000,1.000000,0.064000);
rgb(609.000000 pt)=(0.940000,1.000000,0.060000);
rgb(610.000000 pt)=(0.944000,1.000000,0.056000);
rgb(611.000000 pt)=(0.948000,1.000000,0.052000);
rgb(612.000000 pt)=(0.952000,1.000000,0.048000);
rgb(613.000000 pt)=(0.956000,1.000000,0.044000);
rgb(614.000000 pt)=(0.960000,1.000000,0.040000);
rgb(615.000000 pt)=(0.964000,1.000000,0.036000);
rgb(616.000000 pt)=(0.968000,1.000000,0.032000);
rgb(617.000000 pt)=(0.972000,1.000000,0.028000);
rgb(618.000000 pt)=(0.976000,1.000000,0.024000);
rgb(619.000000 pt)=(0.980000,1.000000,0.020000);
rgb(620.000000 pt)=(0.984000,1.000000,0.016000);
rgb(621.000000 pt)=(0.988000,1.000000,0.012000);
rgb(622.000000 pt)=(0.992000,1.000000,0.008000);
rgb(623.000000 pt)=(0.996000,1.000000,0.004000);
rgb(624.000000 pt)=(1.000000,1.000000,0.000000);
rgb(625.000000 pt)=(1.000000,0.996000,0.000000);
rgb(626.000000 pt)=(1.000000,0.992000,0.000000);
rgb(627.000000 pt)=(1.000000,0.988000,0.000000);
rgb(628.000000 pt)=(1.000000,0.984000,0.000000);
rgb(629.000000 pt)=(1.000000,0.980000,0.000000);
rgb(630.000000 pt)=(1.000000,0.976000,0.000000);
rgb(631.000000 pt)=(1.000000,0.972000,0.000000);
rgb(632.000000 pt)=(1.000000,0.968000,0.000000);
rgb(633.000000 pt)=(1.000000,0.964000,0.000000);
rgb(634.000000 pt)=(1.000000,0.960000,0.000000);
rgb(635.000000 pt)=(1.000000,0.956000,0.000000);
rgb(636.000000 pt)=(1.000000,0.952000,0.000000);
rgb(637.000000 pt)=(1.000000,0.948000,0.000000);
rgb(638.000000 pt)=(1.000000,0.944000,0.000000);
rgb(639.000000 pt)=(1.000000,0.940000,0.000000);
rgb(640.000000 pt)=(1.000000,0.936000,0.000000);
rgb(641.000000 pt)=(1.000000,0.932000,0.000000);
rgb(642.000000 pt)=(1.000000,0.928000,0.000000);
rgb(643.000000 pt)=(1.000000,0.924000,0.000000);
rgb(644.000000 pt)=(1.000000,0.920000,0.000000);
rgb(645.000000 pt)=(1.000000,0.916000,0.000000);
rgb(646.000000 pt)=(1.000000,0.912000,0.000000);
rgb(647.000000 pt)=(1.000000,0.908000,0.000000);
rgb(648.000000 pt)=(1.000000,0.904000,0.000000);
rgb(649.000000 pt)=(1.000000,0.900000,0.000000);
rgb(650.000000 pt)=(1.000000,0.896000,0.000000);
rgb(651.000000 pt)=(1.000000,0.892000,0.000000);
rgb(652.000000 pt)=(1.000000,0.888000,0.000000);
rgb(653.000000 pt)=(1.000000,0.884000,0.000000);
rgb(654.000000 pt)=(1.000000,0.880000,0.000000);
rgb(655.000000 pt)=(1.000000,0.876000,0.000000);
rgb(656.000000 pt)=(1.000000,0.872000,0.000000);
rgb(657.000000 pt)=(1.000000,0.868000,0.000000);
rgb(658.000000 pt)=(1.000000,0.864000,0.000000);
rgb(659.000000 pt)=(1.000000,0.860000,0.000000);
rgb(660.000000 pt)=(1.000000,0.856000,0.000000);
rgb(661.000000 pt)=(1.000000,0.852000,0.000000);
rgb(662.000000 pt)=(1.000000,0.848000,0.000000);
rgb(663.000000 pt)=(1.000000,0.844000,0.000000);
rgb(664.000000 pt)=(1.000000,0.840000,0.000000);
rgb(665.000000 pt)=(1.000000,0.836000,0.000000);
rgb(666.000000 pt)=(1.000000,0.832000,0.000000);
rgb(667.000000 pt)=(1.000000,0.828000,0.000000);
rgb(668.000000 pt)=(1.000000,0.824000,0.000000);
rgb(669.000000 pt)=(1.000000,0.820000,0.000000);
rgb(670.000000 pt)=(1.000000,0.816000,0.000000);
rgb(671.000000 pt)=(1.000000,0.812000,0.000000);
rgb(672.000000 pt)=(1.000000,0.808000,0.000000);
rgb(673.000000 pt)=(1.000000,0.804000,0.000000);
rgb(674.000000 pt)=(1.000000,0.800000,0.000000);
rgb(675.000000 pt)=(1.000000,0.796000,0.000000);
rgb(676.000000 pt)=(1.000000,0.792000,0.000000);
rgb(677.000000 pt)=(1.000000,0.788000,0.000000);
rgb(678.000000 pt)=(1.000000,0.784000,0.000000);
rgb(679.000000 pt)=(1.000000,0.780000,0.000000);
rgb(680.000000 pt)=(1.000000,0.776000,0.000000);
rgb(681.000000 pt)=(1.000000,0.772000,0.000000);
rgb(682.000000 pt)=(1.000000,0.768000,0.000000);
rgb(683.000000 pt)=(1.000000,0.764000,0.000000);
rgb(684.000000 pt)=(1.000000,0.760000,0.000000);
rgb(685.000000 pt)=(1.000000,0.756000,0.000000);
rgb(686.000000 pt)=(1.000000,0.752000,0.000000);
rgb(687.000000 pt)=(1.000000,0.748000,0.000000);
rgb(688.000000 pt)=(1.000000,0.744000,0.000000);
rgb(689.000000 pt)=(1.000000,0.740000,0.000000);
rgb(690.000000 pt)=(1.000000,0.736000,0.000000);
rgb(691.000000 pt)=(1.000000,0.732000,0.000000);
rgb(692.000000 pt)=(1.000000,0.728000,0.000000);
rgb(693.000000 pt)=(1.000000,0.724000,0.000000);
rgb(694.000000 pt)=(1.000000,0.720000,0.000000);
rgb(695.000000 pt)=(1.000000,0.716000,0.000000);
rgb(696.000000 pt)=(1.000000,0.712000,0.000000);
rgb(697.000000 pt)=(1.000000,0.708000,0.000000);
rgb(698.000000 pt)=(1.000000,0.704000,0.000000);
rgb(699.000000 pt)=(1.000000,0.700000,0.000000);
rgb(700.000000 pt)=(1.000000,0.696000,0.000000);
rgb(701.000000 pt)=(1.000000,0.692000,0.000000);
rgb(702.000000 pt)=(1.000000,0.688000,0.000000);
rgb(703.000000 pt)=(1.000000,0.684000,0.000000);
rgb(704.000000 pt)=(1.000000,0.680000,0.000000);
rgb(705.000000 pt)=(1.000000,0.676000,0.000000);
rgb(706.000000 pt)=(1.000000,0.672000,0.000000);
rgb(707.000000 pt)=(1.000000,0.668000,0.000000);
rgb(708.000000 pt)=(1.000000,0.664000,0.000000);
rgb(709.000000 pt)=(1.000000,0.660000,0.000000);
rgb(710.000000 pt)=(1.000000,0.656000,0.000000);
rgb(711.000000 pt)=(1.000000,0.652000,0.000000);
rgb(712.000000 pt)=(1.000000,0.648000,0.000000);
rgb(713.000000 pt)=(1.000000,0.644000,0.000000);
rgb(714.000000 pt)=(1.000000,0.640000,0.000000);
rgb(715.000000 pt)=(1.000000,0.636000,0.000000);
rgb(716.000000 pt)=(1.000000,0.632000,0.000000);
rgb(717.000000 pt)=(1.000000,0.628000,0.000000);
rgb(718.000000 pt)=(1.000000,0.624000,0.000000);
rgb(719.000000 pt)=(1.000000,0.620000,0.000000);
rgb(720.000000 pt)=(1.000000,0.616000,0.000000);
rgb(721.000000 pt)=(1.000000,0.612000,0.000000);
rgb(722.000000 pt)=(1.000000,0.608000,0.000000);
rgb(723.000000 pt)=(1.000000,0.604000,0.000000);
rgb(724.000000 pt)=(1.000000,0.600000,0.000000);
rgb(725.000000 pt)=(1.000000,0.596000,0.000000);
rgb(726.000000 pt)=(1.000000,0.592000,0.000000);
rgb(727.000000 pt)=(1.000000,0.588000,0.000000);
rgb(728.000000 pt)=(1.000000,0.584000,0.000000);
rgb(729.000000 pt)=(1.000000,0.580000,0.000000);
rgb(730.000000 pt)=(1.000000,0.576000,0.000000);
rgb(731.000000 pt)=(1.000000,0.572000,0.000000);
rgb(732.000000 pt)=(1.000000,0.568000,0.000000);
rgb(733.000000 pt)=(1.000000,0.564000,0.000000);
rgb(734.000000 pt)=(1.000000,0.560000,0.000000);
rgb(735.000000 pt)=(1.000000,0.556000,0.000000);
rgb(736.000000 pt)=(1.000000,0.552000,0.000000);
rgb(737.000000 pt)=(1.000000,0.548000,0.000000);
rgb(738.000000 pt)=(1.000000,0.544000,0.000000);
rgb(739.000000 pt)=(1.000000,0.540000,0.000000);
rgb(740.000000 pt)=(1.000000,0.536000,0.000000);
rgb(741.000000 pt)=(1.000000,0.532000,0.000000);
rgb(742.000000 pt)=(1.000000,0.528000,0.000000);
rgb(743.000000 pt)=(1.000000,0.524000,0.000000);
rgb(744.000000 pt)=(1.000000,0.520000,0.000000);
rgb(745.000000 pt)=(1.000000,0.516000,0.000000);
rgb(746.000000 pt)=(1.000000,0.512000,0.000000);
rgb(747.000000 pt)=(1.000000,0.508000,0.000000);
rgb(748.000000 pt)=(1.000000,0.504000,0.000000);
rgb(749.000000 pt)=(1.000000,0.500000,0.000000);
rgb(750.000000 pt)=(1.000000,0.496000,0.000000);
rgb(751.000000 pt)=(1.000000,0.492000,0.000000);
rgb(752.000000 pt)=(1.000000,0.488000,0.000000);
rgb(753.000000 pt)=(1.000000,0.484000,0.000000);
rgb(754.000000 pt)=(1.000000,0.480000,0.000000);
rgb(755.000000 pt)=(1.000000,0.476000,0.000000);
rgb(756.000000 pt)=(1.000000,0.472000,0.000000);
rgb(757.000000 pt)=(1.000000,0.468000,0.000000);
rgb(758.000000 pt)=(1.000000,0.464000,0.000000);
rgb(759.000000 pt)=(1.000000,0.460000,0.000000);
rgb(760.000000 pt)=(1.000000,0.456000,0.000000);
rgb(761.000000 pt)=(1.000000,0.452000,0.000000);
rgb(762.000000 pt)=(1.000000,0.448000,0.000000);
rgb(763.000000 pt)=(1.000000,0.444000,0.000000);
rgb(764.000000 pt)=(1.000000,0.440000,0.000000);
rgb(765.000000 pt)=(1.000000,0.436000,0.000000);
rgb(766.000000 pt)=(1.000000,0.432000,0.000000);
rgb(767.000000 pt)=(1.000000,0.428000,0.000000);
rgb(768.000000 pt)=(1.000000,0.424000,0.000000);
rgb(769.000000 pt)=(1.000000,0.420000,0.000000);
rgb(770.000000 pt)=(1.000000,0.416000,0.000000);
rgb(771.000000 pt)=(1.000000,0.412000,0.000000);
rgb(772.000000 pt)=(1.000000,0.408000,0.000000);
rgb(773.000000 pt)=(1.000000,0.404000,0.000000);
rgb(774.000000 pt)=(1.000000,0.400000,0.000000);
rgb(775.000000 pt)=(1.000000,0.396000,0.000000);
rgb(776.000000 pt)=(1.000000,0.392000,0.000000);
rgb(777.000000 pt)=(1.000000,0.388000,0.000000);
rgb(778.000000 pt)=(1.000000,0.384000,0.000000);
rgb(779.000000 pt)=(1.000000,0.380000,0.000000);
rgb(780.000000 pt)=(1.000000,0.376000,0.000000);
rgb(781.000000 pt)=(1.000000,0.372000,0.000000);
rgb(782.000000 pt)=(1.000000,0.368000,0.000000);
rgb(783.000000 pt)=(1.000000,0.364000,0.000000);
rgb(784.000000 pt)=(1.000000,0.360000,0.000000);
rgb(785.000000 pt)=(1.000000,0.356000,0.000000);
rgb(786.000000 pt)=(1.000000,0.352000,0.000000);
rgb(787.000000 pt)=(1.000000,0.348000,0.000000);
rgb(788.000000 pt)=(1.000000,0.344000,0.000000);
rgb(789.000000 pt)=(1.000000,0.340000,0.000000);
rgb(790.000000 pt)=(1.000000,0.336000,0.000000);
rgb(791.000000 pt)=(1.000000,0.332000,0.000000);
rgb(792.000000 pt)=(1.000000,0.328000,0.000000);
rgb(793.000000 pt)=(1.000000,0.324000,0.000000);
rgb(794.000000 pt)=(1.000000,0.320000,0.000000);
rgb(795.000000 pt)=(1.000000,0.316000,0.000000);
rgb(796.000000 pt)=(1.000000,0.312000,0.000000);
rgb(797.000000 pt)=(1.000000,0.308000,0.000000);
rgb(798.000000 pt)=(1.000000,0.304000,0.000000);
rgb(799.000000 pt)=(1.000000,0.300000,0.000000);
rgb(800.000000 pt)=(1.000000,0.296000,0.000000);
rgb(801.000000 pt)=(1.000000,0.292000,0.000000);
rgb(802.000000 pt)=(1.000000,0.288000,0.000000);
rgb(803.000000 pt)=(1.000000,0.284000,0.000000);
rgb(804.000000 pt)=(1.000000,0.280000,0.000000);
rgb(805.000000 pt)=(1.000000,0.276000,0.000000);
rgb(806.000000 pt)=(1.000000,0.272000,0.000000);
rgb(807.000000 pt)=(1.000000,0.268000,0.000000);
rgb(808.000000 pt)=(1.000000,0.264000,0.000000);
rgb(809.000000 pt)=(1.000000,0.260000,0.000000);
rgb(810.000000 pt)=(1.000000,0.256000,0.000000);
rgb(811.000000 pt)=(1.000000,0.252000,0.000000);
rgb(812.000000 pt)=(1.000000,0.248000,0.000000);
rgb(813.000000 pt)=(1.000000,0.244000,0.000000);
rgb(814.000000 pt)=(1.000000,0.240000,0.000000);
rgb(815.000000 pt)=(1.000000,0.236000,0.000000);
rgb(816.000000 pt)=(1.000000,0.232000,0.000000);
rgb(817.000000 pt)=(1.000000,0.228000,0.000000);
rgb(818.000000 pt)=(1.000000,0.224000,0.000000);
rgb(819.000000 pt)=(1.000000,0.220000,0.000000);
rgb(820.000000 pt)=(1.000000,0.216000,0.000000);
rgb(821.000000 pt)=(1.000000,0.212000,0.000000);
rgb(822.000000 pt)=(1.000000,0.208000,0.000000);
rgb(823.000000 pt)=(1.000000,0.204000,0.000000);
rgb(824.000000 pt)=(1.000000,0.200000,0.000000);
rgb(825.000000 pt)=(1.000000,0.196000,0.000000);
rgb(826.000000 pt)=(1.000000,0.192000,0.000000);
rgb(827.000000 pt)=(1.000000,0.188000,0.000000);
rgb(828.000000 pt)=(1.000000,0.184000,0.000000);
rgb(829.000000 pt)=(1.000000,0.180000,0.000000);
rgb(830.000000 pt)=(1.000000,0.176000,0.000000);
rgb(831.000000 pt)=(1.000000,0.172000,0.000000);
rgb(832.000000 pt)=(1.000000,0.168000,0.000000);
rgb(833.000000 pt)=(1.000000,0.164000,0.000000);
rgb(834.000000 pt)=(1.000000,0.160000,0.000000);
rgb(835.000000 pt)=(1.000000,0.156000,0.000000);
rgb(836.000000 pt)=(1.000000,0.152000,0.000000);
rgb(837.000000 pt)=(1.000000,0.148000,0.000000);
rgb(838.000000 pt)=(1.000000,0.144000,0.000000);
rgb(839.000000 pt)=(1.000000,0.140000,0.000000);
rgb(840.000000 pt)=(1.000000,0.136000,0.000000);
rgb(841.000000 pt)=(1.000000,0.132000,0.000000);
rgb(842.000000 pt)=(1.000000,0.128000,0.000000);
rgb(843.000000 pt)=(1.000000,0.124000,0.000000);
rgb(844.000000 pt)=(1.000000,0.120000,0.000000);
rgb(845.000000 pt)=(1.000000,0.116000,0.000000);
rgb(846.000000 pt)=(1.000000,0.112000,0.000000);
rgb(847.000000 pt)=(1.000000,0.108000,0.000000);
rgb(848.000000 pt)=(1.000000,0.104000,0.000000);
rgb(849.000000 pt)=(1.000000,0.100000,0.000000);
rgb(850.000000 pt)=(1.000000,0.096000,0.000000);
rgb(851.000000 pt)=(1.000000,0.092000,0.000000);
rgb(852.000000 pt)=(1.000000,0.088000,0.000000);
rgb(853.000000 pt)=(1.000000,0.084000,0.000000);
rgb(854.000000 pt)=(1.000000,0.080000,0.000000);
rgb(855.000000 pt)=(1.000000,0.076000,0.000000);
rgb(856.000000 pt)=(1.000000,0.072000,0.000000);
rgb(857.000000 pt)=(1.000000,0.068000,0.000000);
rgb(858.000000 pt)=(1.000000,0.064000,0.000000);
rgb(859.000000 pt)=(1.000000,0.060000,0.000000);
rgb(860.000000 pt)=(1.000000,0.056000,0.000000);
rgb(861.000000 pt)=(1.000000,0.052000,0.000000);
rgb(862.000000 pt)=(1.000000,0.048000,0.000000);
rgb(863.000000 pt)=(1.000000,0.044000,0.000000);
rgb(864.000000 pt)=(1.000000,0.040000,0.000000);
rgb(865.000000 pt)=(1.000000,0.036000,0.000000);
rgb(866.000000 pt)=(1.000000,0.032000,0.000000);
rgb(867.000000 pt)=(1.000000,0.028000,0.000000);
rgb(868.000000 pt)=(1.000000,0.024000,0.000000);
rgb(869.000000 pt)=(1.000000,0.020000,0.000000);
rgb(870.000000 pt)=(1.000000,0.016000,0.000000);
rgb(871.000000 pt)=(1.000000,0.012000,0.000000);
rgb(872.000000 pt)=(1.000000,0.008000,0.000000);
rgb(873.000000 pt)=(1.000000,0.004000,0.000000);
rgb(874.000000 pt)=(1.000000,0.000000,0.000000);
rgb(875.000000 pt)=(0.996000,0.000000,0.000000);
rgb(876.000000 pt)=(0.992000,0.000000,0.000000);
rgb(877.000000 pt)=(0.988000,0.000000,0.000000);
rgb(878.000000 pt)=(0.984000,0.000000,0.000000);
rgb(879.000000 pt)=(0.980000,0.000000,0.000000);
rgb(880.000000 pt)=(0.976000,0.000000,0.000000);
rgb(881.000000 pt)=(0.972000,0.000000,0.000000);
rgb(882.000000 pt)=(0.968000,0.000000,0.000000);
rgb(883.000000 pt)=(0.964000,0.000000,0.000000);
rgb(884.000000 pt)=(0.960000,0.000000,0.000000);
rgb(885.000000 pt)=(0.956000,0.000000,0.000000);
rgb(886.000000 pt)=(0.952000,0.000000,0.000000);
rgb(887.000000 pt)=(0.948000,0.000000,0.000000);
rgb(888.000000 pt)=(0.944000,0.000000,0.000000);
rgb(889.000000 pt)=(0.940000,0.000000,0.000000);
rgb(890.000000 pt)=(0.936000,0.000000,0.000000);
rgb(891.000000 pt)=(0.932000,0.000000,0.000000);
rgb(892.000000 pt)=(0.928000,0.000000,0.000000);
rgb(893.000000 pt)=(0.924000,0.000000,0.000000);
rgb(894.000000 pt)=(0.920000,0.000000,0.000000);
rgb(895.000000 pt)=(0.916000,0.000000,0.000000);
rgb(896.000000 pt)=(0.912000,0.000000,0.000000);
rgb(897.000000 pt)=(0.908000,0.000000,0.000000);
rgb(898.000000 pt)=(0.904000,0.000000,0.000000);
rgb(899.000000 pt)=(0.900000,0.000000,0.000000);
rgb(900.000000 pt)=(0.896000,0.000000,0.000000);
rgb(901.000000 pt)=(0.892000,0.000000,0.000000);
rgb(902.000000 pt)=(0.888000,0.000000,0.000000);
rgb(903.000000 pt)=(0.884000,0.000000,0.000000);
rgb(904.000000 pt)=(0.880000,0.000000,0.000000);
rgb(905.000000 pt)=(0.876000,0.000000,0.000000);
rgb(906.000000 pt)=(0.872000,0.000000,0.000000);
rgb(907.000000 pt)=(0.868000,0.000000,0.000000);
rgb(908.000000 pt)=(0.864000,0.000000,0.000000);
rgb(909.000000 pt)=(0.860000,0.000000,0.000000);
rgb(910.000000 pt)=(0.856000,0.000000,0.000000);
rgb(911.000000 pt)=(0.852000,0.000000,0.000000);
rgb(912.000000 pt)=(0.848000,0.000000,0.000000);
rgb(913.000000 pt)=(0.844000,0.000000,0.000000);
rgb(914.000000 pt)=(0.840000,0.000000,0.000000);
rgb(915.000000 pt)=(0.836000,0.000000,0.000000);
rgb(916.000000 pt)=(0.832000,0.000000,0.000000);
rgb(917.000000 pt)=(0.828000,0.000000,0.000000);
rgb(918.000000 pt)=(0.824000,0.000000,0.000000);
rgb(919.000000 pt)=(0.820000,0.000000,0.000000);
rgb(920.000000 pt)=(0.816000,0.000000,0.000000);
rgb(921.000000 pt)=(0.812000,0.000000,0.000000);
rgb(922.000000 pt)=(0.808000,0.000000,0.000000);
rgb(923.000000 pt)=(0.804000,0.000000,0.000000);
rgb(924.000000 pt)=(0.800000,0.000000,0.000000);
rgb(925.000000 pt)=(0.796000,0.000000,0.000000);
rgb(926.000000 pt)=(0.792000,0.000000,0.000000);
rgb(927.000000 pt)=(0.788000,0.000000,0.000000);
rgb(928.000000 pt)=(0.784000,0.000000,0.000000);
rgb(929.000000 pt)=(0.780000,0.000000,0.000000);
rgb(930.000000 pt)=(0.776000,0.000000,0.000000);
rgb(931.000000 pt)=(0.772000,0.000000,0.000000);
rgb(932.000000 pt)=(0.768000,0.000000,0.000000);
rgb(933.000000 pt)=(0.764000,0.000000,0.000000);
rgb(934.000000 pt)=(0.760000,0.000000,0.000000);
rgb(935.000000 pt)=(0.756000,0.000000,0.000000);
rgb(936.000000 pt)=(0.752000,0.000000,0.000000);
rgb(937.000000 pt)=(0.748000,0.000000,0.000000);
rgb(938.000000 pt)=(0.744000,0.000000,0.000000);
rgb(939.000000 pt)=(0.740000,0.000000,0.000000);
rgb(940.000000 pt)=(0.736000,0.000000,0.000000);
rgb(941.000000 pt)=(0.732000,0.000000,0.000000);
rgb(942.000000 pt)=(0.728000,0.000000,0.000000);
rgb(943.000000 pt)=(0.724000,0.000000,0.000000);
rgb(944.000000 pt)=(0.720000,0.000000,0.000000);
rgb(945.000000 pt)=(0.716000,0.000000,0.000000);
rgb(946.000000 pt)=(0.712000,0.000000,0.000000);
rgb(947.000000 pt)=(0.708000,0.000000,0.000000);
rgb(948.000000 pt)=(0.704000,0.000000,0.000000);
rgb(949.000000 pt)=(0.700000,0.000000,0.000000);
rgb(950.000000 pt)=(0.696000,0.000000,0.000000);
rgb(951.000000 pt)=(0.692000,0.000000,0.000000);
rgb(952.000000 pt)=(0.688000,0.000000,0.000000);
rgb(953.000000 pt)=(0.684000,0.000000,0.000000);
rgb(954.000000 pt)=(0.680000,0.000000,0.000000);
rgb(955.000000 pt)=(0.676000,0.000000,0.000000);
rgb(956.000000 pt)=(0.672000,0.000000,0.000000);
rgb(957.000000 pt)=(0.668000,0.000000,0.000000);
rgb(958.000000 pt)=(0.664000,0.000000,0.000000);
rgb(959.000000 pt)=(0.660000,0.000000,0.000000);
rgb(960.000000 pt)=(0.656000,0.000000,0.000000);
rgb(961.000000 pt)=(0.652000,0.000000,0.000000);
rgb(962.000000 pt)=(0.648000,0.000000,0.000000);
rgb(963.000000 pt)=(0.644000,0.000000,0.000000);
rgb(964.000000 pt)=(0.640000,0.000000,0.000000);
rgb(965.000000 pt)=(0.636000,0.000000,0.000000);
rgb(966.000000 pt)=(0.632000,0.000000,0.000000);
rgb(967.000000 pt)=(0.628000,0.000000,0.000000);
rgb(968.000000 pt)=(0.624000,0.000000,0.000000);
rgb(969.000000 pt)=(0.620000,0.000000,0.000000);
rgb(970.000000 pt)=(0.616000,0.000000,0.000000);
rgb(971.000000 pt)=(0.612000,0.000000,0.000000);
rgb(972.000000 pt)=(0.608000,0.000000,0.000000);
rgb(973.000000 pt)=(0.604000,0.000000,0.000000);
rgb(974.000000 pt)=(0.600000,0.000000,0.000000);
rgb(975.000000 pt)=(0.596000,0.000000,0.000000);
rgb(976.000000 pt)=(0.592000,0.000000,0.000000);
rgb(977.000000 pt)=(0.588000,0.000000,0.000000);
rgb(978.000000 pt)=(0.584000,0.000000,0.000000);
rgb(979.000000 pt)=(0.580000,0.000000,0.000000);
rgb(980.000000 pt)=(0.576000,0.000000,0.000000);
rgb(981.000000 pt)=(0.572000,0.000000,0.000000);
rgb(982.000000 pt)=(0.568000,0.000000,0.000000);
rgb(983.000000 pt)=(0.564000,0.000000,0.000000);
rgb(984.000000 pt)=(0.560000,0.000000,0.000000);
rgb(985.000000 pt)=(0.556000,0.000000,0.000000);
rgb(986.000000 pt)=(0.552000,0.000000,0.000000);
rgb(987.000000 pt)=(0.548000,0.000000,0.000000);
rgb(988.000000 pt)=(0.544000,0.000000,0.000000);
rgb(989.000000 pt)=(0.540000,0.000000,0.000000);
rgb(990.000000 pt)=(0.536000,0.000000,0.000000);
rgb(991.000000 pt)=(0.532000,0.000000,0.000000);
rgb(992.000000 pt)=(0.528000,0.000000,0.000000);
rgb(993.000000 pt)=(0.524000,0.000000,0.000000);
rgb(994.000000 pt)=(0.520000,0.000000,0.000000);
rgb(995.000000 pt)=(0.516000,0.000000,0.000000);
rgb(996.000000 pt)=(0.512000,0.000000,0.000000);
rgb(997.000000 pt)=(0.508000,0.000000,0.000000);
rgb(998.000000 pt)=(0.504000,0.000000,0.000000);
rgb(999.000000 pt)=(0.500000,0.000000,0.000000);
}}

\pgfplotsset{
	colormap={myhot}{
		rgb(0pt)=(1.000000,1.000000,1.000000);
		rgb(1pt)=(1.000000,1.000000,0.996000);
		rgb(2pt)=(1.000000,1.000000,0.992000);
		rgb(3pt)=(1.000000,1.000000,0.988000);
		rgb(4pt)=(1.000000,1.000000,0.984000);
		rgb(5pt)=(1.000000,1.000000,0.980000);
		rgb(6pt)=(1.000000,1.000000,0.976000);
		rgb(7pt)=(1.000000,1.000000,0.972000);
		rgb(8pt)=(1.000000,1.000000,0.968000);
		rgb(9pt)=(1.000000,1.000000,0.964000);
		rgb(10pt)=(1.000000,1.000000,0.960000);
		rgb(11pt)=(1.000000,1.000000,0.956000);
		rgb(12pt)=(1.000000,1.000000,0.952000);
		rgb(13pt)=(1.000000,1.000000,0.948000);
		rgb(14pt)=(1.000000,1.000000,0.944000);
		rgb(15pt)=(1.000000,1.000000,0.940000);
		rgb(16pt)=(1.000000,1.000000,0.936000);
		rgb(17pt)=(1.000000,1.000000,0.932000);
		rgb(18pt)=(1.000000,1.000000,0.928000);
		rgb(19pt)=(1.000000,1.000000,0.924000);
		rgb(20pt)=(1.000000,1.000000,0.920000);
		rgb(21pt)=(1.000000,1.000000,0.916000);
		rgb(22pt)=(1.000000,1.000000,0.912000);
		rgb(23pt)=(1.000000,1.000000,0.908000);
		rgb(24pt)=(1.000000,1.000000,0.904000);
		rgb(25pt)=(1.000000,1.000000,0.900000);
		rgb(26pt)=(1.000000,1.000000,0.896000);
		rgb(27pt)=(1.000000,1.000000,0.892000);
		rgb(28pt)=(1.000000,1.000000,0.888000);
		rgb(29pt)=(1.000000,1.000000,0.884000);
		rgb(30pt)=(1.000000,1.000000,0.880000);
		rgb(31pt)=(1.000000,1.000000,0.876000);
		rgb(32pt)=(1.000000,1.000000,0.872000);
		rgb(33pt)=(1.000000,1.000000,0.868000);
		rgb(34pt)=(1.000000,1.000000,0.864000);
		rgb(35pt)=(1.000000,1.000000,0.860000);
		rgb(36pt)=(1.000000,1.000000,0.856000);
		rgb(37pt)=(1.000000,1.000000,0.852000);
		rgb(38pt)=(1.000000,1.000000,0.848000);
		rgb(39pt)=(1.000000,1.000000,0.844000);
		rgb(40pt)=(1.000000,1.000000,0.840000);
		rgb(41pt)=(1.000000,1.000000,0.836000);
		rgb(42pt)=(1.000000,1.000000,0.832000);
		rgb(43pt)=(1.000000,1.000000,0.828000);
		rgb(44pt)=(1.000000,1.000000,0.824000);
		rgb(45pt)=(1.000000,1.000000,0.820000);
		rgb(46pt)=(1.000000,1.000000,0.816000);
		rgb(47pt)=(1.000000,1.000000,0.812000);
		rgb(48pt)=(1.000000,1.000000,0.808000);
		rgb(49pt)=(1.000000,1.000000,0.804000);
		rgb(50pt)=(1.000000,1.000000,0.800000);
		rgb(51pt)=(1.000000,1.000000,0.796000);
		rgb(52pt)=(1.000000,1.000000,0.792000);
		rgb(53pt)=(1.000000,1.000000,0.788000);
		rgb(54pt)=(1.000000,1.000000,0.784000);
		rgb(55pt)=(1.000000,1.000000,0.780000);
		rgb(56pt)=(1.000000,1.000000,0.776000);
		rgb(57pt)=(1.000000,1.000000,0.772000);
		rgb(58pt)=(1.000000,1.000000,0.768000);
		rgb(59pt)=(1.000000,1.000000,0.764000);
		rgb(60pt)=(1.000000,1.000000,0.760000);
		rgb(61pt)=(1.000000,1.000000,0.756000);
		rgb(62pt)=(1.000000,1.000000,0.752000);
		rgb(63pt)=(1.000000,1.000000,0.748000);
		rgb(64pt)=(1.000000,1.000000,0.744000);
		rgb(65pt)=(1.000000,1.000000,0.740000);
		rgb(66pt)=(1.000000,1.000000,0.736000);
		rgb(67pt)=(1.000000,1.000000,0.732000);
		rgb(68pt)=(1.000000,1.000000,0.728000);
		rgb(69pt)=(1.000000,1.000000,0.724000);
		rgb(70pt)=(1.000000,1.000000,0.720000);
		rgb(71pt)=(1.000000,1.000000,0.716000);
		rgb(72pt)=(1.000000,1.000000,0.712000);
		rgb(73pt)=(1.000000,1.000000,0.708000);
		rgb(74pt)=(1.000000,1.000000,0.704000);
		rgb(75pt)=(1.000000,1.000000,0.700000);
		rgb(76pt)=(1.000000,1.000000,0.696000);
		rgb(77pt)=(1.000000,1.000000,0.692000);
		rgb(78pt)=(1.000000,1.000000,0.688000);
		rgb(79pt)=(1.000000,1.000000,0.684000);
		rgb(80pt)=(1.000000,1.000000,0.680000);
		rgb(81pt)=(1.000000,1.000000,0.676000);
		rgb(82pt)=(1.000000,1.000000,0.672000);
		rgb(83pt)=(1.000000,1.000000,0.668000);
		rgb(84pt)=(1.000000,1.000000,0.664000);
		rgb(85pt)=(1.000000,1.000000,0.660000);
		rgb(86pt)=(1.000000,1.000000,0.656000);
		rgb(87pt)=(1.000000,1.000000,0.652000);
		rgb(88pt)=(1.000000,1.000000,0.648000);
		rgb(89pt)=(1.000000,1.000000,0.644000);
		rgb(90pt)=(1.000000,1.000000,0.640000);
		rgb(91pt)=(1.000000,1.000000,0.636000);
		rgb(92pt)=(1.000000,1.000000,0.632000);
		rgb(93pt)=(1.000000,1.000000,0.628000);
		rgb(94pt)=(1.000000,1.000000,0.624000);
		rgb(95pt)=(1.000000,1.000000,0.620000);
		rgb(96pt)=(1.000000,1.000000,0.616000);
		rgb(97pt)=(1.000000,1.000000,0.612000);
		rgb(98pt)=(1.000000,1.000000,0.608000);
		rgb(99pt)=(1.000000,1.000000,0.604000);
		rgb(100pt)=(1.000000,1.000000,0.600000);
		rgb(101pt)=(1.000000,1.000000,0.596000);
		rgb(102pt)=(1.000000,1.000000,0.592000);
		rgb(103pt)=(1.000000,1.000000,0.588000);
		rgb(104pt)=(1.000000,1.000000,0.584000);
		rgb(105pt)=(1.000000,1.000000,0.580000);
		rgb(106pt)=(1.000000,1.000000,0.576000);
		rgb(107pt)=(1.000000,1.000000,0.572000);
		rgb(108pt)=(1.000000,1.000000,0.568000);
		rgb(109pt)=(1.000000,1.000000,0.564000);
		rgb(110pt)=(1.000000,1.000000,0.560000);
		rgb(111pt)=(1.000000,1.000000,0.556000);
		rgb(112pt)=(1.000000,1.000000,0.552000);
		rgb(113pt)=(1.000000,1.000000,0.548000);
		rgb(114pt)=(1.000000,1.000000,0.544000);
		rgb(115pt)=(1.000000,1.000000,0.540000);
		rgb(116pt)=(1.000000,1.000000,0.536000);
		rgb(117pt)=(1.000000,1.000000,0.532000);
		rgb(118pt)=(1.000000,1.000000,0.528000);
		rgb(119pt)=(1.000000,1.000000,0.524000);
		rgb(120pt)=(1.000000,1.000000,0.520000);
		rgb(121pt)=(1.000000,1.000000,0.516000);
		rgb(122pt)=(1.000000,1.000000,0.512000);
		rgb(123pt)=(1.000000,1.000000,0.508000);
		rgb(124pt)=(1.000000,1.000000,0.504000);
		rgb(125pt)=(1.000000,1.000000,0.500000);
		rgb(126pt)=(1.000000,1.000000,0.496000);
		rgb(127pt)=(1.000000,1.000000,0.492000);
		rgb(128pt)=(1.000000,1.000000,0.488000);
		rgb(129pt)=(1.000000,1.000000,0.484000);
		rgb(130pt)=(1.000000,1.000000,0.480000);
		rgb(131pt)=(1.000000,1.000000,0.476000);
		rgb(132pt)=(1.000000,1.000000,0.472000);
		rgb(133pt)=(1.000000,1.000000,0.468000);
		rgb(134pt)=(1.000000,1.000000,0.464000);
		rgb(135pt)=(1.000000,1.000000,0.460000);
		rgb(136pt)=(1.000000,1.000000,0.456000);
		rgb(137pt)=(1.000000,1.000000,0.452000);
		rgb(138pt)=(1.000000,1.000000,0.448000);
		rgb(139pt)=(1.000000,1.000000,0.444000);
		rgb(140pt)=(1.000000,1.000000,0.440000);
		rgb(141pt)=(1.000000,1.000000,0.436000);
		rgb(142pt)=(1.000000,1.000000,0.432000);
		rgb(143pt)=(1.000000,1.000000,0.428000);
		rgb(144pt)=(1.000000,1.000000,0.424000);
		rgb(145pt)=(1.000000,1.000000,0.420000);
		rgb(146pt)=(1.000000,1.000000,0.416000);
		rgb(147pt)=(1.000000,1.000000,0.412000);
		rgb(148pt)=(1.000000,1.000000,0.408000);
		rgb(149pt)=(1.000000,1.000000,0.404000);
		rgb(150pt)=(1.000000,1.000000,0.400000);
		rgb(151pt)=(1.000000,1.000000,0.396000);
		rgb(152pt)=(1.000000,1.000000,0.392000);
		rgb(153pt)=(1.000000,1.000000,0.388000);
		rgb(154pt)=(1.000000,1.000000,0.384000);
		rgb(155pt)=(1.000000,1.000000,0.380000);
		rgb(156pt)=(1.000000,1.000000,0.376000);
		rgb(157pt)=(1.000000,1.000000,0.372000);
		rgb(158pt)=(1.000000,1.000000,0.368000);
		rgb(159pt)=(1.000000,1.000000,0.364000);
		rgb(160pt)=(1.000000,1.000000,0.360000);
		rgb(161pt)=(1.000000,1.000000,0.356000);
		rgb(162pt)=(1.000000,1.000000,0.352000);
		rgb(163pt)=(1.000000,1.000000,0.348000);
		rgb(164pt)=(1.000000,1.000000,0.344000);
		rgb(165pt)=(1.000000,1.000000,0.340000);
		rgb(166pt)=(1.000000,1.000000,0.336000);
		rgb(167pt)=(1.000000,1.000000,0.332000);
		rgb(168pt)=(1.000000,1.000000,0.328000);
		rgb(169pt)=(1.000000,1.000000,0.324000);
		rgb(170pt)=(1.000000,1.000000,0.320000);
		rgb(171pt)=(1.000000,1.000000,0.316000);
		rgb(172pt)=(1.000000,1.000000,0.312000);
		rgb(173pt)=(1.000000,1.000000,0.308000);
		rgb(174pt)=(1.000000,1.000000,0.304000);
		rgb(175pt)=(1.000000,1.000000,0.300000);
		rgb(176pt)=(1.000000,1.000000,0.296000);
		rgb(177pt)=(1.000000,1.000000,0.292000);
		rgb(178pt)=(1.000000,1.000000,0.288000);
		rgb(179pt)=(1.000000,1.000000,0.284000);
		rgb(180pt)=(1.000000,1.000000,0.280000);
		rgb(181pt)=(1.000000,1.000000,0.276000);
		rgb(182pt)=(1.000000,1.000000,0.272000);
		rgb(183pt)=(1.000000,1.000000,0.268000);
		rgb(184pt)=(1.000000,1.000000,0.264000);
		rgb(185pt)=(1.000000,1.000000,0.260000);
		rgb(186pt)=(1.000000,1.000000,0.256000);
		rgb(187pt)=(1.000000,1.000000,0.252000);
		rgb(188pt)=(1.000000,1.000000,0.248000);
		rgb(189pt)=(1.000000,1.000000,0.244000);
		rgb(190pt)=(1.000000,1.000000,0.240000);
		rgb(191pt)=(1.000000,1.000000,0.236000);
		rgb(192pt)=(1.000000,1.000000,0.232000);
		rgb(193pt)=(1.000000,1.000000,0.228000);
		rgb(194pt)=(1.000000,1.000000,0.224000);
		rgb(195pt)=(1.000000,1.000000,0.220000);
		rgb(196pt)=(1.000000,1.000000,0.216000);
		rgb(197pt)=(1.000000,1.000000,0.212000);
		rgb(198pt)=(1.000000,1.000000,0.208000);
		rgb(199pt)=(1.000000,1.000000,0.204000);
		rgb(200pt)=(1.000000,1.000000,0.200000);
		rgb(201pt)=(1.000000,1.000000,0.196000);
		rgb(202pt)=(1.000000,1.000000,0.192000);
		rgb(203pt)=(1.000000,1.000000,0.188000);
		rgb(204pt)=(1.000000,1.000000,0.184000);
		rgb(205pt)=(1.000000,1.000000,0.180000);
		rgb(206pt)=(1.000000,1.000000,0.176000);
		rgb(207pt)=(1.000000,1.000000,0.172000);
		rgb(208pt)=(1.000000,1.000000,0.168000);
		rgb(209pt)=(1.000000,1.000000,0.164000);
		rgb(210pt)=(1.000000,1.000000,0.160000);
		rgb(211pt)=(1.000000,1.000000,0.156000);
		rgb(212pt)=(1.000000,1.000000,0.152000);
		rgb(213pt)=(1.000000,1.000000,0.148000);
		rgb(214pt)=(1.000000,1.000000,0.144000);
		rgb(215pt)=(1.000000,1.000000,0.140000);
		rgb(216pt)=(1.000000,1.000000,0.136000);
		rgb(217pt)=(1.000000,1.000000,0.132000);
		rgb(218pt)=(1.000000,1.000000,0.128000);
		rgb(219pt)=(1.000000,1.000000,0.124000);
		rgb(220pt)=(1.000000,1.000000,0.120000);
		rgb(221pt)=(1.000000,1.000000,0.116000);
		rgb(222pt)=(1.000000,1.000000,0.112000);
		rgb(223pt)=(1.000000,1.000000,0.108000);
		rgb(224pt)=(1.000000,1.000000,0.104000);
		rgb(225pt)=(1.000000,1.000000,0.100000);
		rgb(226pt)=(1.000000,1.000000,0.096000);
		rgb(227pt)=(1.000000,1.000000,0.092000);
		rgb(228pt)=(1.000000,1.000000,0.088000);
		rgb(229pt)=(1.000000,1.000000,0.084000);
		rgb(230pt)=(1.000000,1.000000,0.080000);
		rgb(231pt)=(1.000000,1.000000,0.076000);
		rgb(232pt)=(1.000000,1.000000,0.072000);
		rgb(233pt)=(1.000000,1.000000,0.068000);
		rgb(234pt)=(1.000000,1.000000,0.064000);
		rgb(235pt)=(1.000000,1.000000,0.060000);
		rgb(236pt)=(1.000000,1.000000,0.056000);
		rgb(237pt)=(1.000000,1.000000,0.052000);
		rgb(238pt)=(1.000000,1.000000,0.048000);
		rgb(239pt)=(1.000000,1.000000,0.044000);
		rgb(240pt)=(1.000000,1.000000,0.040000);
		rgb(241pt)=(1.000000,1.000000,0.036000);
		rgb(242pt)=(1.000000,1.000000,0.032000);
		rgb(243pt)=(1.000000,1.000000,0.028000);
		rgb(244pt)=(1.000000,1.000000,0.024000);
		rgb(245pt)=(1.000000,1.000000,0.020000);
		rgb(246pt)=(1.000000,1.000000,0.016000);
		rgb(247pt)=(1.000000,1.000000,0.012000);
		rgb(248pt)=(1.000000,1.000000,0.008000);
		rgb(249pt)=(1.000000,1.000000,0.004000);
		rgb(250pt)=(1.000000,1.000000,0.000000);
		rgb(251pt)=(1.000000,0.997333,0.000000);
		rgb(252pt)=(1.000000,0.994667,0.000000);
		rgb(253pt)=(1.000000,0.992000,0.000000);
		rgb(254pt)=(1.000000,0.989333,0.000000);
		rgb(255pt)=(1.000000,0.986667,0.000000);
		rgb(256pt)=(1.000000,0.984000,0.000000);
		rgb(257pt)=(1.000000,0.981333,0.000000);
		rgb(258pt)=(1.000000,0.978667,0.000000);
		rgb(259pt)=(1.000000,0.976000,0.000000);
		rgb(260pt)=(1.000000,0.973333,0.000000);
		rgb(261pt)=(1.000000,0.970667,0.000000);
		rgb(262pt)=(1.000000,0.968000,0.000000);
		rgb(263pt)=(1.000000,0.965333,0.000000);
		rgb(264pt)=(1.000000,0.962667,0.000000);
		rgb(265pt)=(1.000000,0.960000,0.000000);
		rgb(266pt)=(1.000000,0.957333,0.000000);
		rgb(267pt)=(1.000000,0.954667,0.000000);
		rgb(268pt)=(1.000000,0.952000,0.000000);
		rgb(269pt)=(1.000000,0.949333,0.000000);
		rgb(270pt)=(1.000000,0.946667,0.000000);
		rgb(271pt)=(1.000000,0.944000,0.000000);
		rgb(272pt)=(1.000000,0.941333,0.000000);
		rgb(273pt)=(1.000000,0.938667,0.000000);
		rgb(274pt)=(1.000000,0.936000,0.000000);
		rgb(275pt)=(1.000000,0.933333,0.000000);
		rgb(276pt)=(1.000000,0.930667,0.000000);
		rgb(277pt)=(1.000000,0.928000,0.000000);
		rgb(278pt)=(1.000000,0.925333,0.000000);
		rgb(279pt)=(1.000000,0.922667,0.000000);
		rgb(280pt)=(1.000000,0.920000,0.000000);
		rgb(281pt)=(1.000000,0.917333,0.000000);
		rgb(282pt)=(1.000000,0.914667,0.000000);
		rgb(283pt)=(1.000000,0.912000,0.000000);
		rgb(284pt)=(1.000000,0.909333,0.000000);
		rgb(285pt)=(1.000000,0.906667,0.000000);
		rgb(286pt)=(1.000000,0.904000,0.000000);
		rgb(287pt)=(1.000000,0.901333,0.000000);
		rgb(288pt)=(1.000000,0.898667,0.000000);
		rgb(289pt)=(1.000000,0.896000,0.000000);
		rgb(290pt)=(1.000000,0.893333,0.000000);
		rgb(291pt)=(1.000000,0.890667,0.000000);
		rgb(292pt)=(1.000000,0.888000,0.000000);
		rgb(293pt)=(1.000000,0.885333,0.000000);
		rgb(294pt)=(1.000000,0.882667,0.000000);
		rgb(295pt)=(1.000000,0.880000,0.000000);
		rgb(296pt)=(1.000000,0.877333,0.000000);
		rgb(297pt)=(1.000000,0.874667,0.000000);
		rgb(298pt)=(1.000000,0.872000,0.000000);
		rgb(299pt)=(1.000000,0.869333,0.000000);
		rgb(300pt)=(1.000000,0.866667,0.000000);
		rgb(301pt)=(1.000000,0.864000,0.000000);
		rgb(302pt)=(1.000000,0.861333,0.000000);
		rgb(303pt)=(1.000000,0.858667,0.000000);
		rgb(304pt)=(1.000000,0.856000,0.000000);
		rgb(305pt)=(1.000000,0.853333,0.000000);
		rgb(306pt)=(1.000000,0.850667,0.000000);
		rgb(307pt)=(1.000000,0.848000,0.000000);
		rgb(308pt)=(1.000000,0.845333,0.000000);
		rgb(309pt)=(1.000000,0.842667,0.000000);
		rgb(310pt)=(1.000000,0.840000,0.000000);
		rgb(311pt)=(1.000000,0.837333,0.000000);
		rgb(312pt)=(1.000000,0.834667,0.000000);
		rgb(313pt)=(1.000000,0.832000,0.000000);
		rgb(314pt)=(1.000000,0.829333,0.000000);
		rgb(315pt)=(1.000000,0.826667,0.000000);
		rgb(316pt)=(1.000000,0.824000,0.000000);
		rgb(317pt)=(1.000000,0.821333,0.000000);
		rgb(318pt)=(1.000000,0.818667,0.000000);
		rgb(319pt)=(1.000000,0.816000,0.000000);
		rgb(320pt)=(1.000000,0.813333,0.000000);
		rgb(321pt)=(1.000000,0.810667,0.000000);
		rgb(322pt)=(1.000000,0.808000,0.000000);
		rgb(323pt)=(1.000000,0.805333,0.000000);
		rgb(324pt)=(1.000000,0.802667,0.000000);
		rgb(325pt)=(1.000000,0.800000,0.000000);
		rgb(326pt)=(1.000000,0.797333,0.000000);
		rgb(327pt)=(1.000000,0.794667,0.000000);
		rgb(328pt)=(1.000000,0.792000,0.000000);
		rgb(329pt)=(1.000000,0.789333,0.000000);
		rgb(330pt)=(1.000000,0.786667,0.000000);
		rgb(331pt)=(1.000000,0.784000,0.000000);
		rgb(332pt)=(1.000000,0.781333,0.000000);
		rgb(333pt)=(1.000000,0.778667,0.000000);
		rgb(334pt)=(1.000000,0.776000,0.000000);
		rgb(335pt)=(1.000000,0.773333,0.000000);
		rgb(336pt)=(1.000000,0.770667,0.000000);
		rgb(337pt)=(1.000000,0.768000,0.000000);
		rgb(338pt)=(1.000000,0.765333,0.000000);
		rgb(339pt)=(1.000000,0.762667,0.000000);
		rgb(340pt)=(1.000000,0.760000,0.000000);
		rgb(341pt)=(1.000000,0.757333,0.000000);
		rgb(342pt)=(1.000000,0.754667,0.000000);
		rgb(343pt)=(1.000000,0.752000,0.000000);
		rgb(344pt)=(1.000000,0.749333,0.000000);
		rgb(345pt)=(1.000000,0.746667,0.000000);
		rgb(346pt)=(1.000000,0.744000,0.000000);
		rgb(347pt)=(1.000000,0.741333,0.000000);
		rgb(348pt)=(1.000000,0.738667,0.000000);
		rgb(349pt)=(1.000000,0.736000,0.000000);
		rgb(350pt)=(1.000000,0.733333,0.000000);
		rgb(351pt)=(1.000000,0.730667,0.000000);
		rgb(352pt)=(1.000000,0.728000,0.000000);
		rgb(353pt)=(1.000000,0.725333,0.000000);
		rgb(354pt)=(1.000000,0.722667,0.000000);
		rgb(355pt)=(1.000000,0.720000,0.000000);
		rgb(356pt)=(1.000000,0.717333,0.000000);
		rgb(357pt)=(1.000000,0.714667,0.000000);
		rgb(358pt)=(1.000000,0.712000,0.000000);
		rgb(359pt)=(1.000000,0.709333,0.000000);
		rgb(360pt)=(1.000000,0.706667,0.000000);
		rgb(361pt)=(1.000000,0.704000,0.000000);
		rgb(362pt)=(1.000000,0.701333,0.000000);
		rgb(363pt)=(1.000000,0.698667,0.000000);
		rgb(364pt)=(1.000000,0.696000,0.000000);
		rgb(365pt)=(1.000000,0.693333,0.000000);
		rgb(366pt)=(1.000000,0.690667,0.000000);
		rgb(367pt)=(1.000000,0.688000,0.000000);
		rgb(368pt)=(1.000000,0.685333,0.000000);
		rgb(369pt)=(1.000000,0.682667,0.000000);
		rgb(370pt)=(1.000000,0.680000,0.000000);
		rgb(371pt)=(1.000000,0.677333,0.000000);
		rgb(372pt)=(1.000000,0.674667,0.000000);
		rgb(373pt)=(1.000000,0.672000,0.000000);
		rgb(374pt)=(1.000000,0.669333,0.000000);
		rgb(375pt)=(1.000000,0.666667,0.000000);
		rgb(376pt)=(1.000000,0.664000,0.000000);
		rgb(377pt)=(1.000000,0.661333,0.000000);
		rgb(378pt)=(1.000000,0.658667,0.000000);
		rgb(379pt)=(1.000000,0.656000,0.000000);
		rgb(380pt)=(1.000000,0.653333,0.000000);
		rgb(381pt)=(1.000000,0.650667,0.000000);
		rgb(382pt)=(1.000000,0.648000,0.000000);
		rgb(383pt)=(1.000000,0.645333,0.000000);
		rgb(384pt)=(1.000000,0.642667,0.000000);
		rgb(385pt)=(1.000000,0.640000,0.000000);
		rgb(386pt)=(1.000000,0.637333,0.000000);
		rgb(387pt)=(1.000000,0.634667,0.000000);
		rgb(388pt)=(1.000000,0.632000,0.000000);
		rgb(389pt)=(1.000000,0.629333,0.000000);
		rgb(390pt)=(1.000000,0.626667,0.000000);
		rgb(391pt)=(1.000000,0.624000,0.000000);
		rgb(392pt)=(1.000000,0.621333,0.000000);
		rgb(393pt)=(1.000000,0.618667,0.000000);
		rgb(394pt)=(1.000000,0.616000,0.000000);
		rgb(395pt)=(1.000000,0.613333,0.000000);
		rgb(396pt)=(1.000000,0.610667,0.000000);
		rgb(397pt)=(1.000000,0.608000,0.000000);
		rgb(398pt)=(1.000000,0.605333,0.000000);
		rgb(399pt)=(1.000000,0.602667,0.000000);
		rgb(400pt)=(1.000000,0.600000,0.000000);
		rgb(401pt)=(1.000000,0.597333,0.000000);
		rgb(402pt)=(1.000000,0.594667,0.000000);
		rgb(403pt)=(1.000000,0.592000,0.000000);
		rgb(404pt)=(1.000000,0.589333,0.000000);
		rgb(405pt)=(1.000000,0.586667,0.000000);
		rgb(406pt)=(1.000000,0.584000,0.000000);
		rgb(407pt)=(1.000000,0.581333,0.000000);
		rgb(408pt)=(1.000000,0.578667,0.000000);
		rgb(409pt)=(1.000000,0.576000,0.000000);
		rgb(410pt)=(1.000000,0.573333,0.000000);
		rgb(411pt)=(1.000000,0.570667,0.000000);
		rgb(412pt)=(1.000000,0.568000,0.000000);
		rgb(413pt)=(1.000000,0.565333,0.000000);
		rgb(414pt)=(1.000000,0.562667,0.000000);
		rgb(415pt)=(1.000000,0.560000,0.000000);
		rgb(416pt)=(1.000000,0.557333,0.000000);
		rgb(417pt)=(1.000000,0.554667,0.000000);
		rgb(418pt)=(1.000000,0.552000,0.000000);
		rgb(419pt)=(1.000000,0.549333,0.000000);
		rgb(420pt)=(1.000000,0.546667,0.000000);
		rgb(421pt)=(1.000000,0.544000,0.000000);
		rgb(422pt)=(1.000000,0.541333,0.000000);
		rgb(423pt)=(1.000000,0.538667,0.000000);
		rgb(424pt)=(1.000000,0.536000,0.000000);
		rgb(425pt)=(1.000000,0.533333,0.000000);
		rgb(426pt)=(1.000000,0.530667,0.000000);
		rgb(427pt)=(1.000000,0.528000,0.000000);
		rgb(428pt)=(1.000000,0.525333,0.000000);
		rgb(429pt)=(1.000000,0.522667,0.000000);
		rgb(430pt)=(1.000000,0.520000,0.000000);
		rgb(431pt)=(1.000000,0.517333,0.000000);
		rgb(432pt)=(1.000000,0.514667,0.000000);
		rgb(433pt)=(1.000000,0.512000,0.000000);
		rgb(434pt)=(1.000000,0.509333,0.000000);
		rgb(435pt)=(1.000000,0.506667,0.000000);
		rgb(436pt)=(1.000000,0.504000,0.000000);
		rgb(437pt)=(1.000000,0.501333,0.000000);
		rgb(438pt)=(1.000000,0.498667,0.000000);
		rgb(439pt)=(1.000000,0.496000,0.000000);
		rgb(440pt)=(1.000000,0.493333,0.000000);
		rgb(441pt)=(1.000000,0.490667,0.000000);
		rgb(442pt)=(1.000000,0.488000,0.000000);
		rgb(443pt)=(1.000000,0.485333,0.000000);
		rgb(444pt)=(1.000000,0.482667,0.000000);
		rgb(445pt)=(1.000000,0.480000,0.000000);
		rgb(446pt)=(1.000000,0.477333,0.000000);
		rgb(447pt)=(1.000000,0.474667,0.000000);
		rgb(448pt)=(1.000000,0.472000,0.000000);
		rgb(449pt)=(1.000000,0.469333,0.000000);
		rgb(450pt)=(1.000000,0.466667,0.000000);
		rgb(451pt)=(1.000000,0.464000,0.000000);
		rgb(452pt)=(1.000000,0.461333,0.000000);
		rgb(453pt)=(1.000000,0.458667,0.000000);
		rgb(454pt)=(1.000000,0.456000,0.000000);
		rgb(455pt)=(1.000000,0.453333,0.000000);
		rgb(456pt)=(1.000000,0.450667,0.000000);
		rgb(457pt)=(1.000000,0.448000,0.000000);
		rgb(458pt)=(1.000000,0.445333,0.000000);
		rgb(459pt)=(1.000000,0.442667,0.000000);
		rgb(460pt)=(1.000000,0.440000,0.000000);
		rgb(461pt)=(1.000000,0.437333,0.000000);
		rgb(462pt)=(1.000000,0.434667,0.000000);
		rgb(463pt)=(1.000000,0.432000,0.000000);
		rgb(464pt)=(1.000000,0.429333,0.000000);
		rgb(465pt)=(1.000000,0.426667,0.000000);
		rgb(466pt)=(1.000000,0.424000,0.000000);
		rgb(467pt)=(1.000000,0.421333,0.000000);
		rgb(468pt)=(1.000000,0.418667,0.000000);
		rgb(469pt)=(1.000000,0.416000,0.000000);
		rgb(470pt)=(1.000000,0.413333,0.000000);
		rgb(471pt)=(1.000000,0.410667,0.000000);
		rgb(472pt)=(1.000000,0.408000,0.000000);
		rgb(473pt)=(1.000000,0.405333,0.000000);
		rgb(474pt)=(1.000000,0.402667,0.000000);
		rgb(475pt)=(1.000000,0.400000,0.000000);
		rgb(476pt)=(1.000000,0.397333,0.000000);
		rgb(477pt)=(1.000000,0.394667,0.000000);
		rgb(478pt)=(1.000000,0.392000,0.000000);
		rgb(479pt)=(1.000000,0.389333,0.000000);
		rgb(480pt)=(1.000000,0.386667,0.000000);
		rgb(481pt)=(1.000000,0.384000,0.000000);
		rgb(482pt)=(1.000000,0.381333,0.000000);
		rgb(483pt)=(1.000000,0.378667,0.000000);
		rgb(484pt)=(1.000000,0.376000,0.000000);
		rgb(485pt)=(1.000000,0.373333,0.000000);
		rgb(486pt)=(1.000000,0.370667,0.000000);
		rgb(487pt)=(1.000000,0.368000,0.000000);
		rgb(488pt)=(1.000000,0.365333,0.000000);
		rgb(489pt)=(1.000000,0.362667,0.000000);
		rgb(490pt)=(1.000000,0.360000,0.000000);
		rgb(491pt)=(1.000000,0.357333,0.000000);
		rgb(492pt)=(1.000000,0.354667,0.000000);
		rgb(493pt)=(1.000000,0.352000,0.000000);
		rgb(494pt)=(1.000000,0.349333,0.000000);
		rgb(495pt)=(1.000000,0.346667,0.000000);
		rgb(496pt)=(1.000000,0.344000,0.000000);
		rgb(497pt)=(1.000000,0.341333,0.000000);
		rgb(498pt)=(1.000000,0.338667,0.000000);
		rgb(499pt)=(1.000000,0.336000,0.000000);
		rgb(500pt)=(1.000000,0.333333,0.000000);
		rgb(501pt)=(1.000000,0.330667,0.000000);
		rgb(502pt)=(1.000000,0.328000,0.000000);
		rgb(503pt)=(1.000000,0.325333,0.000000);
		rgb(504pt)=(1.000000,0.322667,0.000000);
		rgb(505pt)=(1.000000,0.320000,0.000000);
		rgb(506pt)=(1.000000,0.317333,0.000000);
		rgb(507pt)=(1.000000,0.314667,0.000000);
		rgb(508pt)=(1.000000,0.312000,0.000000);
		rgb(509pt)=(1.000000,0.309333,0.000000);
		rgb(510pt)=(1.000000,0.306667,0.000000);
		rgb(511pt)=(1.000000,0.304000,0.000000);
		rgb(512pt)=(1.000000,0.301333,0.000000);
		rgb(513pt)=(1.000000,0.298667,0.000000);
		rgb(514pt)=(1.000000,0.296000,0.000000);
		rgb(515pt)=(1.000000,0.293333,0.000000);
		rgb(516pt)=(1.000000,0.290667,0.000000);
		rgb(517pt)=(1.000000,0.288000,0.000000);
		rgb(518pt)=(1.000000,0.285333,0.000000);
		rgb(519pt)=(1.000000,0.282667,0.000000);
		rgb(520pt)=(1.000000,0.280000,0.000000);
		rgb(521pt)=(1.000000,0.277333,0.000000);
		rgb(522pt)=(1.000000,0.274667,0.000000);
		rgb(523pt)=(1.000000,0.272000,0.000000);
		rgb(524pt)=(1.000000,0.269333,0.000000);
		rgb(525pt)=(1.000000,0.266667,0.000000);
		rgb(526pt)=(1.000000,0.264000,0.000000);
		rgb(527pt)=(1.000000,0.261333,0.000000);
		rgb(528pt)=(1.000000,0.258667,0.000000);
		rgb(529pt)=(1.000000,0.256000,0.000000);
		rgb(530pt)=(1.000000,0.253333,0.000000);
		rgb(531pt)=(1.000000,0.250667,0.000000);
		rgb(532pt)=(1.000000,0.248000,0.000000);
		rgb(533pt)=(1.000000,0.245333,0.000000);
		rgb(534pt)=(1.000000,0.242667,0.000000);
		rgb(535pt)=(1.000000,0.240000,0.000000);
		rgb(536pt)=(1.000000,0.237333,0.000000);
		rgb(537pt)=(1.000000,0.234667,0.000000);
		rgb(538pt)=(1.000000,0.232000,0.000000);
		rgb(539pt)=(1.000000,0.229333,0.000000);
		rgb(540pt)=(1.000000,0.226667,0.000000);
		rgb(541pt)=(1.000000,0.224000,0.000000);
		rgb(542pt)=(1.000000,0.221333,0.000000);
		rgb(543pt)=(1.000000,0.218667,0.000000);
		rgb(544pt)=(1.000000,0.216000,0.000000);
		rgb(545pt)=(1.000000,0.213333,0.000000);
		rgb(546pt)=(1.000000,0.210667,0.000000);
		rgb(547pt)=(1.000000,0.208000,0.000000);
		rgb(548pt)=(1.000000,0.205333,0.000000);
		rgb(549pt)=(1.000000,0.202667,0.000000);
		rgb(550pt)=(1.000000,0.200000,0.000000);
		rgb(551pt)=(1.000000,0.197333,0.000000);
		rgb(552pt)=(1.000000,0.194667,0.000000);
		rgb(553pt)=(1.000000,0.192000,0.000000);
		rgb(554pt)=(1.000000,0.189333,0.000000);
		rgb(555pt)=(1.000000,0.186667,0.000000);
		rgb(556pt)=(1.000000,0.184000,0.000000);
		rgb(557pt)=(1.000000,0.181333,0.000000);
		rgb(558pt)=(1.000000,0.178667,0.000000);
		rgb(559pt)=(1.000000,0.176000,0.000000);
		rgb(560pt)=(1.000000,0.173333,0.000000);
		rgb(561pt)=(1.000000,0.170667,0.000000);
		rgb(562pt)=(1.000000,0.168000,0.000000);
		rgb(563pt)=(1.000000,0.165333,0.000000);
		rgb(564pt)=(1.000000,0.162667,0.000000);
		rgb(565pt)=(1.000000,0.160000,0.000000);
		rgb(566pt)=(1.000000,0.157333,0.000000);
		rgb(567pt)=(1.000000,0.154667,0.000000);
		rgb(568pt)=(1.000000,0.152000,0.000000);
		rgb(569pt)=(1.000000,0.149333,0.000000);
		rgb(570pt)=(1.000000,0.146667,0.000000);
		rgb(571pt)=(1.000000,0.144000,0.000000);
		rgb(572pt)=(1.000000,0.141333,0.000000);
		rgb(573pt)=(1.000000,0.138667,0.000000);
		rgb(574pt)=(1.000000,0.136000,0.000000);
		rgb(575pt)=(1.000000,0.133333,0.000000);
		rgb(576pt)=(1.000000,0.130667,0.000000);
		rgb(577pt)=(1.000000,0.128000,0.000000);
		rgb(578pt)=(1.000000,0.125333,0.000000);
		rgb(579pt)=(1.000000,0.122667,0.000000);
		rgb(580pt)=(1.000000,0.120000,0.000000);
		rgb(581pt)=(1.000000,0.117333,0.000000);
		rgb(582pt)=(1.000000,0.114667,0.000000);
		rgb(583pt)=(1.000000,0.112000,0.000000);
		rgb(584pt)=(1.000000,0.109333,0.000000);
		rgb(585pt)=(1.000000,0.106667,0.000000);
		rgb(586pt)=(1.000000,0.104000,0.000000);
		rgb(587pt)=(1.000000,0.101333,0.000000);
		rgb(588pt)=(1.000000,0.098667,0.000000);
		rgb(589pt)=(1.000000,0.096000,0.000000);
		rgb(590pt)=(1.000000,0.093333,0.000000);
		rgb(591pt)=(1.000000,0.090667,0.000000);
		rgb(592pt)=(1.000000,0.088000,0.000000);
		rgb(593pt)=(1.000000,0.085333,0.000000);
		rgb(594pt)=(1.000000,0.082667,0.000000);
		rgb(595pt)=(1.000000,0.080000,0.000000);
		rgb(596pt)=(1.000000,0.077333,0.000000);
		rgb(597pt)=(1.000000,0.074667,0.000000);
		rgb(598pt)=(1.000000,0.072000,0.000000);
		rgb(599pt)=(1.000000,0.069333,0.000000);
		rgb(600pt)=(1.000000,0.066667,0.000000);
		rgb(601pt)=(1.000000,0.064000,0.000000);
		rgb(602pt)=(1.000000,0.061333,0.000000);
		rgb(603pt)=(1.000000,0.058667,0.000000);
		rgb(604pt)=(1.000000,0.056000,0.000000);
		rgb(605pt)=(1.000000,0.053333,0.000000);
		rgb(606pt)=(1.000000,0.050667,0.000000);
		rgb(607pt)=(1.000000,0.048000,0.000000);
		rgb(608pt)=(1.000000,0.045333,0.000000);
		rgb(609pt)=(1.000000,0.042667,0.000000);
		rgb(610pt)=(1.000000,0.040000,0.000000);
		rgb(611pt)=(1.000000,0.037333,0.000000);
		rgb(612pt)=(1.000000,0.034667,0.000000);
		rgb(613pt)=(1.000000,0.032000,0.000000);
		rgb(614pt)=(1.000000,0.029333,0.000000);
		rgb(615pt)=(1.000000,0.026667,0.000000);
		rgb(616pt)=(1.000000,0.024000,0.000000);
		rgb(617pt)=(1.000000,0.021333,0.000000);
		rgb(618pt)=(1.000000,0.018667,0.000000);
		rgb(619pt)=(1.000000,0.016000,0.000000);
		rgb(620pt)=(1.000000,0.013333,0.000000);
		rgb(621pt)=(1.000000,0.010667,0.000000);
		rgb(622pt)=(1.000000,0.008000,0.000000);
		rgb(623pt)=(1.000000,0.005333,0.000000);
		rgb(624pt)=(1.000000,0.002667,0.000000);
		rgb(625pt)=(1.000000,0.000000,0.000000);
		rgb(626pt)=(0.997333,0.000000,0.000000);
		rgb(627pt)=(0.994667,0.000000,0.000000);
		rgb(628pt)=(0.992000,0.000000,0.000000);
		rgb(629pt)=(0.989333,0.000000,0.000000);
		rgb(630pt)=(0.986667,0.000000,0.000000);
		rgb(631pt)=(0.984000,0.000000,0.000000);
		rgb(632pt)=(0.981333,0.000000,0.000000);
		rgb(633pt)=(0.978667,0.000000,0.000000);
		rgb(634pt)=(0.976000,0.000000,0.000000);
		rgb(635pt)=(0.973333,0.000000,0.000000);
		rgb(636pt)=(0.970667,0.000000,0.000000);
		rgb(637pt)=(0.968000,0.000000,0.000000);
		rgb(638pt)=(0.965333,0.000000,0.000000);
		rgb(639pt)=(0.962667,0.000000,0.000000);
		rgb(640pt)=(0.960000,0.000000,0.000000);
		rgb(641pt)=(0.957333,0.000000,0.000000);
		rgb(642pt)=(0.954667,0.000000,0.000000);
		rgb(643pt)=(0.952000,0.000000,0.000000);
		rgb(644pt)=(0.949333,0.000000,0.000000);
		rgb(645pt)=(0.946667,0.000000,0.000000);
		rgb(646pt)=(0.944000,0.000000,0.000000);
		rgb(647pt)=(0.941333,0.000000,0.000000);
		rgb(648pt)=(0.938667,0.000000,0.000000);
		rgb(649pt)=(0.936000,0.000000,0.000000);
		rgb(650pt)=(0.933333,0.000000,0.000000);
		rgb(651pt)=(0.930667,0.000000,0.000000);
		rgb(652pt)=(0.928000,0.000000,0.000000);
		rgb(653pt)=(0.925333,0.000000,0.000000);
		rgb(654pt)=(0.922667,0.000000,0.000000);
		rgb(655pt)=(0.920000,0.000000,0.000000);
		rgb(656pt)=(0.917333,0.000000,0.000000);
		rgb(657pt)=(0.914667,0.000000,0.000000);
		rgb(658pt)=(0.912000,0.000000,0.000000);
		rgb(659pt)=(0.909333,0.000000,0.000000);
		rgb(660pt)=(0.906667,0.000000,0.000000);
		rgb(661pt)=(0.904000,0.000000,0.000000);
		rgb(662pt)=(0.901333,0.000000,0.000000);
		rgb(663pt)=(0.898667,0.000000,0.000000);
		rgb(664pt)=(0.896000,0.000000,0.000000);
		rgb(665pt)=(0.893333,0.000000,0.000000);
		rgb(666pt)=(0.890667,0.000000,0.000000);
		rgb(667pt)=(0.888000,0.000000,0.000000);
		rgb(668pt)=(0.885333,0.000000,0.000000);
		rgb(669pt)=(0.882667,0.000000,0.000000);
		rgb(670pt)=(0.880000,0.000000,0.000000);
		rgb(671pt)=(0.877333,0.000000,0.000000);
		rgb(672pt)=(0.874667,0.000000,0.000000);
		rgb(673pt)=(0.872000,0.000000,0.000000);
		rgb(674pt)=(0.869333,0.000000,0.000000);
		rgb(675pt)=(0.866667,0.000000,0.000000);
		rgb(676pt)=(0.864000,0.000000,0.000000);
		rgb(677pt)=(0.861333,0.000000,0.000000);
		rgb(678pt)=(0.858667,0.000000,0.000000);
		rgb(679pt)=(0.856000,0.000000,0.000000);
		rgb(680pt)=(0.853333,0.000000,0.000000);
		rgb(681pt)=(0.850667,0.000000,0.000000);
		rgb(682pt)=(0.848000,0.000000,0.000000);
		rgb(683pt)=(0.845333,0.000000,0.000000);
		rgb(684pt)=(0.842667,0.000000,0.000000);
		rgb(685pt)=(0.840000,0.000000,0.000000);
		rgb(686pt)=(0.837333,0.000000,0.000000);
		rgb(687pt)=(0.834667,0.000000,0.000000);
		rgb(688pt)=(0.832000,0.000000,0.000000);
		rgb(689pt)=(0.829333,0.000000,0.000000);
		rgb(690pt)=(0.826667,0.000000,0.000000);
		rgb(691pt)=(0.824000,0.000000,0.000000);
		rgb(692pt)=(0.821333,0.000000,0.000000);
		rgb(693pt)=(0.818667,0.000000,0.000000);
		rgb(694pt)=(0.816000,0.000000,0.000000);
		rgb(695pt)=(0.813333,0.000000,0.000000);
		rgb(696pt)=(0.810667,0.000000,0.000000);
		rgb(697pt)=(0.808000,0.000000,0.000000);
		rgb(698pt)=(0.805333,0.000000,0.000000);
		rgb(699pt)=(0.802667,0.000000,0.000000);
		rgb(700pt)=(0.800000,0.000000,0.000000);
		rgb(701pt)=(0.797333,0.000000,0.000000);
		rgb(702pt)=(0.794667,0.000000,0.000000);
		rgb(703pt)=(0.792000,0.000000,0.000000);
		rgb(704pt)=(0.789333,0.000000,0.000000);
		rgb(705pt)=(0.786667,0.000000,0.000000);
		rgb(706pt)=(0.784000,0.000000,0.000000);
		rgb(707pt)=(0.781333,0.000000,0.000000);
		rgb(708pt)=(0.778667,0.000000,0.000000);
		rgb(709pt)=(0.776000,0.000000,0.000000);
		rgb(710pt)=(0.773333,0.000000,0.000000);
		rgb(711pt)=(0.770667,0.000000,0.000000);
		rgb(712pt)=(0.768000,0.000000,0.000000);
		rgb(713pt)=(0.765333,0.000000,0.000000);
		rgb(714pt)=(0.762667,0.000000,0.000000);
		rgb(715pt)=(0.760000,0.000000,0.000000);
		rgb(716pt)=(0.757333,0.000000,0.000000);
		rgb(717pt)=(0.754667,0.000000,0.000000);
		rgb(718pt)=(0.752000,0.000000,0.000000);
		rgb(719pt)=(0.749333,0.000000,0.000000);
		rgb(720pt)=(0.746667,0.000000,0.000000);
		rgb(721pt)=(0.744000,0.000000,0.000000);
		rgb(722pt)=(0.741333,0.000000,0.000000);
		rgb(723pt)=(0.738667,0.000000,0.000000);
		rgb(724pt)=(0.736000,0.000000,0.000000);
		rgb(725pt)=(0.733333,0.000000,0.000000);
		rgb(726pt)=(0.730667,0.000000,0.000000);
		rgb(727pt)=(0.728000,0.000000,0.000000);
		rgb(728pt)=(0.725333,0.000000,0.000000);
		rgb(729pt)=(0.722667,0.000000,0.000000);
		rgb(730pt)=(0.720000,0.000000,0.000000);
		rgb(731pt)=(0.717333,0.000000,0.000000);
		rgb(732pt)=(0.714667,0.000000,0.000000);
		rgb(733pt)=(0.712000,0.000000,0.000000);
		rgb(734pt)=(0.709333,0.000000,0.000000);
		rgb(735pt)=(0.706667,0.000000,0.000000);
		rgb(736pt)=(0.704000,0.000000,0.000000);
		rgb(737pt)=(0.701333,0.000000,0.000000);
		rgb(738pt)=(0.698667,0.000000,0.000000);
		rgb(739pt)=(0.696000,0.000000,0.000000);
		rgb(740pt)=(0.693333,0.000000,0.000000);
		rgb(741pt)=(0.690667,0.000000,0.000000);
		rgb(742pt)=(0.688000,0.000000,0.000000);
		rgb(743pt)=(0.685333,0.000000,0.000000);
		rgb(744pt)=(0.682667,0.000000,0.000000);
		rgb(745pt)=(0.680000,0.000000,0.000000);
		rgb(746pt)=(0.677333,0.000000,0.000000);
		rgb(747pt)=(0.674667,0.000000,0.000000);
		rgb(748pt)=(0.672000,0.000000,0.000000);
		rgb(749pt)=(0.669333,0.000000,0.000000);
		rgb(750pt)=(0.666667,0.000000,0.000000);
		rgb(751pt)=(0.664000,0.000000,0.000000);
		rgb(752pt)=(0.661333,0.000000,0.000000);
		rgb(753pt)=(0.658667,0.000000,0.000000);
		rgb(754pt)=(0.656000,0.000000,0.000000);
		rgb(755pt)=(0.653333,0.000000,0.000000);
		rgb(756pt)=(0.650667,0.000000,0.000000);
		rgb(757pt)=(0.648000,0.000000,0.000000);
		rgb(758pt)=(0.645333,0.000000,0.000000);
		rgb(759pt)=(0.642667,0.000000,0.000000);
		rgb(760pt)=(0.640000,0.000000,0.000000);
		rgb(761pt)=(0.637333,0.000000,0.000000);
		rgb(762pt)=(0.634667,0.000000,0.000000);
		rgb(763pt)=(0.632000,0.000000,0.000000);
		rgb(764pt)=(0.629333,0.000000,0.000000);
		rgb(765pt)=(0.626667,0.000000,0.000000);
		rgb(766pt)=(0.624000,0.000000,0.000000);
		rgb(767pt)=(0.621333,0.000000,0.000000);
		rgb(768pt)=(0.618667,0.000000,0.000000);
		rgb(769pt)=(0.616000,0.000000,0.000000);
		rgb(770pt)=(0.613333,0.000000,0.000000);
		rgb(771pt)=(0.610667,0.000000,0.000000);
		rgb(772pt)=(0.608000,0.000000,0.000000);
		rgb(773pt)=(0.605333,0.000000,0.000000);
		rgb(774pt)=(0.602667,0.000000,0.000000);
		rgb(775pt)=(0.600000,0.000000,0.000000);
		rgb(776pt)=(0.597333,0.000000,0.000000);
		rgb(777pt)=(0.594667,0.000000,0.000000);
		rgb(778pt)=(0.592000,0.000000,0.000000);
		rgb(779pt)=(0.589333,0.000000,0.000000);
		rgb(780pt)=(0.586667,0.000000,0.000000);
		rgb(781pt)=(0.584000,0.000000,0.000000);
		rgb(782pt)=(0.581333,0.000000,0.000000);
		rgb(783pt)=(0.578667,0.000000,0.000000);
		rgb(784pt)=(0.576000,0.000000,0.000000);
		rgb(785pt)=(0.573333,0.000000,0.000000);
		rgb(786pt)=(0.570667,0.000000,0.000000);
		rgb(787pt)=(0.568000,0.000000,0.000000);
		rgb(788pt)=(0.565333,0.000000,0.000000);
		rgb(789pt)=(0.562667,0.000000,0.000000);
		rgb(790pt)=(0.560000,0.000000,0.000000);
		rgb(791pt)=(0.557333,0.000000,0.000000);
		rgb(792pt)=(0.554667,0.000000,0.000000);
		rgb(793pt)=(0.552000,0.000000,0.000000);
		rgb(794pt)=(0.549333,0.000000,0.000000);
		rgb(795pt)=(0.546667,0.000000,0.000000);
		rgb(796pt)=(0.544000,0.000000,0.000000);
		rgb(797pt)=(0.541333,0.000000,0.000000);
		rgb(798pt)=(0.538667,0.000000,0.000000);
		rgb(799pt)=(0.536000,0.000000,0.000000);
		rgb(800pt)=(0.533333,0.000000,0.000000);
		rgb(801pt)=(0.530667,0.000000,0.000000);
		rgb(802pt)=(0.528000,0.000000,0.000000);
		rgb(803pt)=(0.525333,0.000000,0.000000);
		rgb(804pt)=(0.522667,0.000000,0.000000);
		rgb(805pt)=(0.520000,0.000000,0.000000);
		rgb(806pt)=(0.517333,0.000000,0.000000);
		rgb(807pt)=(0.514667,0.000000,0.000000);
		rgb(808pt)=(0.512000,0.000000,0.000000);
		rgb(809pt)=(0.509333,0.000000,0.000000);
		rgb(810pt)=(0.506667,0.000000,0.000000);
		rgb(811pt)=(0.504000,0.000000,0.000000);
		rgb(812pt)=(0.501333,0.000000,0.000000);
		rgb(813pt)=(0.498667,0.000000,0.000000);
		rgb(814pt)=(0.496000,0.000000,0.000000);
		rgb(815pt)=(0.493333,0.000000,0.000000);
		rgb(816pt)=(0.490667,0.000000,0.000000);
		rgb(817pt)=(0.488000,0.000000,0.000000);
		rgb(818pt)=(0.485333,0.000000,0.000000);
		rgb(819pt)=(0.482667,0.000000,0.000000);
		rgb(820pt)=(0.480000,0.000000,0.000000);
		rgb(821pt)=(0.477333,0.000000,0.000000);
		rgb(822pt)=(0.474667,0.000000,0.000000);
		rgb(823pt)=(0.472000,0.000000,0.000000);
		rgb(824pt)=(0.469333,0.000000,0.000000);
		rgb(825pt)=(0.466667,0.000000,0.000000);
		rgb(826pt)=(0.464000,0.000000,0.000000);
		rgb(827pt)=(0.461333,0.000000,0.000000);
		rgb(828pt)=(0.458667,0.000000,0.000000);
		rgb(829pt)=(0.456000,0.000000,0.000000);
		rgb(830pt)=(0.453333,0.000000,0.000000);
		rgb(831pt)=(0.450667,0.000000,0.000000);
		rgb(832pt)=(0.448000,0.000000,0.000000);
		rgb(833pt)=(0.445333,0.000000,0.000000);
		rgb(834pt)=(0.442667,0.000000,0.000000);
		rgb(835pt)=(0.440000,0.000000,0.000000);
		rgb(836pt)=(0.437333,0.000000,0.000000);
		rgb(837pt)=(0.434667,0.000000,0.000000);
		rgb(838pt)=(0.432000,0.000000,0.000000);
		rgb(839pt)=(0.429333,0.000000,0.000000);
		rgb(840pt)=(0.426667,0.000000,0.000000);
		rgb(841pt)=(0.424000,0.000000,0.000000);
		rgb(842pt)=(0.421333,0.000000,0.000000);
		rgb(843pt)=(0.418667,0.000000,0.000000);
		rgb(844pt)=(0.416000,0.000000,0.000000);
		rgb(845pt)=(0.413333,0.000000,0.000000);
		rgb(846pt)=(0.410667,0.000000,0.000000);
		rgb(847pt)=(0.408000,0.000000,0.000000);
		rgb(848pt)=(0.405333,0.000000,0.000000);
		rgb(849pt)=(0.402667,0.000000,0.000000);
		rgb(850pt)=(0.400000,0.000000,0.000000);
		rgb(851pt)=(0.397333,0.000000,0.000000);
		rgb(852pt)=(0.394667,0.000000,0.000000);
		rgb(853pt)=(0.392000,0.000000,0.000000);
		rgb(854pt)=(0.389333,0.000000,0.000000);
		rgb(855pt)=(0.386667,0.000000,0.000000);
		rgb(856pt)=(0.384000,0.000000,0.000000);
		rgb(857pt)=(0.381333,0.000000,0.000000);
		rgb(858pt)=(0.378667,0.000000,0.000000);
		rgb(859pt)=(0.376000,0.000000,0.000000);
		rgb(860pt)=(0.373333,0.000000,0.000000);
		rgb(861pt)=(0.370667,0.000000,0.000000);
		rgb(862pt)=(0.368000,0.000000,0.000000);
		rgb(863pt)=(0.365333,0.000000,0.000000);
		rgb(864pt)=(0.362667,0.000000,0.000000);
		rgb(865pt)=(0.360000,0.000000,0.000000);
		rgb(866pt)=(0.357333,0.000000,0.000000);
		rgb(867pt)=(0.354667,0.000000,0.000000);
		rgb(868pt)=(0.352000,0.000000,0.000000);
		rgb(869pt)=(0.349333,0.000000,0.000000);
		rgb(870pt)=(0.346667,0.000000,0.000000);
		rgb(871pt)=(0.344000,0.000000,0.000000);
		rgb(872pt)=(0.341333,0.000000,0.000000);
		rgb(873pt)=(0.338667,0.000000,0.000000);
		rgb(874pt)=(0.336000,0.000000,0.000000);
		rgb(875pt)=(0.333333,0.000000,0.000000);
		rgb(876pt)=(0.330667,0.000000,0.000000);
		rgb(877pt)=(0.328000,0.000000,0.000000);
		rgb(878pt)=(0.325333,0.000000,0.000000);
		rgb(879pt)=(0.322667,0.000000,0.000000);
		rgb(880pt)=(0.320000,0.000000,0.000000);
		rgb(881pt)=(0.317333,0.000000,0.000000);
		rgb(882pt)=(0.314667,0.000000,0.000000);
		rgb(883pt)=(0.312000,0.000000,0.000000);
		rgb(884pt)=(0.309333,0.000000,0.000000);
		rgb(885pt)=(0.306667,0.000000,0.000000);
		rgb(886pt)=(0.304000,0.000000,0.000000);
		rgb(887pt)=(0.301333,0.000000,0.000000);
		rgb(888pt)=(0.298667,0.000000,0.000000);
		rgb(889pt)=(0.296000,0.000000,0.000000);
		rgb(890pt)=(0.293333,0.000000,0.000000);
		rgb(891pt)=(0.290667,0.000000,0.000000);
		rgb(892pt)=(0.288000,0.000000,0.000000);
		rgb(893pt)=(0.285333,0.000000,0.000000);
		rgb(894pt)=(0.282667,0.000000,0.000000);
		rgb(895pt)=(0.280000,0.000000,0.000000);
		rgb(896pt)=(0.277333,0.000000,0.000000);
		rgb(897pt)=(0.274667,0.000000,0.000000);
		rgb(898pt)=(0.272000,0.000000,0.000000);
		rgb(899pt)=(0.269333,0.000000,0.000000);
		rgb(900pt)=(0.266667,0.000000,0.000000);
		rgb(901pt)=(0.264000,0.000000,0.000000);
		rgb(902pt)=(0.261333,0.000000,0.000000);
		rgb(903pt)=(0.258667,0.000000,0.000000);
		rgb(904pt)=(0.256000,0.000000,0.000000);
		rgb(905pt)=(0.253333,0.000000,0.000000);
		rgb(906pt)=(0.250667,0.000000,0.000000);
		rgb(907pt)=(0.248000,0.000000,0.000000);
		rgb(908pt)=(0.245333,0.000000,0.000000);
		rgb(909pt)=(0.242667,0.000000,0.000000);
		rgb(910pt)=(0.240000,0.000000,0.000000);
		rgb(911pt)=(0.237333,0.000000,0.000000);
		rgb(912pt)=(0.234667,0.000000,0.000000);
		rgb(913pt)=(0.232000,0.000000,0.000000);
		rgb(914pt)=(0.229333,0.000000,0.000000);
		rgb(915pt)=(0.226667,0.000000,0.000000);
		rgb(916pt)=(0.224000,0.000000,0.000000);
		rgb(917pt)=(0.221333,0.000000,0.000000);
		rgb(918pt)=(0.218667,0.000000,0.000000);
		rgb(919pt)=(0.216000,0.000000,0.000000);
		rgb(920pt)=(0.213333,0.000000,0.000000);
		rgb(921pt)=(0.210667,0.000000,0.000000);
		rgb(922pt)=(0.208000,0.000000,0.000000);
		rgb(923pt)=(0.205333,0.000000,0.000000);
		rgb(924pt)=(0.202667,0.000000,0.000000);
		rgb(925pt)=(0.200000,0.000000,0.000000);
		rgb(926pt)=(0.197333,0.000000,0.000000);
		rgb(927pt)=(0.194667,0.000000,0.000000);
		rgb(928pt)=(0.192000,0.000000,0.000000);
		rgb(929pt)=(0.189333,0.000000,0.000000);
		rgb(930pt)=(0.186667,0.000000,0.000000);
		rgb(931pt)=(0.184000,0.000000,0.000000);
		rgb(932pt)=(0.181333,0.000000,0.000000);
		rgb(933pt)=(0.178667,0.000000,0.000000);
		rgb(934pt)=(0.176000,0.000000,0.000000);
		rgb(935pt)=(0.173333,0.000000,0.000000);
		rgb(936pt)=(0.170667,0.000000,0.000000);
		rgb(937pt)=(0.168000,0.000000,0.000000);
		rgb(938pt)=(0.165333,0.000000,0.000000);
		rgb(939pt)=(0.162667,0.000000,0.000000);
		rgb(940pt)=(0.160000,0.000000,0.000000);
		rgb(941pt)=(0.157333,0.000000,0.000000);
		rgb(942pt)=(0.154667,0.000000,0.000000);
		rgb(943pt)=(0.152000,0.000000,0.000000);
		rgb(944pt)=(0.149333,0.000000,0.000000);
		rgb(945pt)=(0.146667,0.000000,0.000000);
		rgb(946pt)=(0.144000,0.000000,0.000000);
		rgb(947pt)=(0.141333,0.000000,0.000000);
		rgb(948pt)=(0.138667,0.000000,0.000000);
		rgb(949pt)=(0.136000,0.000000,0.000000);
		rgb(950pt)=(0.133333,0.000000,0.000000);
		rgb(951pt)=(0.130667,0.000000,0.000000);
		rgb(952pt)=(0.128000,0.000000,0.000000);
		rgb(953pt)=(0.125333,0.000000,0.000000);
		rgb(954pt)=(0.122667,0.000000,0.000000);
		rgb(955pt)=(0.120000,0.000000,0.000000);
		rgb(956pt)=(0.117333,0.000000,0.000000);
		rgb(957pt)=(0.114667,0.000000,0.000000);
		rgb(958pt)=(0.112000,0.000000,0.000000);
		rgb(959pt)=(0.109333,0.000000,0.000000);
		rgb(960pt)=(0.106667,0.000000,0.000000);
		rgb(961pt)=(0.104000,0.000000,0.000000);
		rgb(962pt)=(0.101333,0.000000,0.000000);
		rgb(963pt)=(0.098667,0.000000,0.000000);
		rgb(964pt)=(0.096000,0.000000,0.000000);
		rgb(965pt)=(0.093333,0.000000,0.000000);
		rgb(966pt)=(0.090667,0.000000,0.000000);
		rgb(967pt)=(0.088000,0.000000,0.000000);
		rgb(968pt)=(0.085333,0.000000,0.000000);
		rgb(969pt)=(0.082667,0.000000,0.000000);
		rgb(970pt)=(0.080000,0.000000,0.000000);
		rgb(971pt)=(0.077333,0.000000,0.000000);
		rgb(972pt)=(0.074667,0.000000,0.000000);
		rgb(973pt)=(0.072000,0.000000,0.000000);
		rgb(974pt)=(0.069333,0.000000,0.000000);
		rgb(975pt)=(0.066667,0.000000,0.000000);
		rgb(976pt)=(0.064000,0.000000,0.000000);
		rgb(977pt)=(0.061333,0.000000,0.000000);
		rgb(978pt)=(0.058667,0.000000,0.000000);
		rgb(979pt)=(0.056000,0.000000,0.000000);
		rgb(980pt)=(0.053333,0.000000,0.000000);
		rgb(981pt)=(0.050667,0.000000,0.000000);
		rgb(982pt)=(0.048000,0.000000,0.000000);
		rgb(983pt)=(0.045333,0.000000,0.000000);
		rgb(984pt)=(0.042667,0.000000,0.000000);
		rgb(985pt)=(0.040000,0.000000,0.000000);
		rgb(986pt)=(0.037333,0.000000,0.000000);
		rgb(987pt)=(0.034667,0.000000,0.000000);
		rgb(988pt)=(0.032000,0.000000,0.000000);
		rgb(989pt)=(0.029333,0.000000,0.000000);
		rgb(990pt)=(0.026667,0.000000,0.000000);
		rgb(991pt)=(0.024000,0.000000,0.000000);
		rgb(992pt)=(0.021333,0.000000,0.000000);
		rgb(993pt)=(0.018667,0.000000,0.000000);
		rgb(994pt)=(0.016000,0.000000,0.000000);
		rgb(995pt)=(0.013333,0.000000,0.000000);
		rgb(996pt)=(0.010667,0.000000,0.000000);
		rgb(997pt)=(0.008000,0.000000,0.000000);
		rgb(998pt)=(0.005333,0.000000,0.000000);
		rgb(999pt)=(0.002667,0.000000,0.000000);
}}

\newlength{\figureheight}
\newlength{\figurewidth}
\pgfplotsset{
	colormap={parula}{
		rgb(0pt)=(0.2081,0.1663,0.5292);
		rgb(1pt)=(0.208355,0.16778,0.532238);
		rgb(2pt)=(0.208611,0.169261,0.535275);
		rgb(3pt)=(0.208866,0.170741,0.538313);
		rgb(4pt)=(0.209121,0.172222,0.54135);
		rgb(5pt)=(0.209376,0.173702,0.544388);
		rgb(6pt)=(0.209632,0.175183,0.547425);
		rgb(7pt)=(0.209887,0.176663,0.550463);
		rgb(8pt)=(0.210134,0.178144,0.553505);
		rgb(9pt)=(0.210338,0.179624,0.556568);
		rgb(10pt)=(0.210542,0.181105,0.559631);
		rgb(11pt)=(0.210746,0.182585,0.562694);
		rgb(12pt)=(0.210944,0.184066,0.565763);
		rgb(13pt)=(0.211123,0.185546,0.568852);
		rgb(14pt)=(0.211302,0.187027,0.57194);
		rgb(15pt)=(0.21148,0.188507,0.575029);
		rgb(16pt)=(0.211642,0.189996,0.578117);
		rgb(17pt)=(0.21177,0.191502,0.581206);
		rgb(18pt)=(0.211897,0.193008,0.584295);
		rgb(19pt)=(0.212025,0.194514,0.587383);
		rgb(20pt)=(0.212132,0.19602,0.590472);
		rgb(21pt)=(0.212208,0.197526,0.59356);
		rgb(22pt)=(0.212285,0.199032,0.596649);
		rgb(23pt)=(0.212361,0.200538,0.599738);
		rgb(24pt)=(0.212413,0.202044,0.602839);
		rgb(25pt)=(0.212438,0.20355,0.605953);
		rgb(26pt)=(0.212464,0.205056,0.609067);
		rgb(27pt)=(0.212489,0.206562,0.612181);
		rgb(28pt)=(0.212471,0.208083,0.61531);
		rgb(29pt)=(0.21242,0.209614,0.61845);
		rgb(30pt)=(0.212368,0.211146,0.621589);
		rgb(31pt)=(0.212317,0.212677,0.624729);
		rgb(32pt)=(0.212216,0.214209,0.627868);
		rgb(33pt)=(0.212088,0.215741,0.631008);
		rgb(34pt)=(0.211961,0.217272,0.634148);
		rgb(35pt)=(0.211833,0.218804,0.637287);
		rgb(36pt)=(0.211668,0.220354,0.640446);
		rgb(37pt)=(0.211489,0.221911,0.643611);
		rgb(38pt)=(0.21131,0.223468,0.646776);
		rgb(39pt)=(0.211132,0.225025,0.649941);
		rgb(40pt)=(0.210848,0.226603,0.653107);
		rgb(41pt)=(0.210541,0.228186,0.656272);
		rgb(42pt)=(0.210235,0.229768,0.659437);
		rgb(43pt)=(0.209929,0.231351,0.662602);
		rgb(44pt)=(0.209553,0.232934,0.665767);
		rgb(45pt)=(0.20917,0.234516,0.668932);
		rgb(46pt)=(0.208787,0.236099,0.672098);
		rgb(47pt)=(0.208405,0.237681,0.675263);
		rgb(48pt)=(0.20787,0.239289,0.678453);
		rgb(49pt)=(0.207334,0.240897,0.681644);
		rgb(50pt)=(0.206798,0.242505,0.684835);
		rgb(51pt)=(0.206255,0.244114,0.688025);
		rgb(52pt)=(0.205617,0.245722,0.691216);
		rgb(53pt)=(0.204979,0.24733,0.694407);
		rgb(54pt)=(0.204341,0.248938,0.697597);
		rgb(55pt)=(0.203675,0.250554,0.700792);
		rgb(56pt)=(0.202858,0.252213,0.704008);
		rgb(57pt)=(0.202041,0.253872,0.707224);
		rgb(58pt)=(0.201225,0.255531,0.710441);
		rgb(59pt)=(0.200372,0.257184,0.713657);
		rgb(60pt)=(0.199402,0.258818,0.716873);
		rgb(61pt)=(0.198432,0.260452,0.720089);
		rgb(62pt)=(0.197462,0.262085,0.723305);
		rgb(63pt)=(0.196419,0.263735,0.726522);
		rgb(64pt)=(0.195219,0.26542,0.729738);
		rgb(65pt)=(0.19402,0.267105,0.732954);
		rgb(66pt)=(0.19282,0.268789,0.73617);
		rgb(67pt)=(0.191549,0.270474,0.739386);
		rgb(68pt)=(0.19017,0.272159,0.742603);
		rgb(69pt)=(0.188792,0.273843,0.745819);
		rgb(70pt)=(0.187414,0.275528,0.749035);
		rgb(71pt)=(0.1859,0.277237,0.752264);
		rgb(72pt)=(0.184241,0.278973,0.755505);
		rgb(73pt)=(0.182581,0.280709,0.758747);
		rgb(74pt)=(0.180922,0.282444,0.761989);
		rgb(75pt)=(0.179133,0.284209,0.765245);
		rgb(76pt)=(0.177244,0.285996,0.768512);
		rgb(77pt)=(0.175356,0.287783,0.77178);
		rgb(78pt)=(0.173467,0.289569,0.775047);
		rgb(79pt)=(0.171363,0.291406,0.778314);
		rgb(80pt)=(0.169142,0.293269,0.781581);
		rgb(81pt)=(0.166922,0.295132,0.784849);
		rgb(82pt)=(0.164701,0.296996,0.788116);
		rgb(83pt)=(0.162238,0.298934,0.791365);
		rgb(84pt)=(0.159686,0.300899,0.794606);
		rgb(85pt)=(0.157133,0.302865,0.797848);
		rgb(86pt)=(0.15458,0.30483,0.80109);
		rgb(87pt)=(0.151738,0.306858,0.804352);
		rgb(88pt)=(0.148828,0.3089,0.80762);
		rgb(89pt)=(0.145918,0.310942,0.810887);
		rgb(90pt)=(0.143008,0.312984,0.814154);
		rgb(91pt)=(0.139687,0.31514,0.81733);
		rgb(92pt)=(0.136318,0.31731,0.820495);
		rgb(93pt)=(0.132949,0.319479,0.82366);
		rgb(94pt)=(0.129579,0.321649,0.826826);
		rgb(95pt)=(0.125811,0.323918,0.829841);
		rgb(96pt)=(0.122033,0.32619,0.832853);
		rgb(97pt)=(0.118256,0.328462,0.835865);
		rgb(98pt)=(0.114458,0.330737,0.838862);
		rgb(99pt)=(0.110349,0.333059,0.841619);
		rgb(100pt)=(0.106239,0.335382,0.844376);
		rgb(101pt)=(0.102129,0.337705,0.847132);
		rgb(102pt)=(0.0979874,0.340021,0.849835);
		rgb(103pt)=(0.093648,0.342292,0.852209);
		rgb(104pt)=(0.0893087,0.344564,0.854583);
		rgb(105pt)=(0.0849694,0.346836,0.856957);
		rgb(106pt)=(0.08063,0.349091,0.859234);
		rgb(107pt)=(0.0762907,0.351286,0.861174);
		rgb(108pt)=(0.0719514,0.353481,0.863114);
		rgb(109pt)=(0.067612,0.355676,0.865053);
		rgb(110pt)=(0.0633195,0.357817,0.866853);
		rgb(111pt)=(0.0591333,0.359833,0.868333);
		rgb(112pt)=(0.0549471,0.36185,0.869814);
		rgb(113pt)=(0.050761,0.363866,0.871294);
		rgb(114pt)=(0.0466838,0.365823,0.872626);
		rgb(115pt)=(0.0427784,0.367687,0.873724);
		rgb(116pt)=(0.038873,0.36955,0.874821);
		rgb(117pt)=(0.0349676,0.371414,0.875919);
		rgb(118pt)=(0.0315066,0.373217,0.876872);
		rgb(119pt)=(0.0285456,0.374953,0.877664);
		rgb(120pt)=(0.0255847,0.376688,0.878455);
		rgb(121pt)=(0.0226237,0.378424,0.879246);
		rgb(122pt)=(0.0202132,0.380061,0.879868);
		rgb(123pt)=(0.0182477,0.381618,0.880353);
		rgb(124pt)=(0.0162823,0.383175,0.880838);
		rgb(125pt)=(0.0143168,0.384732,0.881323);
		rgb(126pt)=(0.0127892,0.386241,0.881695);
		rgb(127pt)=(0.0115129,0.387721,0.882001);
		rgb(128pt)=(0.0102366,0.389202,0.882307);
		rgb(129pt)=(0.00896036,0.390682,0.882614);
		rgb(130pt)=(0.00812372,0.392089,0.88281);
		rgb(131pt)=(0.00746006,0.393468,0.882963);
		rgb(132pt)=(0.0067964,0.394846,0.883116);
		rgb(133pt)=(0.00613273,0.396224,0.883269);
		rgb(134pt)=(0.00581622,0.397562,0.88332);
		rgb(135pt)=(0.00558649,0.398889,0.883346);
		rgb(136pt)=(0.00535676,0.400217,0.883371);
		rgb(137pt)=(0.00512703,0.401544,0.883397);
		rgb(138pt)=(0.00516757,0.402804,0.883332);
		rgb(139pt)=(0.00524414,0.404054,0.883256);
		rgb(140pt)=(0.00532072,0.405305,0.883179);
		rgb(141pt)=(0.0053973,0.406556,0.883103);
		rgb(142pt)=(0.00572012,0.407757,0.882952);
		rgb(143pt)=(0.00605195,0.408957,0.882799);
		rgb(144pt)=(0.00638378,0.410157,0.882646);
		rgb(145pt)=(0.00672643,0.411355,0.882489);
		rgb(146pt)=(0.00728799,0.412529,0.882259);
		rgb(147pt)=(0.00784955,0.413704,0.88203);
		rgb(148pt)=(0.00841111,0.414878,0.8818);
		rgb(149pt)=(0.00898919,0.416045,0.881564);
		rgb(150pt)=(0.00967838,0.417168,0.881283);
		rgb(151pt)=(0.0103676,0.418292,0.881002);
		rgb(152pt)=(0.0110568,0.419415,0.880721);
		rgb(153pt)=(0.011773,0.420532,0.880435);
		rgb(154pt)=(0.0125898,0.42163,0.880129);
		rgb(155pt)=(0.0134066,0.422728,0.879823);
		rgb(156pt)=(0.0142234,0.423825,0.879516);
		rgb(157pt)=(0.0150703,0.424915,0.879195);
		rgb(158pt)=(0.0159892,0.425987,0.878838);
		rgb(159pt)=(0.0169081,0.427059,0.87848);
		rgb(160pt)=(0.017827,0.428132,0.878123);
		rgb(161pt)=(0.0187748,0.429194,0.877746);
		rgb(162pt)=(0.0197703,0.430241,0.877338);
		rgb(163pt)=(0.0207658,0.431287,0.876929);
		rgb(164pt)=(0.0217613,0.432334,0.876521);
		rgb(165pt)=(0.0227802,0.43338,0.876113);
		rgb(166pt)=(0.0238267,0.434427,0.875704);
		rgb(167pt)=(0.0248733,0.435473,0.875296);
		rgb(168pt)=(0.0259198,0.43652,0.874887);
		rgb(169pt)=(0.0269802,0.437553,0.874451);
		rgb(170pt)=(0.0280523,0.438574,0.873992);
		rgb(171pt)=(0.0291243,0.439595,0.873532);
		rgb(172pt)=(0.0301964,0.440616,0.873073);
		rgb(173pt)=(0.0312844,0.441621,0.872614);
		rgb(174pt)=(0.032382,0.442616,0.872154);
		rgb(175pt)=(0.0334796,0.443612,0.871695);
		rgb(176pt)=(0.0345772,0.444607,0.871235);
		rgb(177pt)=(0.0357108,0.445603,0.870758);
		rgb(178pt)=(0.0368595,0.446598,0.870273);
		rgb(179pt)=(0.0380081,0.447594,0.869788);
		rgb(180pt)=(0.0391568,0.448589,0.869303);
		rgb(181pt)=(0.0402652,0.449565,0.868798);
		rgb(182pt)=(0.0413628,0.450535,0.868287);
		rgb(183pt)=(0.0424604,0.451505,0.867777);
		rgb(184pt)=(0.043558,0.452474,0.867266);
		rgb(185pt)=(0.0445889,0.453444,0.866756);
		rgb(186pt)=(0.0456099,0.454414,0.866245);
		rgb(187pt)=(0.0466309,0.455384,0.865735);
		rgb(188pt)=(0.047652,0.456354,0.865224);
		rgb(189pt)=(0.0486,0.457324,0.864714);
		rgb(190pt)=(0.0495444,0.458294,0.864203);
		rgb(191pt)=(0.0504889,0.459264,0.863692);
		rgb(192pt)=(0.0514315,0.460234,0.863181);
		rgb(193pt)=(0.0523249,0.461204,0.862645);
		rgb(194pt)=(0.0532183,0.462174,0.862109);
		rgb(195pt)=(0.0541117,0.463144,0.861573);
		rgb(196pt)=(0.0549991,0.464111,0.861034);
		rgb(197pt)=(0.0558414,0.465056,0.860472);
		rgb(198pt)=(0.0566838,0.466,0.859911);
		rgb(199pt)=(0.0575261,0.466944,0.859349);
		rgb(200pt)=(0.0583532,0.467889,0.858793);
		rgb(201pt)=(0.0591189,0.468833,0.858257);
		rgb(202pt)=(0.0598847,0.469778,0.857721);
		rgb(203pt)=(0.0606505,0.470722,0.857185);
		rgb(204pt)=(0.0614018,0.471667,0.856641);
		rgb(205pt)=(0.0621165,0.472611,0.85608);
		rgb(206pt)=(0.0628312,0.473556,0.855518);
		rgb(207pt)=(0.0635459,0.4745,0.854957);
		rgb(208pt)=(0.064242,0.475444,0.854405);
		rgb(209pt)=(0.0649057,0.476389,0.853868);
		rgb(210pt)=(0.0655694,0.477333,0.853332);
		rgb(211pt)=(0.066233,0.478278,0.852796);
		rgb(212pt)=(0.0668625,0.479222,0.852249);
		rgb(213pt)=(0.0674495,0.480167,0.851687);
		rgb(214pt)=(0.0680366,0.481111,0.851126);
		rgb(215pt)=(0.0686237,0.482056,0.850564);
		rgb(216pt)=(0.0691838,0.483,0.850003);
		rgb(217pt)=(0.0697198,0.483944,0.849441);
		rgb(218pt)=(0.0702559,0.484889,0.84888);
		rgb(219pt)=(0.0707919,0.485833,0.848318);
		rgb(220pt)=(0.0712967,0.486778,0.847772);
		rgb(221pt)=(0.0717817,0.487722,0.847236);
		rgb(222pt)=(0.0722667,0.488667,0.8467);
		rgb(223pt)=(0.0727517,0.489611,0.846164);
		rgb(224pt)=(0.0732012,0.490573,0.845628);
		rgb(225pt)=(0.0736351,0.491543,0.845092);
		rgb(226pt)=(0.0740691,0.492513,0.844556);
		rgb(227pt)=(0.074503,0.493483,0.84402);
		rgb(228pt)=(0.0748973,0.494433,0.843484);
		rgb(229pt)=(0.0752802,0.495378,0.842948);
		rgb(230pt)=(0.0756631,0.496322,0.842412);
		rgb(231pt)=(0.0760459,0.497267,0.841876);
		rgb(232pt)=(0.0763631,0.498233,0.841362);
		rgb(233pt)=(0.0766694,0.499203,0.840851);
		rgb(234pt)=(0.0769757,0.500173,0.840341);
		rgb(235pt)=(0.077282,0.501143,0.83983);
		rgb(236pt)=(0.0775162,0.502137,0.83932);
		rgb(237pt)=(0.0777459,0.503132,0.838809);
		rgb(238pt)=(0.0779757,0.504128,0.838298);
		rgb(239pt)=(0.0782042,0.505123,0.837789);
		rgb(240pt)=(0.0783829,0.506093,0.837304);
		rgb(241pt)=(0.0785616,0.507063,0.836819);
		rgb(242pt)=(0.0787402,0.508033,0.836334);
		rgb(243pt)=(0.0789135,0.509008,0.835851);
		rgb(244pt)=(0.0790411,0.510029,0.835392);
		rgb(245pt)=(0.0791688,0.51105,0.834932);
		rgb(246pt)=(0.0792964,0.512071,0.834473);
		rgb(247pt)=(0.0794048,0.513092,0.834018);
		rgb(248pt)=(0.0794303,0.514113,0.833584);
		rgb(249pt)=(0.0794559,0.515134,0.83315);
		rgb(250pt)=(0.0794814,0.516155,0.832717);
		rgb(251pt)=(0.0794862,0.517183,0.832289);
		rgb(252pt)=(0.0794351,0.51823,0.831881);
		rgb(253pt)=(0.0793841,0.519276,0.831473);
		rgb(254pt)=(0.079333,0.520323,0.831064);
		rgb(255pt)=(0.079255,0.521369,0.830665);
		rgb(256pt)=(0.0791273,0.522416,0.830282);
		rgb(257pt)=(0.0789997,0.523462,0.829899);
		rgb(258pt)=(0.0788721,0.524509,0.829516);
		rgb(259pt)=(0.0786889,0.525589,0.829156);
		rgb(260pt)=(0.0784336,0.526712,0.828824);
		rgb(261pt)=(0.0781784,0.527835,0.828492);
		rgb(262pt)=(0.0779231,0.528958,0.82816);
		rgb(263pt)=(0.077615,0.530081,0.827868);
		rgb(264pt)=(0.0772577,0.531205,0.827613);
		rgb(265pt)=(0.0769003,0.532328,0.827357);
		rgb(266pt)=(0.0765429,0.533451,0.827102);
		rgb(267pt)=(0.0761243,0.534589,0.826862);
		rgb(268pt)=(0.0756649,0.535738,0.826632);
		rgb(269pt)=(0.0752054,0.536886,0.826403);
		rgb(270pt)=(0.0747459,0.538035,0.826173);
		rgb(271pt)=(0.0742168,0.539219,0.825961);
		rgb(272pt)=(0.0736553,0.540418,0.825756);
		rgb(273pt)=(0.0730937,0.541618,0.825552);
		rgb(274pt)=(0.0725321,0.542818,0.825348);
		rgb(275pt)=(0.0718925,0.544037,0.825183);
		rgb(276pt)=(0.0712288,0.545262,0.82503);
		rgb(277pt)=(0.0705652,0.546487,0.824877);
		rgb(278pt)=(0.0699015,0.547713,0.824723);
		rgb(279pt)=(0.0691514,0.548938,0.824614);
		rgb(280pt)=(0.0683856,0.550163,0.824511);
		rgb(281pt)=(0.0676198,0.551388,0.824409);
		rgb(282pt)=(0.0668541,0.552614,0.824307);
		rgb(283pt)=(0.0660408,0.553886,0.824205);
		rgb(284pt)=(0.065224,0.555162,0.824103);
		rgb(285pt)=(0.0644072,0.556439,0.824001);
		rgb(286pt)=(0.0635892,0.557715,0.823899);
		rgb(287pt)=(0.0626703,0.558991,0.823848);
		rgb(288pt)=(0.0617514,0.560268,0.823797);
		rgb(289pt)=(0.0608324,0.561544,0.823746);
		rgb(290pt)=(0.0599087,0.56282,0.823693);
		rgb(291pt)=(0.0589387,0.564096,0.823616);
		rgb(292pt)=(0.0579688,0.565373,0.82354);
		rgb(293pt)=(0.0569988,0.566649,0.823463);
		rgb(294pt)=(0.0560243,0.567925,0.823386);
		rgb(295pt)=(0.0550288,0.569202,0.82331);
		rgb(296pt)=(0.0540333,0.570478,0.823233);
		rgb(297pt)=(0.0530378,0.571754,0.823157);
		rgb(298pt)=(0.0520423,0.57303,0.82308);
		rgb(299pt)=(0.0510468,0.574307,0.823004);
		rgb(300pt)=(0.0500514,0.575583,0.822927);
		rgb(301pt)=(0.0490559,0.576859,0.82285);
		rgb(302pt)=(0.0480604,0.578127,0.822756);
		rgb(303pt)=(0.0470649,0.579377,0.822629);
		rgb(304pt)=(0.0460694,0.580628,0.822501);
		rgb(305pt)=(0.0450739,0.581879,0.822374);
		rgb(306pt)=(0.0441,0.583119,0.822235);
		rgb(307pt)=(0.0431556,0.584344,0.822082);
		rgb(308pt)=(0.0422111,0.585569,0.821929);
		rgb(309pt)=(0.0412667,0.586795,0.821776);
		rgb(310pt)=(0.0403351,0.58802,0.821597);
		rgb(311pt)=(0.0394162,0.589245,0.821392);
		rgb(312pt)=(0.0384973,0.59047,0.821188);
		rgb(313pt)=(0.0375784,0.591695,0.820984);
		rgb(314pt)=(0.0367495,0.592891,0.820735);
		rgb(315pt)=(0.0359838,0.594065,0.820454);
		rgb(316pt)=(0.035218,0.595239,0.820173);
		rgb(317pt)=(0.0344523,0.596413,0.819892);
		rgb(318pt)=(0.0337721,0.597553,0.819595);
		rgb(319pt)=(0.0331339,0.598676,0.819288);
		rgb(320pt)=(0.0324958,0.599799,0.818982);
		rgb(321pt)=(0.0318577,0.600923,0.818676);
		rgb(322pt)=(0.0312964,0.602026,0.818312);
		rgb(323pt)=(0.0307604,0.603124,0.817929);
		rgb(324pt)=(0.0302243,0.604222,0.817546);
		rgb(325pt)=(0.0296883,0.605319,0.817163);
		rgb(326pt)=(0.0292375,0.606395,0.816738);
		rgb(327pt)=(0.0288036,0.607468,0.816304);
		rgb(328pt)=(0.0283697,0.60854,0.81587);
		rgb(329pt)=(0.0279357,0.609612,0.815436);
		rgb(330pt)=(0.0275721,0.610637,0.814955);
		rgb(331pt)=(0.0272147,0.611658,0.81447);
		rgb(332pt)=(0.0268574,0.612679,0.813985);
		rgb(333pt)=(0.0265,0.6137,0.8135);
		rgb(334pt)=(0.0262447,0.614695,0.812964);
		rgb(335pt)=(0.0259895,0.615691,0.812428);
		rgb(336pt)=(0.0257342,0.616686,0.811892);
		rgb(337pt)=(0.0254853,0.61768,0.811352);
		rgb(338pt)=(0.0253066,0.61865,0.810765);
		rgb(339pt)=(0.0251279,0.61962,0.810177);
		rgb(340pt)=(0.0249492,0.62059,0.80959);
		rgb(341pt)=(0.024779,0.621551,0.808995);
		rgb(342pt)=(0.0246514,0.62247,0.808357);
		rgb(343pt)=(0.0245237,0.623389,0.807719);
		rgb(344pt)=(0.0243961,0.624308,0.80708);
		rgb(345pt)=(0.0242748,0.625221,0.80643);
		rgb(346pt)=(0.0241727,0.626114,0.805741);
		rgb(347pt)=(0.0240706,0.627008,0.805051);
		rgb(348pt)=(0.0239685,0.627901,0.804362);
		rgb(349pt)=(0.0238832,0.628786,0.803656);
		rgb(350pt)=(0.0238321,0.629654,0.802916);
		rgb(351pt)=(0.0237811,0.630522,0.802176);
		rgb(352pt)=(0.02373,0.631389,0.801435);
		rgb(353pt)=(0.023679,0.632247,0.800685);
		rgb(354pt)=(0.0236279,0.633089,0.799919);
		rgb(355pt)=(0.0235769,0.633932,0.799153);
		rgb(356pt)=(0.0235258,0.634774,0.798387);
		rgb(357pt)=(0.0234748,0.635604,0.797596);
		rgb(358pt)=(0.0234237,0.63642,0.79678);
		rgb(359pt)=(0.0233727,0.637237,0.795963);
		rgb(360pt)=(0.0233216,0.638054,0.795146);
		rgb(361pt)=(0.0232706,0.638856,0.794329);
		rgb(362pt)=(0.0232195,0.639647,0.793512);
		rgb(363pt)=(0.0231685,0.640439,0.792695);
		rgb(364pt)=(0.0231174,0.64123,0.791879);
		rgb(365pt)=(0.0230832,0.642005,0.791011);
		rgb(366pt)=(0.0230577,0.64277,0.790118);
		rgb(367pt)=(0.0230321,0.643536,0.789225);
		rgb(368pt)=(0.0230066,0.644302,0.788331);
		rgb(369pt)=(0.0229811,0.645049,0.787438);
		rgb(370pt)=(0.0229556,0.645789,0.786544);
		rgb(371pt)=(0.02293,0.646529,0.785651);
		rgb(372pt)=(0.0229045,0.647269,0.784758);
		rgb(373pt)=(0.022858,0.64801,0.783843);
		rgb(374pt)=(0.0228069,0.64875,0.782924);
		rgb(375pt)=(0.0227559,0.64949,0.782005);
		rgb(376pt)=(0.0227048,0.65023,0.781086);
		rgb(377pt)=(0.0227,0.650947,0.780144);
		rgb(378pt)=(0.0227,0.651662,0.7792);
		rgb(379pt)=(0.0227,0.652377,0.778256);
		rgb(380pt)=(0.0227,0.653092,0.777311);
		rgb(381pt)=(0.0228261,0.653781,0.776341);
		rgb(382pt)=(0.0229538,0.65447,0.775371);
		rgb(383pt)=(0.0230814,0.655159,0.774402);
		rgb(384pt)=(0.0232108,0.655849,0.77343);
		rgb(385pt)=(0.023364,0.656538,0.772434);
		rgb(386pt)=(0.0235171,0.657227,0.771439);
		rgb(387pt)=(0.0236703,0.657916,0.770443);
		rgb(388pt)=(0.0238312,0.658602,0.769444);
		rgb(389pt)=(0.0240354,0.659265,0.768423);
		rgb(390pt)=(0.0242396,0.659929,0.767402);
		rgb(391pt)=(0.0244438,0.660592,0.766381);
		rgb(392pt)=(0.0247021,0.661256,0.765354);
		rgb(393pt)=(0.025136,0.66192,0.764307);
		rgb(394pt)=(0.02557,0.662583,0.763261);
		rgb(395pt)=(0.0260039,0.663247,0.762214);
		rgb(396pt)=(0.0264541,0.663911,0.761168);
		rgb(397pt)=(0.026939,0.664574,0.760121);
		rgb(398pt)=(0.027424,0.665238,0.759074);
		rgb(399pt)=(0.027909,0.665902,0.758028);
		rgb(400pt)=(0.028445,0.666555,0.756971);
		rgb(401pt)=(0.0290577,0.667193,0.755899);
		rgb(402pt)=(0.0296703,0.667832,0.754827);
		rgb(403pt)=(0.0302829,0.66847,0.753755);
		rgb(404pt)=(0.030994,0.669095,0.752683);
		rgb(405pt)=(0.0318108,0.669708,0.751611);
		rgb(406pt)=(0.0326276,0.670321,0.750539);
		rgb(407pt)=(0.0334444,0.670933,0.749467);
		rgb(408pt)=(0.0343045,0.67156,0.748366);
		rgb(409pt)=(0.0351979,0.672198,0.747243);
		rgb(410pt)=(0.0360913,0.672837,0.74612);
		rgb(411pt)=(0.0369847,0.673475,0.744996);
		rgb(412pt)=(0.0380432,0.674096,0.743873);
		rgb(413pt)=(0.0391919,0.674709,0.74275);
		rgb(414pt)=(0.0403405,0.675322,0.741627);
		rgb(415pt)=(0.0414892,0.675934,0.740504);
		rgb(416pt)=(0.0427123,0.676528,0.739381);
		rgb(417pt)=(0.0439631,0.677115,0.738258);
		rgb(418pt)=(0.0452138,0.677702,0.737135);
		rgb(419pt)=(0.0464646,0.678289,0.736011);
		rgb(420pt)=(0.0477153,0.678897,0.734868);
		rgb(421pt)=(0.0489661,0.67951,0.733719);
		rgb(422pt)=(0.0502168,0.680123,0.73257);
		rgb(423pt)=(0.0514676,0.680735,0.731422);
		rgb(424pt)=(0.0529237,0.681325,0.73025);
		rgb(425pt)=(0.0544042,0.681912,0.729076);
		rgb(426pt)=(0.0558847,0.682499,0.727902);
		rgb(427pt)=(0.0573652,0.683086,0.726728);
		rgb(428pt)=(0.0587709,0.683673,0.725553);
		rgb(429pt)=(0.0601748,0.68426,0.724379);
		rgb(430pt)=(0.0615787,0.684847,0.723205);
		rgb(431pt)=(0.0629946,0.685435,0.722028);
		rgb(432pt)=(0.0646027,0.686022,0.720803);
		rgb(433pt)=(0.0662108,0.686609,0.719577);
		rgb(434pt)=(0.0678189,0.687196,0.718352);
		rgb(435pt)=(0.069427,0.687779,0.717131);
		rgb(436pt)=(0.0710351,0.688341,0.715931);
		rgb(437pt)=(0.0726432,0.688902,0.714731);
		rgb(438pt)=(0.0742514,0.689464,0.713532);
		rgb(439pt)=(0.0758709,0.690026,0.712326);
		rgb(440pt)=(0.07753,0.690587,0.711101);
		rgb(441pt)=(0.0791892,0.691149,0.709876);
		rgb(442pt)=(0.0808483,0.69171,0.70865);
		rgb(443pt)=(0.0825387,0.692272,0.707417);
		rgb(444pt)=(0.0843,0.692833,0.706167);
		rgb(445pt)=(0.0860613,0.693395,0.704916);
		rgb(446pt)=(0.0878225,0.693956,0.703665);
		rgb(447pt)=(0.089564,0.694518,0.702405);
		rgb(448pt)=(0.0912742,0.69508,0.701128);
		rgb(449pt)=(0.0929844,0.695641,0.699852);
		rgb(450pt)=(0.0946946,0.696203,0.698576);
		rgb(451pt)=(0.0965009,0.696752,0.697299);
		rgb(452pt)=(0.0984153,0.697288,0.696023);
		rgb(453pt)=(0.10033,0.697824,0.694747);
		rgb(454pt)=(0.102244,0.69836,0.693471);
		rgb(455pt)=(0.10413,0.698896,0.69218);
		rgb(456pt)=(0.105994,0.699432,0.690878);
		rgb(457pt)=(0.107857,0.699968,0.689577);
		rgb(458pt)=(0.10972,0.700505,0.688275);
		rgb(459pt)=(0.111632,0.701041,0.686973);
		rgb(460pt)=(0.113572,0.701577,0.685671);
		rgb(461pt)=(0.115512,0.702113,0.684369);
		rgb(462pt)=(0.117452,0.702649,0.683068);
		rgb(463pt)=(0.119429,0.703185,0.681747);
		rgb(464pt)=(0.12142,0.703721,0.68042);
		rgb(465pt)=(0.123411,0.704257,0.679093);
		rgb(466pt)=(0.125402,0.704793,0.677765);
		rgb(467pt)=(0.127372,0.705308,0.676438);
		rgb(468pt)=(0.129338,0.705819,0.675111);
		rgb(469pt)=(0.131303,0.706329,0.673783);
		rgb(470pt)=(0.133269,0.70684,0.672456);
		rgb(471pt)=(0.135369,0.70735,0.671084);
		rgb(472pt)=(0.137488,0.707861,0.669705);
		rgb(473pt)=(0.139607,0.708371,0.668327);
		rgb(474pt)=(0.141725,0.708882,0.666949);
		rgb(475pt)=(0.143795,0.709392,0.665595);
		rgb(476pt)=(0.145862,0.709903,0.664242);
		rgb(477pt)=(0.14793,0.710414,0.662889);
		rgb(478pt)=(0.150003,0.710924,0.661534);
		rgb(479pt)=(0.152198,0.711435,0.66013);
		rgb(480pt)=(0.154394,0.711945,0.658726);
		rgb(481pt)=(0.156589,0.712456,0.657322);
		rgb(482pt)=(0.158784,0.712963,0.655922);
		rgb(483pt)=(0.160979,0.713448,0.654543);
		rgb(484pt)=(0.163174,0.713933,0.653165);
		rgb(485pt)=(0.16537,0.714418,0.651786);
		rgb(486pt)=(0.16757,0.714908,0.650397);
		rgb(487pt)=(0.169791,0.715419,0.648968);
		rgb(488pt)=(0.172012,0.715929,0.647538);
		rgb(489pt)=(0.174232,0.71644,0.646109);
		rgb(490pt)=(0.176483,0.716935,0.64468);
		rgb(491pt)=(0.178806,0.717395,0.64325);
		rgb(492pt)=(0.181129,0.717854,0.641821);
		rgb(493pt)=(0.183452,0.718314,0.640391);
		rgb(494pt)=(0.185755,0.718783,0.638952);
		rgb(495pt)=(0.188027,0.719268,0.637497);
		rgb(496pt)=(0.190299,0.719753,0.636042);
		rgb(497pt)=(0.192571,0.720238,0.634587);
		rgb(498pt)=(0.194913,0.720711,0.633132);
		rgb(499pt)=(0.197338,0.72117,0.631677);
		rgb(500pt)=(0.199762,0.72163,0.630223);
		rgb(501pt)=(0.202187,0.722089,0.628768);
		rgb(502pt)=(0.204612,0.722549,0.627299);
		rgb(503pt)=(0.207037,0.723008,0.625818);
		rgb(504pt)=(0.209462,0.723468,0.624338);
		rgb(505pt)=(0.211887,0.723927,0.622857);
		rgb(506pt)=(0.214328,0.724386,0.621377);
		rgb(507pt)=(0.216778,0.724846,0.619896);
		rgb(508pt)=(0.219229,0.725305,0.618416);
		rgb(509pt)=(0.221679,0.725765,0.616935);
		rgb(510pt)=(0.224202,0.726188,0.615455);
		rgb(511pt)=(0.226754,0.726597,0.613974);
		rgb(512pt)=(0.229307,0.727005,0.612494);
		rgb(513pt)=(0.231859,0.727414,0.611014);
		rgb(514pt)=(0.234392,0.727842,0.609513);
		rgb(515pt)=(0.236919,0.728276,0.608007);
		rgb(516pt)=(0.239446,0.72871,0.606501);
		rgb(517pt)=(0.241973,0.729144,0.604995);
		rgb(518pt)=(0.244611,0.729556,0.603467);
		rgb(519pt)=(0.247266,0.729964,0.601935);
		rgb(520pt)=(0.24992,0.730372,0.600404);
		rgb(521pt)=(0.252575,0.730781,0.598872);
		rgb(522pt)=(0.25523,0.731189,0.597365);
		rgb(523pt)=(0.257884,0.731598,0.595859);
		rgb(524pt)=(0.260539,0.732006,0.594353);
		rgb(525pt)=(0.263194,0.732414,0.592846);
		rgb(526pt)=(0.265848,0.732796,0.591314);
		rgb(527pt)=(0.268503,0.733179,0.589783);
		rgb(528pt)=(0.271158,0.733562,0.588251);
		rgb(529pt)=(0.27383,0.733945,0.58672);
		rgb(530pt)=(0.276638,0.734328,0.585188);
		rgb(531pt)=(0.279446,0.734711,0.583657);
		rgb(532pt)=(0.282254,0.735094,0.582125);
		rgb(533pt)=(0.285051,0.735471,0.580594);
		rgb(534pt)=(0.287808,0.735829,0.579062);
		rgb(535pt)=(0.290565,0.736186,0.577531);
		rgb(536pt)=(0.293322,0.736544,0.575999);
		rgb(537pt)=(0.2961,0.736894,0.574468);
		rgb(538pt)=(0.298933,0.737226,0.572936);
		rgb(539pt)=(0.301767,0.737557,0.571405);
		rgb(540pt)=(0.3046,0.737889,0.569873);
		rgb(541pt)=(0.307452,0.738221,0.568351);
		rgb(542pt)=(0.310336,0.738553,0.566845);
		rgb(543pt)=(0.313221,0.738885,0.565339);
		rgb(544pt)=(0.316105,0.739217,0.563833);
		rgb(545pt)=(0.318978,0.739537,0.562315);
		rgb(546pt)=(0.321837,0.739843,0.560784);
		rgb(547pt)=(0.324696,0.74015,0.559252);
		rgb(548pt)=(0.327555,0.740456,0.557721);
		rgb(549pt)=(0.330468,0.740749,0.556216);
		rgb(550pt)=(0.333429,0.741029,0.554736);
		rgb(551pt)=(0.336389,0.74131,0.553255);
		rgb(552pt)=(0.33935,0.741591,0.551775);
		rgb(553pt)=(0.342296,0.741872,0.550279);
		rgb(554pt)=(0.345231,0.742153,0.548773);
		rgb(555pt)=(0.348167,0.742433,0.547267);
		rgb(556pt)=(0.351102,0.742714,0.545761);
		rgb(557pt)=(0.354038,0.742977,0.54429);
		rgb(558pt)=(0.356973,0.743232,0.542835);
		rgb(559pt)=(0.359908,0.743488,0.54138);
		rgb(560pt)=(0.362844,0.743743,0.539925);
		rgb(561pt)=(0.365839,0.743959,0.53847);
		rgb(562pt)=(0.368851,0.744163,0.537015);
		rgb(563pt)=(0.371863,0.744367,0.53556);
		rgb(564pt)=(0.374875,0.744571,0.534105);
		rgb(565pt)=(0.377843,0.744775,0.532672);
		rgb(566pt)=(0.380804,0.74498,0.531243);
		rgb(567pt)=(0.383765,0.745184,0.529814);
		rgb(568pt)=(0.386726,0.745388,0.528384);
		rgb(569pt)=(0.389711,0.745568,0.527003);
		rgb(570pt)=(0.392697,0.745747,0.525624);
		rgb(571pt)=(0.395684,0.745926,0.524246);
		rgb(572pt)=(0.39867,0.746104,0.522868);
		rgb(573pt)=(0.401657,0.746257,0.521489);
		rgb(574pt)=(0.404643,0.74641,0.520111);
		rgb(575pt)=(0.40763,0.746563,0.518732);
		rgb(576pt)=(0.410611,0.746716,0.517359);
		rgb(577pt)=(0.413546,0.746869,0.516032);
		rgb(578pt)=(0.416482,0.747023,0.514705);
		rgb(579pt)=(0.419417,0.747176,0.513377);
		rgb(580pt)=(0.422357,0.747319,0.512055);
		rgb(581pt)=(0.425318,0.747421,0.510753);
		rgb(582pt)=(0.428279,0.747523,0.509451);
		rgb(583pt)=(0.43124,0.747626,0.50815);
		rgb(584pt)=(0.43418,0.747735,0.506848);
		rgb(585pt)=(0.437065,0.747862,0.505546);
		rgb(586pt)=(0.439949,0.74799,0.504244);
		rgb(587pt)=(0.442834,0.748117,0.502942);
		rgb(588pt)=(0.445727,0.748227,0.501659);
		rgb(589pt)=(0.448637,0.748304,0.500408);
		rgb(590pt)=(0.451547,0.74838,0.499157);
		rgb(591pt)=(0.454457,0.748457,0.497906);
		rgb(592pt)=(0.457333,0.748522,0.496667);
		rgb(593pt)=(0.460167,0.748573,0.495441);
		rgb(594pt)=(0.463,0.748624,0.494216);
		rgb(595pt)=(0.465833,0.748675,0.492991);
		rgb(596pt)=(0.468667,0.748726,0.491779);
		rgb(597pt)=(0.4715,0.748777,0.490579);
		rgb(598pt)=(0.474333,0.748829,0.48938);
		rgb(599pt)=(0.477167,0.74888,0.48818);
		rgb(600pt)=(0.479969,0.748931,0.486995);
		rgb(601pt)=(0.482752,0.748982,0.485821);
		rgb(602pt)=(0.485534,0.749033,0.484647);
		rgb(603pt)=(0.488316,0.749084,0.483473);
		rgb(604pt)=(0.491081,0.7491,0.482316);
		rgb(605pt)=(0.493838,0.7491,0.481168);
		rgb(606pt)=(0.496595,0.7491,0.480019);
		rgb(607pt)=(0.499351,0.7491,0.47887);
		rgb(608pt)=(0.502069,0.74912,0.477722);
		rgb(609pt)=(0.504775,0.749145,0.476573);
		rgb(610pt)=(0.50748,0.749171,0.475424);
		rgb(611pt)=(0.510186,0.749196,0.474276);
		rgb(612pt)=(0.512892,0.7492,0.47317);
		rgb(613pt)=(0.515598,0.7492,0.472073);
		rgb(614pt)=(0.518303,0.7492,0.470975);
		rgb(615pt)=(0.521009,0.7492,0.469877);
		rgb(616pt)=(0.523644,0.749176,0.46878);
		rgb(617pt)=(0.526273,0.749151,0.467682);
		rgb(618pt)=(0.528902,0.749125,0.466585);
		rgb(619pt)=(0.531531,0.7491,0.465487);
		rgb(620pt)=(0.53416,0.749074,0.464415);
		rgb(621pt)=(0.536789,0.749049,0.463343);
		rgb(622pt)=(0.539418,0.749023,0.462271);
		rgb(623pt)=(0.542043,0.748998,0.461199);
		rgb(624pt)=(0.544621,0.748972,0.460127);
		rgb(625pt)=(0.547199,0.748947,0.459055);
		rgb(626pt)=(0.549777,0.748921,0.457983);
		rgb(627pt)=(0.55235,0.748891,0.45692);
		rgb(628pt)=(0.554903,0.74884,0.455899);
		rgb(629pt)=(0.557456,0.748789,0.454878);
		rgb(630pt)=(0.560008,0.748738,0.453857);
		rgb(631pt)=(0.562554,0.748687,0.452829);
		rgb(632pt)=(0.565081,0.748636,0.451783);
		rgb(633pt)=(0.567608,0.748585,0.450736);
		rgb(634pt)=(0.570135,0.748534,0.449689);
		rgb(635pt)=(0.572653,0.748474,0.44866);
		rgb(636pt)=(0.575155,0.748397,0.447665);
		rgb(637pt)=(0.577656,0.748321,0.446669);
		rgb(638pt)=(0.580158,0.748244,0.445674);
		rgb(639pt)=(0.582649,0.748168,0.444678);
		rgb(640pt)=(0.585125,0.748091,0.443683);
		rgb(641pt)=(0.587601,0.748014,0.442687);
		rgb(642pt)=(0.590077,0.747938,0.441692);
		rgb(643pt)=(0.59254,0.747861,0.440709);
		rgb(644pt)=(0.59499,0.747785,0.439739);
		rgb(645pt)=(0.597441,0.747708,0.438769);
		rgb(646pt)=(0.599891,0.747632,0.437799);
		rgb(647pt)=(0.602311,0.747555,0.436814);
		rgb(648pt)=(0.604711,0.747478,0.435819);
		rgb(649pt)=(0.60711,0.747402,0.434823);
		rgb(650pt)=(0.60951,0.747325,0.433828);
		rgb(651pt)=(0.611909,0.747232,0.432867);
		rgb(652pt)=(0.614308,0.747129,0.431922);
		rgb(653pt)=(0.616708,0.747027,0.430978);
		rgb(654pt)=(0.619107,0.746925,0.430033);
		rgb(655pt)=(0.621487,0.746823,0.429089);
		rgb(656pt)=(0.623861,0.746721,0.428144);
		rgb(657pt)=(0.626235,0.746619,0.4272);
		rgb(658pt)=(0.628609,0.746517,0.426256);
		rgb(659pt)=(0.630962,0.746393,0.425311);
		rgb(660pt)=(0.63331,0.746266,0.424367);
		rgb(661pt)=(0.635658,0.746138,0.423422);
		rgb(662pt)=(0.638007,0.746011,0.422478);
		rgb(663pt)=(0.640332,0.745906,0.421557);
		rgb(664pt)=(0.642654,0.745804,0.420638);
		rgb(665pt)=(0.644977,0.745702,0.419719);
		rgb(666pt)=(0.6473,0.7456,0.4188);
		rgb(667pt)=(0.649623,0.745472,0.417881);
		rgb(668pt)=(0.651946,0.745345,0.416962);
		rgb(669pt)=(0.654268,0.745217,0.416043);
		rgb(670pt)=(0.656587,0.745089,0.415124);
		rgb(671pt)=(0.658859,0.744962,0.414205);
		rgb(672pt)=(0.661131,0.744834,0.413286);
		rgb(673pt)=(0.663402,0.744707,0.412368);
		rgb(674pt)=(0.665674,0.744579,0.411453);
		rgb(675pt)=(0.667946,0.744451,0.410559);
		rgb(676pt)=(0.670218,0.744324,0.409666);
		rgb(677pt)=(0.672489,0.744196,0.408773);
		rgb(678pt)=(0.674755,0.744062,0.407879);
		rgb(679pt)=(0.677001,0.743909,0.406986);
		rgb(680pt)=(0.679247,0.743756,0.406092);
		rgb(681pt)=(0.681494,0.743603,0.405199);
		rgb(682pt)=(0.68374,0.743458,0.404306);
		rgb(683pt)=(0.685986,0.74333,0.403412);
		rgb(684pt)=(0.688232,0.743203,0.402519);
		rgb(685pt)=(0.690479,0.743075,0.401626);
		rgb(686pt)=(0.692704,0.742937,0.400732);
		rgb(687pt)=(0.694899,0.742784,0.399839);
		rgb(688pt)=(0.697094,0.742631,0.398945);
		rgb(689pt)=(0.699289,0.742477,0.398052);
		rgb(690pt)=(0.701497,0.742324,0.397171);
		rgb(691pt)=(0.703718,0.742171,0.396303);
		rgb(692pt)=(0.705939,0.742018,0.395435);
		rgb(693pt)=(0.708159,0.741865,0.394568);
		rgb(694pt)=(0.710351,0.741712,0.3937);
		rgb(695pt)=(0.71252,0.741559,0.392832);
		rgb(696pt)=(0.71469,0.741405,0.391964);
		rgb(697pt)=(0.71686,0.741252,0.391096);
		rgb(698pt)=(0.719029,0.741082,0.390228);
		rgb(699pt)=(0.721199,0.740904,0.38936);
		rgb(700pt)=(0.723369,0.740725,0.388492);
		rgb(701pt)=(0.725538,0.740546,0.387625);
		rgb(702pt)=(0.727708,0.740386,0.386757);
		rgb(703pt)=(0.729878,0.740233,0.385889);
		rgb(704pt)=(0.732047,0.74008,0.385021);
		rgb(705pt)=(0.734217,0.739927,0.384153);
		rgb(706pt)=(0.736366,0.739753,0.383285);
		rgb(707pt)=(0.73851,0.739574,0.382417);
		rgb(708pt)=(0.740654,0.739395,0.38155);
		rgb(709pt)=(0.742798,0.739217,0.380682);
		rgb(710pt)=(0.744919,0.739038,0.379837);
		rgb(711pt)=(0.747038,0.738859,0.378995);
		rgb(712pt)=(0.749156,0.738681,0.378152);
		rgb(713pt)=(0.751275,0.738502,0.37731);
		rgb(714pt)=(0.753394,0.738323,0.376442);
		rgb(715pt)=(0.755512,0.738145,0.375574);
		rgb(716pt)=(0.757631,0.737966,0.374707);
		rgb(717pt)=(0.75975,0.737789,0.373841);
		rgb(718pt)=(0.761868,0.737636,0.372998);
		rgb(719pt)=(0.763987,0.737483,0.372156);
		rgb(720pt)=(0.766105,0.73733,0.371314);
		rgb(721pt)=(0.76822,0.737169,0.370471);
		rgb(722pt)=(0.770313,0.736965,0.369629);
		rgb(723pt)=(0.772406,0.73676,0.368786);
		rgb(724pt)=(0.774499,0.736556,0.367944);
		rgb(725pt)=(0.776592,0.736358,0.367096);
		rgb(726pt)=(0.778686,0.736179,0.366228);
		rgb(727pt)=(0.780779,0.736001,0.36536);
		rgb(728pt)=(0.782872,0.735822,0.364492);
		rgb(729pt)=(0.784957,0.735643,0.363632);
		rgb(730pt)=(0.787024,0.735465,0.36279);
		rgb(731pt)=(0.789092,0.735286,0.361948);
		rgb(732pt)=(0.791159,0.735107,0.361105);
		rgb(733pt)=(0.793227,0.734929,0.360263);
		rgb(734pt)=(0.795295,0.73475,0.359421);
		rgb(735pt)=(0.797362,0.734571,0.358578);
		rgb(736pt)=(0.79943,0.734392,0.357736);
		rgb(737pt)=(0.801485,0.734214,0.356881);
		rgb(738pt)=(0.803527,0.734035,0.356014);
		rgb(739pt)=(0.805569,0.733856,0.355146);
		rgb(740pt)=(0.807611,0.733678,0.354278);
		rgb(741pt)=(0.809668,0.733499,0.353424);
		rgb(742pt)=(0.811735,0.73332,0.352582);
		rgb(743pt)=(0.813803,0.733142,0.35174);
		rgb(744pt)=(0.81587,0.732963,0.350897);
		rgb(745pt)=(0.817921,0.732784,0.350038);
		rgb(746pt)=(0.819963,0.732606,0.349171);
		rgb(747pt)=(0.822005,0.732427,0.348303);
		rgb(748pt)=(0.824047,0.732248,0.347435);
		rgb(749pt)=(0.826071,0.73207,0.346567);
		rgb(750pt)=(0.828087,0.731891,0.345699);
		rgb(751pt)=(0.830104,0.731712,0.344831);
		rgb(752pt)=(0.83212,0.731534,0.343963);
		rgb(753pt)=(0.834158,0.731355,0.343095);
		rgb(754pt)=(0.8362,0.731176,0.342228);
		rgb(755pt)=(0.838242,0.730998,0.34136);
		rgb(756pt)=(0.840284,0.730819,0.340492);
		rgb(757pt)=(0.842303,0.73064,0.339624);
		rgb(758pt)=(0.84432,0.730462,0.338756);
		rgb(759pt)=(0.846336,0.730283,0.337888);
		rgb(760pt)=(0.848353,0.730104,0.33702);
		rgb(761pt)=(0.850369,0.729926,0.336153);
		rgb(762pt)=(0.852386,0.729747,0.335285);
		rgb(763pt)=(0.854402,0.729568,0.334417);
		rgb(764pt)=(0.856419,0.729391,0.333546);
		rgb(765pt)=(0.858435,0.729238,0.332627);
		rgb(766pt)=(0.860452,0.729085,0.331708);
		rgb(767pt)=(0.862468,0.728932,0.330789);
		rgb(768pt)=(0.864481,0.728778,0.329874);
		rgb(769pt)=(0.866472,0.728625,0.32898);
		rgb(770pt)=(0.868463,0.728472,0.328087);
		rgb(771pt)=(0.870454,0.728319,0.327194);
		rgb(772pt)=(0.872445,0.728166,0.326295);
		rgb(773pt)=(0.874436,0.728013,0.325376);
		rgb(774pt)=(0.876427,0.727859,0.324457);
		rgb(775pt)=(0.878418,0.727706,0.323538);
		rgb(776pt)=(0.880417,0.727561,0.322619);
		rgb(777pt)=(0.882433,0.727433,0.3217);
		rgb(778pt)=(0.88445,0.727306,0.320781);
		rgb(779pt)=(0.886466,0.727178,0.319862);
		rgb(780pt)=(0.888463,0.72705,0.318933);
		rgb(781pt)=(0.890429,0.726923,0.317989);
		rgb(782pt)=(0.892394,0.726795,0.317044);
		rgb(783pt)=(0.894359,0.726668,0.3161);
		rgb(784pt)=(0.896337,0.726552,0.315132);
		rgb(785pt)=(0.898328,0.72645,0.314136);
		rgb(786pt)=(0.900319,0.726348,0.313141);
		rgb(787pt)=(0.90231,0.726246,0.312145);
		rgb(788pt)=(0.904301,0.726158,0.31115);
		rgb(789pt)=(0.906292,0.726081,0.310154);
		rgb(790pt)=(0.908283,0.726005,0.309159);
		rgb(791pt)=(0.910274,0.725928,0.308163);
		rgb(792pt)=(0.912249,0.725851,0.307151);
		rgb(793pt)=(0.914214,0.725775,0.30613);
		rgb(794pt)=(0.91618,0.725698,0.305109);
		rgb(795pt)=(0.918145,0.725622,0.304088);
		rgb(796pt)=(0.920111,0.7256,0.303031);
		rgb(797pt)=(0.922076,0.7256,0.301959);
		rgb(798pt)=(0.924041,0.7256,0.300886);
		rgb(799pt)=(0.926007,0.7256,0.299814);
		rgb(800pt)=(0.927972,0.7256,0.298722);
		rgb(801pt)=(0.929938,0.7256,0.297624);
		rgb(802pt)=(0.931903,0.7256,0.296527);
		rgb(803pt)=(0.933869,0.7256,0.295429);
		rgb(804pt)=(0.935812,0.725668,0.294264);
		rgb(805pt)=(0.937752,0.725744,0.29309);
		rgb(806pt)=(0.939692,0.725821,0.291916);
		rgb(807pt)=(0.941632,0.725897,0.290741);
		rgb(808pt)=(0.943571,0.726023,0.289518);
		rgb(809pt)=(0.945511,0.726151,0.288293);
		rgb(810pt)=(0.947451,0.726278,0.287068);
		rgb(811pt)=(0.949389,0.726411,0.285839);
		rgb(812pt)=(0.951278,0.726641,0.284537);
		rgb(813pt)=(0.953167,0.72687,0.283235);
		rgb(814pt)=(0.955056,0.7271,0.281933);
		rgb(815pt)=(0.956938,0.72734,0.280622);
		rgb(816pt)=(0.958776,0.727646,0.279243);
		rgb(817pt)=(0.960614,0.727952,0.277865);
		rgb(818pt)=(0.962451,0.728259,0.276486);
		rgb(819pt)=(0.964273,0.728597,0.275086);
		rgb(820pt)=(0.966034,0.729057,0.273606);
		rgb(821pt)=(0.967795,0.729516,0.272126);
		rgb(822pt)=(0.969557,0.729976,0.270645);
		rgb(823pt)=(0.971288,0.730473,0.269135);
		rgb(824pt)=(0.972947,0.73106,0.267552);
		rgb(825pt)=(0.974606,0.731647,0.265969);
		rgb(826pt)=(0.976265,0.732234,0.264387);
		rgb(827pt)=(0.977857,0.732879,0.262785);
		rgb(828pt)=(0.979338,0.733619,0.261151);
		rgb(829pt)=(0.980818,0.734359,0.259518);
		rgb(830pt)=(0.982299,0.735099,0.257884);
		rgb(831pt)=(0.983697,0.73591,0.256227);
		rgb(832pt)=(0.984999,0.736803,0.254542);
		rgb(833pt)=(0.986301,0.737697,0.252858);
		rgb(834pt)=(0.987603,0.73859,0.251173);
		rgb(835pt)=(0.988753,0.739566,0.249474);
		rgb(836pt)=(0.989774,0.740613,0.247764);
		rgb(837pt)=(0.990795,0.741659,0.246054);
		rgb(838pt)=(0.991816,0.742706,0.244344);
		rgb(839pt)=(0.992677,0.743816,0.242681);
		rgb(840pt)=(0.993443,0.744965,0.241048);
		rgb(841pt)=(0.994209,0.746114,0.239414);
		rgb(842pt)=(0.994975,0.747262,0.23778);
		rgb(843pt)=(0.995578,0.748465,0.236165);
		rgb(844pt)=(0.996114,0.74969,0.234557);
		rgb(845pt)=(0.99665,0.750915,0.232949);
		rgb(846pt)=(0.997186,0.752141,0.231341);
		rgb(847pt)=(0.997562,0.753386,0.229813);
		rgb(848pt)=(0.997893,0.754637,0.228307);
		rgb(849pt)=(0.998225,0.755887,0.226801);
		rgb(850pt)=(0.998557,0.757138,0.225295);
		rgb(851pt)=(0.998711,0.758433,0.223856);
		rgb(852pt)=(0.998839,0.759735,0.222426);
		rgb(853pt)=(0.998966,0.761037,0.220997);
		rgb(854pt)=(0.999094,0.762339,0.219567);
		rgb(855pt)=(0.999076,0.763641,0.218186);
		rgb(856pt)=(0.99905,0.764942,0.216808);
		rgb(857pt)=(0.999025,0.766244,0.21543);
		rgb(858pt)=(0.998995,0.767546,0.214054);
		rgb(859pt)=(0.998868,0.768848,0.212752);
		rgb(860pt)=(0.99874,0.77015,0.21145);
		rgb(861pt)=(0.998613,0.771451,0.210149);
		rgb(862pt)=(0.998473,0.772756,0.208856);
		rgb(863pt)=(0.998243,0.774083,0.207631);
		rgb(864pt)=(0.998014,0.775411,0.206405);
		rgb(865pt)=(0.997784,0.776738,0.20518);
		rgb(866pt)=(0.997539,0.77806,0.20396);
		rgb(867pt)=(0.997232,0.779362,0.20276);
		rgb(868pt)=(0.996926,0.780664,0.201561);
		rgb(869pt)=(0.99662,0.781966,0.200361);
		rgb(870pt)=(0.996299,0.783268,0.199168);
		rgb(871pt)=(0.995942,0.784569,0.197994);
		rgb(872pt)=(0.995584,0.785871,0.19682);
		rgb(873pt)=(0.995227,0.787173,0.195646);
		rgb(874pt)=(0.994842,0.788475,0.19449);
		rgb(875pt)=(0.994408,0.789777,0.193367);
		rgb(876pt)=(0.993974,0.791078,0.192244);
		rgb(877pt)=(0.99354,0.79238,0.191121);
		rgb(878pt)=(0.993083,0.793671,0.190021);
		rgb(879pt)=(0.992598,0.794947,0.188949);
		rgb(880pt)=(0.992113,0.796223,0.187877);
		rgb(881pt)=(0.991628,0.797499,0.186805);
		rgb(882pt)=(0.99113,0.798789,0.185732);
		rgb(883pt)=(0.990619,0.800091,0.18466);
		rgb(884pt)=(0.990109,0.801393,0.183588);
		rgb(885pt)=(0.989598,0.802695,0.182516);
		rgb(886pt)=(0.989072,0.803996,0.18146);
		rgb(887pt)=(0.988536,0.805298,0.180413);
		rgb(888pt)=(0.988,0.8066,0.179367);
		rgb(889pt)=(0.987464,0.807902,0.17832);
		rgb(890pt)=(0.98691,0.809186,0.177291);
		rgb(891pt)=(0.986349,0.810462,0.17627);
		rgb(892pt)=(0.985787,0.811738,0.175249);
		rgb(893pt)=(0.985226,0.813015,0.174228);
		rgb(894pt)=(0.984644,0.814311,0.173207);
		rgb(895pt)=(0.984057,0.815613,0.172186);
		rgb(896pt)=(0.98347,0.816914,0.171165);
		rgb(897pt)=(0.982883,0.818216,0.170144);
		rgb(898pt)=(0.982296,0.819518,0.169145);
		rgb(899pt)=(0.981709,0.82082,0.16815);
		rgb(900pt)=(0.981122,0.822122,0.167154);
		rgb(901pt)=(0.980535,0.823423,0.166159);
		rgb(902pt)=(0.979947,0.824725,0.165163);
		rgb(903pt)=(0.97936,0.826027,0.164168);
		rgb(904pt)=(0.978773,0.827329,0.163172);
		rgb(905pt)=(0.978186,0.828631,0.162177);
		rgb(906pt)=(0.977599,0.829932,0.161181);
		rgb(907pt)=(0.977012,0.831234,0.160186);
		rgb(908pt)=(0.976425,0.832536,0.15919);
		rgb(909pt)=(0.975838,0.833841,0.158195);
		rgb(910pt)=(0.975251,0.835168,0.157199);
		rgb(911pt)=(0.974664,0.836495,0.156204);
		rgb(912pt)=(0.974077,0.837823,0.155208);
		rgb(913pt)=(0.973489,0.83915,0.154213);
		rgb(914pt)=(0.972902,0.840477,0.153217);
		rgb(915pt)=(0.972315,0.841805,0.152222);
		rgb(916pt)=(0.971728,0.843132,0.151226);
		rgb(917pt)=(0.971155,0.844466,0.150224);
		rgb(918pt)=(0.970619,0.845819,0.149203);
		rgb(919pt)=(0.970083,0.847172,0.148182);
		rgb(920pt)=(0.969547,0.848525,0.147161);
		rgb(921pt)=(0.96902,0.849886,0.14614);
		rgb(922pt)=(0.968509,0.851265,0.145119);
		rgb(923pt)=(0.967999,0.852643,0.144098);
		rgb(924pt)=(0.967488,0.854022,0.143077);
		rgb(925pt)=(0.967,0.855411,0.142056);
		rgb(926pt)=(0.966541,0.856815,0.141035);
		rgb(927pt)=(0.966081,0.858219,0.140014);
		rgb(928pt)=(0.965622,0.859623,0.138992);
		rgb(929pt)=(0.965189,0.86104,0.137945);
		rgb(930pt)=(0.96478,0.862469,0.136873);
		rgb(931pt)=(0.964372,0.863899,0.135801);
		rgb(932pt)=(0.963963,0.865328,0.134729);
		rgb(933pt)=(0.96357,0.866773,0.133657);
		rgb(934pt)=(0.963187,0.868228,0.132585);
		rgb(935pt)=(0.962805,0.869683,0.131513);
		rgb(936pt)=(0.962422,0.871138,0.130441);
		rgb(937pt)=(0.962091,0.87261,0.129351);
		rgb(938pt)=(0.961785,0.874091,0.128253);
		rgb(939pt)=(0.961478,0.875571,0.127156);
		rgb(940pt)=(0.961172,0.877052,0.126058);
		rgb(941pt)=(0.960885,0.878571,0.124961);
		rgb(942pt)=(0.960605,0.880103,0.123863);
		rgb(943pt)=(0.960324,0.881634,0.122765);
		rgb(944pt)=(0.960043,0.883166,0.121668);
		rgb(945pt)=(0.959849,0.884719,0.120549);
		rgb(946pt)=(0.95967,0.886276,0.119426);
		rgb(947pt)=(0.959491,0.887833,0.118302);
		rgb(948pt)=(0.959313,0.88939,0.117179);
		rgb(949pt)=(0.959181,0.890995,0.116032);
		rgb(950pt)=(0.959054,0.892603,0.114884);
		rgb(951pt)=(0.958926,0.894211,0.113735);
		rgb(952pt)=(0.958799,0.895819,0.112587);
		rgb(953pt)=(0.958748,0.897453,0.111464);
		rgb(954pt)=(0.958697,0.899086,0.110341);
		rgb(955pt)=(0.958646,0.90072,0.109217);
		rgb(956pt)=(0.958602,0.902359,0.108089);
		rgb(957pt)=(0.958628,0.904043,0.106915);
		rgb(958pt)=(0.958653,0.905728,0.105741);
		rgb(959pt)=(0.958679,0.907413,0.104567);
		rgb(960pt)=(0.958718,0.909102,0.103393);
		rgb(961pt)=(0.95882,0.910812,0.102219);
		rgb(962pt)=(0.958922,0.912522,0.101044);
		rgb(963pt)=(0.959024,0.914232,0.0998703);
		rgb(964pt)=(0.959153,0.915962,0.0986961);
		rgb(965pt)=(0.959357,0.917749,0.0975219);
		rgb(966pt)=(0.959561,0.919536,0.0963477);
		rgb(967pt)=(0.959765,0.921323,0.0951736);
		rgb(968pt)=(0.959996,0.923118,0.0939907);
		rgb(969pt)=(0.960277,0.924931,0.092791);
		rgb(970pt)=(0.960557,0.926743,0.0915913);
		rgb(971pt)=(0.960838,0.928555,0.0903916);
		rgb(972pt)=(0.961151,0.930378,0.0891919);
		rgb(973pt)=(0.961509,0.932216,0.0879922);
		rgb(974pt)=(0.961866,0.934054,0.0867925);
		rgb(975pt)=(0.962223,0.935892,0.0855928);
		rgb(976pt)=(0.96262,0.937768,0.0843802);
		rgb(977pt)=(0.963053,0.939683,0.083155);
		rgb(978pt)=(0.963487,0.941597,0.0819297);
		rgb(979pt)=(0.963921,0.943512,0.0807045);
		rgb(980pt)=(0.964415,0.945426,0.0794643);
		rgb(981pt)=(0.964951,0.947341,0.0782135);
		rgb(982pt)=(0.965487,0.949255,0.0769628);
		rgb(983pt)=(0.966023,0.951169,0.075712);
		rgb(984pt)=(0.966594,0.953118,0.0744441);
		rgb(985pt)=(0.967181,0.955083,0.0731679);
		rgb(986pt)=(0.967768,0.957049,0.0718916);
		rgb(987pt)=(0.968355,0.959014,0.0706153);
		rgb(988pt)=(0.96898,0.960999,0.0693006);
		rgb(989pt)=(0.969619,0.96299,0.0679733);
		rgb(990pt)=(0.970257,0.964981,0.0666459);
		rgb(991pt)=(0.970895,0.966972,0.0653186);
		rgb(992pt)=(0.971554,0.968984,0.0639486);
		rgb(993pt)=(0.972218,0.971001,0.0625703);
		rgb(994pt)=(0.972882,0.973017,0.0611919);
		rgb(995pt)=(0.973545,0.975034,0.0598135);
		rgb(996pt)=(0.974232,0.97705,0.058318);
		rgb(997pt)=(0.974922,0.979067,0.056812);
		rgb(998pt)=(0.975611,0.981083,0.055306);
		rgb(999pt)=(0.9763,0.9831,0.0538)
	}
}

\usetikzlibrary{arrows,shapes,positioning}
\usetikzlibrary{decorations.markings}
\tikzstyle arrowstyle=[scale=1]
\tikzstyle directed=[postaction={decorate,decoration={markings,
		mark=at position .65 with {\arrow[arrowstyle]{stealth}}}}]
\tikzstyle reverse directed=[postaction={decorate,decoration={markings,
		mark=at position .65 with {\arrowreversed[arrowstyle]{stealth};}}}]

\newcommand{\secref}[1]{Section~\ref{#1}}

\newcommand{\rmref}[1]{Remark~\ref{#1}}
\newcommand{\figref}[1]{Figure~\ref{#1}}
\newcommand{\tabref}[1]{Table~\ref{#1}}
\newcommand{\algoref}[1]{Algorithm~\ref{#1}}

\usepackage{mathtools}
\DeclarePairedDelimiter{\ceil}{\lceil}{\rceil}

\newcommand{\R}{\mathbb{R}}

\newcommand{\Domain}{\ensuremath{X}} 
\newcommand{\timevar}{\ensuremath{t}} 

\newcommand{\Lp}[1]{\ensuremath{L_{#1}}}

\newcommand{\flux}{\ensuremath{f}}

\newcommand{\x}{\ensuremath{x}}

\newcommand{\dt}{\partial_\timevar}

\newcommand{\quadxpoint}[1]{\hat{\x}_{#1}}
\newcommand{\quadxweight}[1]{\hat{w}_{#1}}
\newcommand{\quadRpoint}[1]{\hat{\uncertainty}_{#1}}
\newcommand{\quadRweight}[1]{\hat{\omega}_{#1}}

\newcommand{\ncells}{\ensuremath{N_x}}

\newcommand{\timeind}{\ensuremath{n}}

\newcommand{\cellind}{\ensuremath{i}}
\newcommand{\cellindR}{\ensuremath{j}}

\newcommand{\cell}[1]{\ensuremath{C_{#1}}}
\newcommand{\cellR}[1]{\ensuremath{D_{#1}}}

\newcommand{\interface}[1]{\ensuremath{\x_{#1}}}

\newcommand{\polybasis}[1]{\ensuremath{\varphi^{#1}}}

\newcommand{\spatialorder}{\ensuremath{{K_\Domain}}}

\newcommand{\numericalFlux}{\ensuremath{\widehat{\flux}}}
\newcommand{\viscosityconstant}{\ensuremath{c}}

\newcommand{\uncertainty}{\ensuremath{\xi}}
\newcommand{\xiPDF}{\ensuremath{f_\Xi}}
\newcommand{\xiBasisPoly}[1]{\ensuremath{\phi^{#1}}}
\newcommand{\SGsumIndex}{\ensuremath{k}}
\newcommand{\SGsumIndexvar}{\ensuremath{\tilde{k}}}
\newcommand{\xsumIndex}{\ensuremath{h}}
\newcommand{\xsumIndexvar}{\ensuremath{\tilde{h}}}
\newcommand{\SGeqIndex}{\ensuremath{l}}
\newcommand{\SGapproach}{\ensuremath{\sum_{\SGsumIndex=0}^\SGtruncorder \solution^\SGsumIndex \xiBasisPoly{\SGsumIndex}}}
\newcommand{\xiPDFdxi}{\ensuremath{\xiPDF \mathrm{d}\uncertainty}}

\newcommand{\intRS}{\ensuremath{\int_{\randomSpace}}}

\newcommand{\SGtruncorder}{\ensuremath{{K_\randomSpace}}}

\newcommand{\nbxnodes}{\ensuremath{Q_\Domain}}
\newcommand{\nbRnodes}{\ensuremath{Q_\randomSpace}}

\newcommand{\xiQuadIndex}{\ensuremath{\rho}}
\newcommand{\xQuadIndex}{\ensuremath{q}}

\newcommand{\sampleSpace}{\ensuremath{\Omega}}

\newcommand{\probabilityMeasure}{\ensuremath{\mathcal{P}}}
\newcommand{\randomSpace}{\ensuremath{\Xi}}

\newcommand{\solution}{\ensuremath{u}}

\newcommand{\localsolution}[2]{\ensuremath{u}_{#1,#2}}

\newcommand{\MEDGmoment}[3]{\ensuremath{u^{#1,#2}_{#3}}}

\newcommand{\randomElement}[1]{\ensuremath{D_{#1}}}
\newcommand{\indicatorVar}[2]{\ensuremath{\chi_{#1}(#2)}}

\newcommand{\MEIndex}{\ensuremath{j}}
\newcommand{\MEElements}{\ensuremath{{N_\randomSpace}}}
\newcommand{\localxiBasisPoly}[2]{\ensuremath{\phi^{#1}_{#2}}}

\newcommand{\xdegree}{\ensuremath{K_\Domain}}

\newcommand{\ME}{Multielement}

\newcommand{\SD}{\Omega}
\newcommand{\advec}{a}
\newcommand{\PD}{X}
\newcommand{\RR}{\mathbb{R}}
\newcommand{\pd}[1]{\frac{\partial}{\partial #1}}
\newcommand{\dint}[1]{{\,\text{d}#1}}

\newcommand{\tend}{T}
\newcommand{\intT}{\int_0^\tend}
\newcommand{\udata}{\bu_D}
\newcommand{\costfun}{J}
\newcommand{\redcostfun}{\hat{J}}
\newcommand{\costfunxi}{j}
\newcommand{\uncertaintyR}{\uncertainty^R}
\newcommand{\uncertaintyL}{\uncertainty^L}
\newcommand{\uncertaintyInt}{\boldsymbol{\uncertainty}^I}

\newcommand{\lagrange}{\mathcal{L}}
\newcommand{\WF}{\bF}
\newcommand{\adj}{\bp}
\newcommand{\solvec}{\mathbf{\solution}}
\newcommand{\direction}{\bv_\solution}
\newcommand{\directionhk}{\bv_\solution^{\xsumIndex,\SGsumIndex}}
\newcommand{\directionhkvar}{\bv_\solution^{\xsumIndexvar,\SGsumIndexvar}}
\newcommand{\directionxi}{\bv_{\uncertainty}}
\newcommand{\refcell}{\widehat{\cell{}}}
\newcommand{\refME}{\widehat{\randomElement{}}}
\newcommand{\nbtime}{{N_\tend}}
\newcommand{\xidelta}{\delta}
\newcommand{\Fvar}{\hat{F}}
\newcommand{\controltostate}{\mathcal{G}}

\pgfplotscreateplotcyclelist{color parula}{%
	royalblue!70,every mark/.append style={solid,line width = 0pt,fill=royalblue!60!black},mark=ball\\%
	goldenrod!50!yelloworange,every mark/.append style={solid,fill=goldenrod!30!black},mark=*\\%
	green,every mark/.append style={black,fill=yellowgreen},mark=diamond*\\%
	bluegreen,every mark/.append style={fill=black},mark=triangle*\\%
	royalblue!70,densely dashed,every mark/.append style={solid,line width = 0pt,fill=royalblue!60!black},mark=ball\\%
	goldenrod!50!yelloworange,densely dashed,every mark/.append style={solid,black,fill=goldenrod},mark=diamond*\\%
	yellowgreen,densely dashed,every mark/.append style={solid,fill=yellowgreen!60!royalblue},mark=square*, mark size= 1pt\\%
	bluegreen,densely dashed,mark size = 1.5pt,every mark/.append style={solid,fill=white},mark=otimes*\\
	royalblue,densely dashed,mark size = 1.5pt,every mark/.append style={solid,fill=royalblue!30!black},mark=pentagon*\\
	goldenrod!40!yelloworange,every mark/.append style={solid,fill=goldenrod!30!black},mark=|\\%
}

\pgfplotscreateplotcyclelist{color louisa}{%
	blue,every mark/.append style={solid,line width = 0pt,fill=royalblue!60!black},mark=ball\\%
	red,every mark/.append style={solid,fill=orangered!30!black},mark=*\\%
	green,every mark/.append style={black,fill=yellowgreen},mark=diamond*\\%
	goldenrod!50!yelloworange,every mark/.append style={fill=black},mark=triangle*\\%
	royalblue!70,densely dashed,every mark/.append style={solid,line width = 0pt,fill=royalblue!60!black},mark=ball\\%
	goldenrod!50!yelloworange,densely dashed,every mark/.append style={solid,black,fill=goldenrod},mark=diamond*\\%
	yellowgreen,densely dashed,every mark/.append style={solid,fill=yellowgreen!60!royalblue},mark=square*, mark size= 1pt\\%
	bluegreen,densely dashed,mark size = 1.5pt,every mark/.append style={solid,fill=white},mark=otimes*\\
	royalblue,densely dashed,mark size = 1.5pt,every mark/.append style={solid,fill=royalblue!30!black},mark=pentagon*\\
	goldenrod!40!yelloworange,every mark/.append style={solid,fill=goldenrod!30!black},mark=|\\%
}

\pgfplotscreateplotcyclelist{color fabian}{%
	black,dashed\\%
	red,every mark/.append style={solid},mark=diamond*\\%
	blue,every mark/.append style={solid},mark=*\\%
	goldenrod!50!yelloworange,every mark/.append style={fill=black},mark=triangle*\\%
	royalblue!70,densely dashed,every mark/.append style={solid,line width = 0pt,fill=royalblue!60!black},mark=ball\\%
	goldenrod!50!yelloworange,densely dashed,every mark/.append style={solid,black,fill=goldenrod},mark=diamond*\\%
	yellowgreen,densely dashed,every mark/.append style={solid,fill=yellowgreen!60!royalblue},mark=square*, mark size= 1pt\\%
	bluegreen,densely dashed,mark size = 1.5pt,every mark/.append style={solid,fill=white},mark=otimes*\\
	royalblue,densely dashed,mark size = 1.5pt,every mark/.append style={solid,fill=royalblue!30!black},mark=pentagon*\\
	goldenrod!40!yelloworange,every mark/.append style={solid,fill=goldenrod!30!black},mark=|\\%
}

\usepackage[colorinlistoftodos,obeyFinal]{todonotes}

\usepackage{media9}

\begin{document}

\begin{abstract}
	We study an identification problem which estimates the parameters of the underlying random distribution for uncertain scalar conservation laws. The hyperbolic equations are discretized with the so-called discontinuous stochastic Galerkin method, i.e., using a spatial discontinuous Galerkin scheme and a \ME{} stochastic Galerkin ansatz in the random space. We assume an uncertain flux or uncertain initial conditions and that a data set of an observed solution is given. The uncertainty is assumed to be uniformly distributed on an unknown interval and we focus on identifying the correct endpoints of this interval.
The first-order optimality conditions from the discontinuous stochastic Galerkin discretization are computed on the time-continuous level. Then, we solve the resulting semi-discrete forward and backward schemes with the \refone{Runge-Kutta} method. To illustrate the feasibility of the approach, we apply the method to a stochastic advection and a stochastic equation of Burgers' type. 
The results show that the method is able to identify the distribution parameters of the random variable in the uncertain differential equation even if discontinuities are present.
\end{abstract}
\maketitle

\noindent


\section{Introduction}
Uncertainties play a role in many socio-economic, biological or physical phenomena which can be modelled with the help of partial differential equations (PDE)\cite{HertyTosinViscontiZanella,TosinZanella,CarrilloZanella,AlbiPareschiZanella,GoettlichKnapp,DegondJinZhu,BorziSchulzSchillingsvonWinckel,CarrilloPareschiZanella}. In the recent years many  approaches, e.g., methods based on Bayesian inversion, Monte Carlo algorithms or stochastic Galerkin schemes, were proposed to quantify the uncertainties in order to account for them in predictions and simulations \refone{\cite{StuartSchillings,Stuart,Abgrall2013,Abgrall2007,Gottlieb2008,Meyer2019a,Poette2009,JinPareschi,Bos2018,Yan2019}}. 

In this article we focus on hyperbolic conservation laws having uncertainties in parameters which arise, i.e., due to measurement errors and 
thus have non-deterministic effects on the approximation of the deterministic problem. These uncertainties in the parameters can be modelled by random variables that follow an appropriate distribution type. In general, it is difficult to prescribe the exact distribution parameters if measurement errors are present. Therefore, our goal is to identify these distribution parameters from observed data that forms a solution to the uncertain conservation law. This yields the formulation of an optimization problem, whereas the hyperbolic partial differential equation constraint poses severe difficulties both at the continuous and discrete level since they typically form shocks even for smooth initial data when the flux function is non-linear. 

\refone{Another difficulty is raised by the fact that the typical solution spaces for hyperbolic equations have no Hilbert space structure, therefore standard techniques for Optimal Control with PDE constraints are not applicable. We therefore} pursue a discretization within the spatial and stochastic variable to formulate the identification problem in time-continuous form.
\refone{On the ODE level we follow the approach 'first optimize, then discretize' motivated by the findings of \cite{HinzeRoesch}. This has the advantage that the state and the adjoint problem can be solved with different techniques leading to higher efficiency.}

Discretizing the conservation law within the stochastic domain, we consider Uncertainty Quantification (UQ) methods \refone{\cite{abgrall2017uncertainty,Kolb2018,le2010spectral,pettersson2015polynomial,Smith2014,Kusch2018, gerster2020entropies}} that aim to model the propagation of the uncertainty into the solution of uncertain equations. We distinguish between the so-called non-intrusive and intrusive schemes. The most widely known non-intrusive UQ method is (Multi-Level) Monte Carlo \cite{lye2016multilevel,Heinrich2001,Giles2008,CarrilloZanella}, which is based on statistical sampling methods that can easily be adopted to our problem setting but comes with potentially high costs due to the repeated application of finite volume schemes. \reftwo{Another approach is to employ a discretization in space which leads to a stochastic differential system. Here one could try to apply an parameter identification in the spirit of \cite{spacemapping}.}
Within this article, we concentrate on an intrusive UQ method, namely the stochastic Galerkin (sG) scheme, that involves modifications of the finite volume solver. The method relies on the generalized Polynomial Chaos (gPC) expansion \cite{Xiu2003,Wan2006a,Abgrall2007,Chen2005,Wiener1938}, thus expands the solution in the stochastic variable and projects it on the space spanned by a truncated orthonormal basis.

The biggest challenge of UQ methods for hyperbolic equations lies in the fact, that discontinuities in the physical space propagate into the solution manifold such that the polynomial expansion of discontinuous data yields huge oscillations \cite{barth2013non,Poette2009,Schlachter2017a}. Therefore, the authors of \cite{Wan2006} introduced the so-called \ME{} approach, where the random space is divided into disjoint elements in order to define local gPC approximations. Further developments of this method can be found in \cite{Tryoen2010,Wan2005,Wan2009}.
Similar to this ansatz, we apply a spatial discontinuous Galerkin discretization \cite{Cockburn1989,Cockburn1989a,Cockburn1991}, where we expand the solution in the spatial variable into piecewise polynomials and perform the Galerkin approach in each physical cell. The resulting semi-discrete weak formulation and the discontinuous stochastic Galerkin method is high-order in both, space and the stochasticity \cite{Meyer2019a}.

Optimal control problems for fully discrete hyperbolic conservation laws have been studied for instance in \cite{Ulbrich2001,Steffensen2013,Hajian2019}. Here, the minimization of a given cost functional subject to a deterministic conservation law is considered, for example in order to identify the correct initial conditions under an observed solution. The authors state the first-order optimality conditions for conservation laws discretized in space and time, using first-order finite volume methods as in Lax Friedrichs or Engquist-Osher discretizations. Moreover, the convergence and stability of the adjoint equation has been studied, as well as the differentiability of the control-to-state map and the reduced cost functional, yielding the introduction of the so-called shift derivative \cite{Bressan1997,Ulbrich2003}. The construction of discrete adjoints for hyperbolic partial differential equations is described for example in \cite{Giles2003,Ulbrich2003b}. In this context, the total variation diminishing properties of \refone{Runge-Kutta} schemes have been analyzed \cite{Hajian2019}, especially in terms of the order of the adjoint approximation. Second-order approximations have been formulated using the so-called relaxation method \cite{Banda2012,Jin1995}. Up to the authors' knowledge, a general framework for optimization in hyperbolic problems without discretization has not been established, yet.
Hence, the intention of this article is the extension of optimal control techniques for conservation laws on uncertain equations using the high-order discontinuous sG scheme on the semi-discrete level. 

The paper is structured as follows. In \secref{sec:motivation} we introduce the problem setting and motivate the arising parameter identification problem. We discretize the conservation law within the stochastic and spatial variable  in \secref{sec:DSG}. This yields an ODE system on which the identification problem is based, we derive its optimality condition and adjoint system on the time-continuous level within \secref{sec:FOcond}. In \secref{sec:scheme} we use these first order optimality conditions to formulate an algorithm for the parameter estimation. Finally, we apply this algorithm to the uncertain linear advection and Burgers' equation. The numerical results in \secref{sec:results} show the impact of discretization parameters on the optimization process. 

\section{Motivation} \label{sec:motivation}
We consider stochastic conservation laws of the form
\begin{subequations}\label{eq:spde}
\begin{equation}\label{eq:spde_eq}
\pd{\timevar} \solution(\timevar,\x,\uncertainty) + \pd{\x} \flux\bigl(\solution(\timevar,\x,\uncertainty),\xi\bigr) = 0,\quad \quad \text{for} \ \ \x\in\Domain,\, \timevar\in [0,\,\tend],\, \uncertainty\in\SD,
\end{equation}
with physical domain $\Domain \subset \RR$, time interval $[0,\,\tend]$, stochastic domain $\SD \subset \RR$ and initial conditions given by
\begin{equation}\label{eq:spde_BC}
\solution(0,\x,\uncertainty) = \solution_{(0)}(\x,\uncertainty), \quad \quad \text{for} \ \ \x\in\Domain,\, \xi\in\SD.
\end{equation}
\end{subequations}
Depending on $\PD$, additional boundary conditions have to be prescribed.
The solution $\solution\in\RR$ depends on a one-dimensional random variable \refone{$\uncertainty\in\sampleSpace\subset\R$.} 
\refone{We denote the random space of this uncertainty by $\randomSpace := \uncertainty(\sampleSpace)$ such that $\uncertainty:\sampleSpace\rightarrow\randomSpace$. The probability density function is given by $\xiPDF(\uncertainty) : \randomSpace \rightarrow \R_+$.}

\refone{We denote by $\uncertaintyInt$ a vector containing the vectorized parameters of the probability distribution. 
The task at hand is to obtain these parameters $\uncertaintyInt$, characterizing the distribution of $\uncertainty$ for the data $u_D$ observed at time $\tend$.}
Thus, we construct an optimization problem of the form
	\begin{equation}\label{eq:OPintro}
	\text{min}_{\uncertaintyInt} \,\costfun(\solution,\uncertaintyInt) = \frac{1}{2}\big\| \hat{u}(T,\uncertaintyInt) - u_D\big\|_{\refone{2}}^2 + \costfunxi(\uncertaintyInt) 
	\qquad \text{subject to } \eqref{eq:spde} \text{ on }\Domain,
	\end{equation}
where $\costfun(\solution,\uncertaintyInt)$ is the cost functional and \refone{$\|\cdot\|_2$ denotes the $L_2$ norm}. Its first term measures the distance of the given data set $u_D$ to the simulated data $\hat{u}(\tend,\uncertaintyInt)$ that solves \eqref{eq:spde} at $\timevar=\tend$ under the parameters $\uncertaintyInt$. 
We will determine the exact structure of this term in \secref{sec:FOcond}.
The second term $\costfunxi(\uncertaintyInt)$ may be used as problem dependent regularizing term. For example, setting \refone{$\costfunxi(\uncertaintyInt) = \frac{1}{2}\big\|\uncertaintyInt-\uncertaintyInt_p\big\|_2^2$} penalizes the distance to a prior $\uncertaintyInt_p$. 

 Up to the authors' knowledge it is not clear how to solve optimization problems constrained by hyperbolic stochastic PDEs on the continuous level. Therefore, we discretize with a discontinuous stochastic Galerkin approximation before we state the identification problem on the ODE level and derive the corresponding first order optimality system.

\section{Discontinuous stochastic Galerkin scheme}\label{sec:DSG}
In this section, we recall the hyperbolicity-preserving discontinuous stochastic Galerkin method introduced in \cite{Meyer2019a}. Our aim is to construct a numerical scheme that is high-order in space, time and the uncertainty. The discretization of the uncertainty is obtained by a stochastic Galerkin scheme with \ME{} ansatz \cite{Wan2006}. For the spatial discretization of the resulting stochastic Galerkin system, we consider a discontinuous Galerkin scheme \cite{Cockburn1989,Cockburn1989a,Cockburn1991}. 

\subsection{\ME{} stochastic Galerkin}
We seek for an approximate solution by a 
generalized Polynomial Chaos (gPC) expansion \cite{Gottlieb2008} with $\SGtruncorder +1$ terms, where $\SGtruncorder \in \mathbb{N}$. Indeed, we have 
\begin{equation}\label{gPCApprox}
\solution(\timevar,\x,\uncertainty) \approx  \sum\limits_{\SGsumIndex=0}^{\SGtruncorder} \solution^\SGsumIndex(\timevar,\x) \, \xiBasisPoly{\SGsumIndex}(\uncertainty) ,
\end{equation}
where the polynomials $\xiBasisPoly{\SGsumIndex}$ of degree $\SGsumIndex$ are supposed to satisfy the orthogonality relation
\begin{equation}\label{eq:orthogonality}
\int\limits_{\randomSpace} \xiBasisPoly{\SGsumIndex}(\uncertainty) \, \xiBasisPoly{\SGeqIndex}(\uncertainty) \, \xiPDF(\uncertainty) \dint{\uncertainty} =
\begin{cases}
1, &\quad\text{if }\SGsumIndex=\SGeqIndex\\
0,&\quad\text{else}
\end{cases} \quad \quad \forall \ \SGsumIndex,\SGeqIndex \in \{0, \dots, \SGtruncorder\}.
\end{equation}
Inserting \eqref{gPCApprox} into \eqref{eq:spde} and applying a Galerkin projection in the stochastic space leads to the so called stochastic Galerkin system
\begin{equation}\label{SGsystem}\pd{t} \solution^\SGeqIndex +  \pd{\x} \intRS \flux\left(\SGapproach,\,\uncertainty\right) \!\xiBasisPoly{\SGeqIndex}\xiPDFdxi = 0, \qquad \SGeqIndex = 0,\ldots,\SGtruncorder.\end{equation}
For discontinuous solutions, the gPC approach may converge slowly or even fail to converge, cf. \cite{Wan2005, Poette2009}.  As presented in \cite{Wan2006,Meyer2019a}, we therefore apply the \ME{} approach, where $\randomSpace$ is divided into disjoint elements with local gPC approximations of  \refone{\eqref{gPCApprox}}.

We assume that $\randomSpace=(\uncertaintyL,\uncertaintyR)$ and define a decomposition of $\randomSpace$ into $\MEElements$ \ME{}s $\randomElement{\MEIndex} = (\uncertainty_{\MEIndex-\frac12}, \uncertainty_{\MEIndex+\frac12})$ of width $\Delta\uncertainty=\frac{\uncertaintyL-\uncertaintyR}{\MEElements}$.  Moreover, we  introduce an indicator variable $\chi_{\MEIndex}: \Omega \to  \{0,1\}$ on every random element
\refone{\begin{equation}\label{def:indicator}
\indicatorVar{\MEIndex}{\uncertainty} :=
\begin{cases}
1&\text{if } \uncertainty \in \randomElement{\MEIndex}, \\
0 &\text{else, }
\end{cases}
\end{equation}}
for $\MEIndex=1,\ldots,\MEElements $ and $\omega\in\sampleSpace$. 
If we let $\{\localxiBasisPoly{\SGsumIndex}{\MEIndex} \}_{\SGsumIndex=0}^\infty$ be  orthonormal polynomials with respect to a conditional probability density function on the \ME{} $\randomElement{\MEIndex}$ as in \cite{Meyer2019a}, 
the global approximation \eqref{gPCApprox} can be written as
\begin{equation}\label{def:globalSG}
\solution(t,x,\uncertainty) = \sum_{\MEIndex=1}^{\MEElements}  \solution_\MEIndex(t,x,\uncertainty) \indicatorVar{\MEIndex}{\uncertainty} \approx \sum_{\MEIndex=1}^{\MEElements}  \sum_{\SGsumIndex=0}^{\SGtruncorder}  \solution^{\SGsumIndex}_{\MEIndex}(\timevar,\x) \localxiBasisPoly{\SGsumIndex}{\MEIndex}(\uncertainty) \indicatorVar{\MEIndex}{\uncertainty}.
\end{equation}
As $\MEElements,\SGtruncorder\to\infty$, the local approximation converges to the global solution in $\Lp{2}(\sampleSpace)$, cf. \cite{Alpert1993}.
\begin{remark}\label{rem:paramDistr}
	We emphasize that this approach is based on the knowledge of the distribution of $\xi$. Most applications assume $\xi$ to be uniformly distributed in a given interval $[\xi^L, \xi^R]$. In the following we propose an ansatz to identify the true values of $\xi^L, \xi^R.$ In the case of a normal distribution, we would estimate the expected value and variance of the random variable.
\end{remark}

\subsection{Discontinuous Galerkin}
Similar to the previous subsection, where we subdivided the random space  $\randomSpace$ into \ME{s} $\cellR{\cellindR}$, $\cellindR =1,\ldots,\MEElements $, we now divide the spatial domain $\Domain \subset\R$ into a uniform rectangular mesh with cells $\cell{\cellind}= [\x_{i-\frac{1}{2}},\x_{i+\frac{1}{2}}]$  of width
$\Delta x:=(\x_{i+\frac{1}{2}}- \x_{i-\frac{1}{2}}),$ where $i=1,\dots, N_x.$ 

The local approximation in the bases of the spaces of piecewise polynomials then reads
\begin{align}\label{localMEDGDGApproximation}
\refone{\localsolution{\cellind}{\cellindR}(\timevar,\x,\uncertainty) 
= \sum_{\xsumIndex=0}^{\xdegree}\sum_{\SGsumIndex=0}^\SGtruncorder \MEDGmoment{\xsumIndex}{\SGsumIndex}{\cellind,\MEIndex}(\timevar)\polybasis{\xsumIndex}_\cellind(\x)\localxiBasisPoly{\SGsumIndex}{\MEIndex}(\uncertainty), }
\end{align}
for $\x\in \cell{\cellind}$ and $\uncertainty\in\cellR{\cellindR}$.
Here, $\{\localxiBasisPoly{\SGsumIndex}{\MEIndex}\}_{\SGsumIndex=0,\ldots,\SGtruncorder}$ are the local basis polynomials on the random element $\cellR{\cellindR}$ and $\{\polybasis{\xsumIndex}_{\cellind}\}_{\xsumIndex=0,\ldots,\xdegree }$ are the basis polynomials on the physical cell $\cell{\cellind}$.
\begin{remark}\label{rem:refCell}
	In our numerical computations we use a reference cell $\refcell=[-\frac12,\,\frac12]$ and reference \ME{} $\refME=[-\frac12,\,\frac12]$, and define our basis polynomials on these reference elements such that the local approximation reads
	$$\localsolution{\cellind}{\cellindR}(\timevar,\x,\uncertainty)  = \sum_{\xsumIndex=0}^{\xdegree}\sum_{\SGsumIndex=0}^\SGtruncorder \MEDGmoment{\xsumIndex}{\SGsumIndex}{\cellind,\MEIndex}(\timevar)\,\polybasis{\xsumIndex}\Big(\frac{\x-\x_{\cellind}}{\Delta\x}\Big)\,\phi^{\SGsumIndex}\Big(\frac{\uncertainty-\uncertainty_{\MEIndex}}{\Delta\uncertainty}\Big).$$
	\refone{Note that for a uniform distribution we have $\xiPDF=\frac{1}{\Delta\uncertainty}$, thus}
	$$\frac{1}{\Delta\x}\int_{\cell{\cellind}}\!\int_{\randomElement{\MEIndex}}\!\polybasis{\xsumIndex}\Big(\frac{\x-\x_{\cellind}}{\Delta\x}\Big)\,\phi^{\SGsumIndex}\Big(\frac{\uncertainty-\uncertainty_{\MEIndex}}{\Delta\uncertainty}\Big)\,\xiPDFdxi\,\mathrm{d}\x 
	=\int_{\refcell}\!\int_{\refME}\!\polybasis{\xsumIndex}(\hat{\x})\phi^{\SGsumIndex}(\hat{\uncertainty})\,\mathrm{d}\hat{\uncertainty}\,\mathrm{d}\hat{\x}.$$
	For the ease of notation, we write $\polybasis{\xsumIndex}_{\cellind}(\x)$ instead of $\polybasis{\xsumIndex}\big(\frac{\x-\x_{\cellind}}{\Delta\x}\big)$ and  $\localxiBasisPoly{\SGsumIndex}{\MEIndex}(\uncertainty)$ for $\phi^{\SGsumIndex}\big(\frac{\uncertainty-\uncertainty_{\MEIndex}}{\Delta\uncertainty}\big)$, respectively. 
\end{remark}

We consider a weak formulation for the solution $\solution$ in the control volume $\cell{\cellind}\times\randomElement{\MEIndex}$ and test \eqref{eq:spde} by a smooth function $\rho(\x,\uncertainty)$ with supp$(\rho)\subseteq \cell{\cellind}\times\randomElement{\MEIndex}$. After formal integration by parts we obtain 
\begin{align}
\label{eq:weakformulation}
\frac{\partial}{\dt} \int_{\cell{\cellind}}\! \int_{\randomElement{\MEIndex}}\! \solution \, \rho \,\xiPDFdxi\mathrm{d}\x 
- \int_{\cell{\cellind}}\! \int_{\randomElement{\MEIndex}}\! \flux(\solution) \, \partial_\x\rho \,\xiPDFdxi\mathrm{d}\x
+ \int_{\randomElement{\MEIndex}}\! \big[\flux(\solution)\rho \big]_{\cellind-\frac12}^{\cellind+\frac12} \xiPDFdxi= 0.
\end{align}
\begin{remark}\label{rem:existunique}
	Solutions to conservation laws \eqref{eq:spde} might develop shocks and must be interpreted in the weak sense as given in \eqref{eq:weakformulation}. The existence and uniqueness of a so called random entropy solution is proved in \cite{Mishra2014}. They correspond $\probabilityMeasure$ almost everywhere to entropy solutions in the deterministic case which is discussed in Kru\v{z}kov's theorem \cite{Kruzov1970}. For more information on this topic, we refer to \cite{Mishra2016,Meyer2019}.
\end{remark}
Since the solution is discontinuous across the interfaces $\x_{i\pm\frac{1}{2}}$, the flux at these points has to be replaced with a numerical flux function $\numericalFlux$, approximately solving the corresponding Riemann problem at the interface. For example, one can choose the global Lax-Friedrichs flux, as we do for the numerical examples. In this case we have
\begin{align}
\label{eq:globalLF}
\numericalFlux(\solution^-, \solution^+) = \dfrac{1}{2} \left( \flux(\solution^-) + \flux(\solution^+) - \viscosityconstant ( \solution^+ - 
\solution^-) \right).
\end{align} 
The numerical viscosity constant $\viscosityconstant$ is taken as the global estimate of the  absolute value of the largest eigenvalue of the Jacobian
$\frac{\partial \flux(\solution)}{\partial \solution}$.
The values 
\begin{equation}\label{eq:limitinterface}
\refone{
\solution^-(\interface{\cellind+\frac12},\uncertainty) := \lim\limits_{\x\uparrow\interface{\cellind+\frac12}} \solution(\x,\uncertainty), \qquad 
\solution^+(\interface{\cellind+\frac12},\uncertainty) := \lim\limits_{\x\downarrow\interface{\cellind+\frac12}} \solution(\x,\uncertainty)}
\end{equation}
denote the left and right limits of the piecewise polynomial solution at the interface $\interface{\cellind+\frac12}$, whereas $\interface{\cellind+\frac12}^\mp$ describe the left and right limits to the this interface, respectively.
\refone{The Lax Friedrichs scheme is stable under the CFL condition
$$\Delta\timevar \leq \frac{\Delta\x}{\viscosityconstant}.$$}

Choosing the test function $\rho = \polybasis{\xsumIndex}_{\cellind}\,\localxiBasisPoly{\SGsumIndex}{\MEIndex}$ together with \eqref{localMEDGDGApproximation} and the orthogonality condition of the basis polynomials, yields the following discontinuous stochastic Galerkin formulation
\begin{subequations}\label{eq:DG}
\begin{align} 
0 = \frac{d}{dt} \MEDGmoment{\xsumIndex}{\SGsumIndex}{\cellind,\MEIndex}
\frac{1}{\Delta\x}\int_{\cell{\cellind}}\!\polybasis{\xsumIndex}_{\cellind}&\polybasis{\xsumIndex}_{\cellind}\,\mathrm{d}\x 
\int_{\randomElement{\MEIndex}}\!\localxiBasisPoly{\SGsumIndex}{\MEIndex}\localxiBasisPoly{\SGsumIndex}{\MEIndex}\,\xiPDFdxi 
\\
&- \frac{1}{\Delta\x} \int_{\cell{\cellind}}\!\int_{\randomElement{\MEIndex}}\! \flux
\big(\localsolution{\cellind}{\cellindR}(\timevar,\x,\uncertainty)\big)
\localxiBasisPoly{\SGsumIndex}{\MEIndex} \,\partial_x\polybasis{\xsumIndex}_{\cellind}\, \xiPDFdxi \mathrm{d}\x \\
& + \frac{1}{\Delta\x}\int_{\randomElement{\MEIndex}}\!\numericalFlux\big(\localsolution{\cellind}{\cellindR}(\timevar,\x_{i+\frac12},\uncertainty), \localsolution{\cellind+1}{\cellindR}(\timevar,\x_{i+\frac12},\uncertainty) \big) \polybasis{\xsumIndex}_{\cellind}(\x_{i+\frac12})\localxiBasisPoly{\SGsumIndex}{\MEIndex} \xiPDFdxi \\
& - \frac{1}{\Delta\x}\int_{\randomElement{\MEIndex}}\! 
\numericalFlux\big(\localsolution{\cellind-1}{\cellindR}(\timevar,\x_{i-\frac12},\uncertainty), \localsolution{\cellind}{\cellindR}(\timevar,\x_{i-\frac12},\uncertainty) \big) \polybasis{\xsumIndex}_{\cellind}(\x_{i-\frac12})\localxiBasisPoly{\SGsumIndex}{\MEIndex} \xiPDFdxi,
\end{align}  
\end{subequations}
for all cells  $\cell{\cellind}\times\randomElement{\MEIndex}$ and $\xsumIndex=0,\ldots,\xdegree$, $\SGsumIndex=0,\ldots,\SGtruncorder$.
The semi-discrete system \eqref{eq:DG} is solved using a $\xdegree$-th order SSP \refone{Runge-Kutta}
method, see \cite{Gottlieb2005}. Moreover, we apply the TVBM minmod slope limiter from \cite{Cockburn2001} to the discontinuous Galerkin polynomial in each \refone{Runge-Kutta} stage. In the case of hyperbolic systems of equations we propose to apply the hyperbolic slope limiter from \cite{Schlachter2017a}. \refone{For convergence and error analyses of the discontinuous stochastic Galerkin scheme \eqref{eq:DG} we refer to \cite{Meyer2019a}.}

The numerical approximation of the integrals in \eqref{eq:DG} is realized by a Gau\ss{}-Lobatto rule on $\cell{\cellind}$ with $\nbxnodes+1 = \ceil*{\frac{\spatialorder+1}{2}}+1$ points and weights $(\quadxpoint{\xQuadIndex},\quadxweight{\xQuadIndex})$, where $\xQuadIndex=0,\ldots,\nbxnodes$. \refone{The Gau\ss{}-Lobatto quadrature rule includes the endpoints, i.e., cell interfaces which we use to evaluate the values \eqref{eq:limitinterface}.} For an uniformly distributed uncertainty and Legendre basis functions, we apply a \refone{Gau\ss{}-Legendre} quadrature rule on $\cellR{\cellindR}$ with order $\SGtruncorder$, i.e., $\nbRnodes+1$  points and weights $(\quadRpoint{\xiQuadIndex},\quadRweight{\xiQuadIndex})$, $\xiQuadIndex=0,\ldots,\nbRnodes$, where $\nbRnodes = \ceil*{\frac{\SGtruncorder+1}{2}}$.  The quadrature weights are scaled such that
\begin{equation}\label{eq:quad}
\sum_{\xQuadIndex=0}^{\nbxnodes} \sum_{\xiQuadIndex=0}^{\nbRnodes}\quadxweight{\xQuadIndex}\quadRweight{\xiQuadIndex}=1, \qquad
\int_{\randomElement{\MEIndex}}\!\bg\,\xiPDFdxi 
\approx  \sum_{\xiQuadIndex=0}^{\nbRnodes}\bg(\quadRpoint{\xiQuadIndex})\quadRweight{\xiQuadIndex}.\end{equation}

\begin{remark}
	For other distributions we propose to use the corresponding Gauß quadrature based on the orthogonal basis polynomials and weighted by the probability density function. 
\end{remark}

Inserting the global Lax-Friedrichs flux \eqref{eq:globalLF} into the discontinuous stochastic Galerkin scheme \eqref{eq:DG}, we define the right hand side of these equations as
\begin{align}\label{eq:DSGM}
\WF(\solvec^{\xsumIndex,\SGsumIndex},\uncertaintyInt)&:=
M_1^{\xsumIndex,\SGsumIndex} \partial_\timevar\solvec^{\xsumIndex,\SGsumIndex} -  
M_2^{\xsumIndex,\SGsumIndex}(\solvec) + 
M_3^{\xsumIndex,\SGsumIndex} (\solvec),
\end{align}
for all $\xsumIndex=0,\ldots,\xdegree$ and $\SGsumIndex=0,\ldots,\SGtruncorder$, where $\solvec^{\xsumIndex,\SGsumIndex}=\solvec^{\xsumIndex,\SGsumIndex}(\timevar)\in\mathbb{R}^{\ncells\cdot\MEElements}$ contains the vectorized coefficients $\MEDGmoment{\xsumIndex}{\SGsumIndex}{\cellind,\MEIndex}(\timevar)$ for all $\cell{\cellind}\times\randomElement{\MEIndex}$. 
Additionally, we denote by $\solvec$ the vectorized values $\solvec^{\xsumIndex,\SGsumIndex}$ and have $\WF(\solvec^{\xsumIndex,\SGsumIndex},\uncertaintyInt)\in\mathbb{R}^{\ncells\cdot\MEElements}$.
\begin{remark}\label{rm:order}
	The chronological order of the elements in these vectors can be chosen arbitrary but consistent throughout the following notation.
\end{remark}
 Moreover, 
 $M_1^{\xsumIndex,\SGsumIndex}\in\mathbb{R}^{\ncells\cdot\MEElements\times\ncells\cdot\MEElements}$ and $M_2^{\xsumIndex,\SGsumIndex}(\solvec),\, M_3^{\xsumIndex,\SGsumIndex}(\solvec) \in\mathbb{R}^{\ncells\cdot\MEElements}$  are given by
\begin{align*}
\refone{M_1^{\xsumIndex,\SGsumIndex}} &\refone{:= \frac{1}{\Delta\x}\text{diag}\Big(\int_{\cell{\cellind}}\!\polybasis{\xsumIndex}_{\cellind}\polybasis{\xsumIndex}_{\cellind}\,\mathrm{d}\x 
\int_{\randomElement{\MEIndex}}\!\localxiBasisPoly{\SGsumIndex}{\MEIndex}\localxiBasisPoly{\SGsumIndex}{\MEIndex}\,\xiPDFdxi \Big)}
\\
M_2^{\xsumIndex,\SGsumIndex}(\solvec) &:= 
\frac{1}{\Delta\x}\text{vec}\Big(\int_{\cell{\cellind}}\!	\int_{\randomElement{\MEIndex}}\! \flux( \localsolution{\cellind}{\cellindR})\,\partial_\x\polybasis{\xsumIndex}_{\cellind}\,
\localxiBasisPoly{\SGsumIndex}{\MEIndex}\,\xiPDFdxi\mathrm{d}\x  \Big)\,,\\
M_3^{\xsumIndex,\SGsumIndex}(\solvec) &:= \frac{1}{2\Delta\x}\text{vec}\begin{pmatrix*}[l]
\hspace*{1.1cm}\int_{\randomElement{\MEIndex}}\!\big(-\flux\big( \localsolution{\cellind-1}{\cellindR}(\x_{\cellind-\frac12})\big)-\viscosityconstant\, \localsolution{\cellind-1}{\cellindR}(\x_{\cellind-\frac12})\big)\,\polybasis{\xsumIndex}_{\cellind}(\x_{\cellind-\frac12})\localxiBasisPoly{\SGsumIndex}{\MEIndex}\,\xiPDFdxi \,\\[.1cm]
+\int_{\randomElement{\MEIndex}}\!\big(\flux\big( \localsolution{\cellind}{\cellindR}(\x_{\cellind+\frac12})\big)+\viscosityconstant\, \localsolution{\cellind}{\cellindR}(\x_{\cellind+\frac12})\big)\,\polybasis{\xsumIndex}_{\cellind}(\x_{\cellind+\frac12})\localxiBasisPoly{\SGsumIndex}{\MEIndex}\,\xiPDFdxi\\
\hspace*{2.4cm}+\int_{\randomElement{\MEIndex}}\!\big(-\flux\big( \localsolution{\cellind}{\cellindR}(\x_{\cellind-\frac12})\big)+\viscosityconstant\, \localsolution{\cellind}{\cellindR}(\x_{\cellind-\frac12})\big)\,\polybasis{\xsumIndex}_{\cellind}(\x_{\cellind-\frac12})\localxiBasisPoly{\SGsumIndex}{\MEIndex}\,\xiPDFdxi\\
\hspace*{1.3cm}+\int_{\randomElement{\MEIndex}}\!\big(\flux\big( \localsolution{\cellind+1}{\cellindR}(\x_{\cellind+\frac12})\big)-\viscosityconstant\, \localsolution{\cellind+1}{\cellindR}(\x_{\cellind+\frac12})\big)\,\polybasis{\xsumIndex}_{\cellind}(\x_{\cellind+\frac12})\localxiBasisPoly{\SGsumIndex}{\MEIndex}\,\xiPDFdxi\\
\end{pmatrix*},\,
\end{align*}
for $\cellind=1,\ldots,\ncells$ and $\MEIndex=1,\ldots,\MEElements$.
Note that the \ME{}s $\randomElement{\MEIndex}$ and $\xiPDF$ are depending on $\uncertaintyInt$. In the notation of $M_3^{\xsumIndex,\SGsumIndex}(\solvec)$, the first line represents the part that corresponds to cell $C_{i-1}\times D_j$, the second and third line correspond to $C_{i}\times D_j$ etc. According to \eqref{localMEDGDGApproximation}, we have  $$\localsolution{\cellind}{\cellindR}(\x_{\cellind\pm\frac12})
= \sum_{\xsumIndex=0}^{\xdegree}\sum_{\SGsumIndex=0}^\SGtruncorder \MEDGmoment{\xsumIndex}{\SGsumIndex}{\cellind,\MEIndex}\,\polybasis{\xsumIndex}_{\cellind} (\x_{\cellind\pm\frac12})\,\localxiBasisPoly{\SGsumIndex}{\MEIndex}.$$
We are now able to formulate the optimization problem within the discretized conservation law \eqref{eq:DSGM}.

\section{Parameter Identification}\label{sec:FOcond}
In this section, we specify the optimization problem for the parameter identification introduced in \eqref{eq:OPintro} and  derive the corresponding first-order optimality conditions. We focus on uncertainties that are uniformly distributed, in fact, we assume $\uncertainty\sim\mathcal{U}(\SD)$. Other distributions can be dealt with analogously.

\subsection{Optimization Problem}
 According to \rmref{rem:paramDistr}, we estimate the space $\randomSpace=[\uncertaintyL,\uncertaintyR]$, i.e., the parameters $\uncertaintyL$ and $\uncertaintyR$ which we write as $\uncertaintyInt = \begin{pmatrix}\uncertaintyL\\\uncertaintyR\end{pmatrix}$ with prior estimates $\uncertaintyInt_p = \begin{pmatrix}\uncertaintyL_p\\\uncertaintyR_p\end{pmatrix}$, and set the penalizing term in the cost functional to
$$\costfunxi(\uncertaintyInt)= \frac{\xidelta}{2}\big\|\uncertaintyInt-\uncertaintyInt_p\big\|_{\refone{2}}^2,$$
where \refone{we denote by $\|\cdot\|_2$ the $L_2$-norm}. Furthermore, $\xidelta\ll1$ diminishes the impact of $\uncertaintyInt_p$ on the cost functional since $\uncertaintyInt_p$ is only a rough guess of the  true value.

\refone{
	The cost functional then reads}
\begin{equation}\label{eq:costfunt}\costfun(\solvec^{\xsumIndex,\SGsumIndex},\uncertaintyInt) =\frac{1}{2}\big\|\solvec^{\xsumIndex,\SGsumIndex}(\tend) - \udata^{\xsumIndex,\SGsumIndex}\big\|_{\refone{2}}^2 + \frac{\xidelta}{2}\big\|\uncertaintyInt-\uncertaintyInt_p\big\|_{\refone{2}}^2,\end{equation}
where the coefficients $\udata^{\xsumIndex,\SGsumIndex}\in\mathbb{R}^{\ncells\cdot\MEElements}$ describe the observed polynomial moments of the solution for each $\xsumIndex=0,\ldots,\xdegree$ and $\SGsumIndex=0,\ldots,\SGtruncorder$. 
Thus, the polynomial degrees $\xdegree$ and $\SGtruncorder$ should be chosen accordingly to the number of moments of $\udata$ that we have observed.
\refone{\begin{remark}
	For other distribution types of $\uncertainty$, the same cost functional \eqref{eq:costfunt} can be used, whereas $\uncertaintyInt$ contains the new parameters of the distribution such as the expected value and variance  in the case of a normal distribution, cf. \rmref{rem:paramDistr}.
\end{remark}}

Furthermore, we define $$\costfun(\solvec,\uncertaintyInt):= \sum_{\xsumIndex=0}^{\xdegree}\sum_{\SGsumIndex=0}^\SGtruncorder \costfun(\solvec^{\xsumIndex,\SGsumIndex},\uncertaintyInt).$$

For the initial conditions, we introduce the term
$$\WF_{(0)}(\solvec^{\xsumIndex,\SGsumIndex}) := \solvec^{\xsumIndex,\SGsumIndex}(0) - \solvec^{\xsumIndex,\SGsumIndex}_{(0)},$$
which is derived in every cell $\cell{\cellind}\times D_\MEIndex$ for $\cellind = 1,\ldots,\ncells$ and $\MEIndex=1,\ldots,\MEElements$ using
\begin{equation}\label{eq:u0}
\big(\solution^{\xsumIndex,\SGsumIndex}_{(0)}\big)_{(\cellind,\MEIndex)} = \refone{\frac{1}{\int_{\cell{\cellind}}\!\int_{\randomElement{\MEIndex}}\!(\localxiBasisPoly{\SGsumIndex}{\MEIndex})^2(\polybasis{\xsumIndex}_{\cellind})^2\, \xiPDFdxi\mathrm{d}\x}}
\,\int_{\cell{\cellind}}\!\int_{\randomElement{\MEIndex}} \solution_{(0)}(\x,\uncertainty)\, \localxiBasisPoly{\SGsumIndex}{\MEIndex}\polybasis{\xsumIndex}_{\cellind} \, \xiPDFdxi \mathrm{d}\x.
\end{equation}

Finally, the optimization problem for the parameter identification of the discretized conservation law \eqref{eq:DSGM} on the time-continuous level is given by\\
\begin{mdframed}
	\begin{equation}\label{eq:OP}\tag{\textbf{opt}}
		\begin{alignedat}{1}
	&\text{min}_{\uncertaintyInt} \,\costfun(\solution,\uncertaintyInt) = \refone{\text{min}_{\uncertaintyInt}} \sum_{\xsumIndex=0}^{\xdegree}\sum_{\SGsumIndex=0}^\SGtruncorder
	\frac{1}{2}\big\|\solvec^{\xsumIndex,\SGsumIndex}(\tend) - \udata^{\xsumIndex,\SGsumIndex}\big\|_{\refone{2}}^2 + \frac{\xidelta}{2}\big\|\uncertaintyInt-\uncertaintyInt_p\big\|_{\refone{2}}^2\\[.2cm]
	\text{subject to } &\sum_{\xsumIndex=0}^{\xdegree}\sum_{\SGsumIndex=0}^\SGtruncorder
	\bigg(\begin{matrix}
	\WF(\solvec^{\xsumIndex,\SGsumIndex},\uncertaintyInt)\\[3pt]
	\WF_{(0)}(\solvec^{\xsumIndex,\SGsumIndex})
	\end{matrix} \bigg) =
	\sum_{\xsumIndex=0}^{\xdegree}\sum_{\SGsumIndex=0}^\SGtruncorder
	\bigg(\begin{matrix}
	M_1^{\xsumIndex,\SGsumIndex} \partial_\timevar\solvec^{\xsumIndex,\SGsumIndex} -  
	M_2^{\xsumIndex,\SGsumIndex}(\solvec) + 
	M_3^{\xsumIndex,\SGsumIndex} (\solvec)\\[3pt]
	\solvec^{\xsumIndex,\SGsumIndex}(0) - \solvec^{\xsumIndex,\SGsumIndex}_{(0)}
	\end{matrix} \bigg)
	=0 \nonumber
	\end{alignedat}
	\end{equation}
\end{mdframed}

We determine the first-order optimality conditions for this optimization problem in the following.

\subsection{First-order optimality conditions}
The Lagrange functional corresponding to \eqref{eq:OP} reads
\begin{equation}\label{eq:lagrangian}
\lagrange(\solvec,\uncertaintyInt,\adj) = \sum_{\xsumIndex=0}^{\xdegree}\sum_{\SGsumIndex=0}^\SGtruncorder\, \costfun(\solvec^{\xsumIndex,\SGsumIndex},\uncertaintyInt) 
+ \sum_{\xsumIndex=0}^{\xdegree}\sum_{\SGsumIndex=0}^\SGtruncorder\, \intT\! 
\bigg(\begin{matrix}
\WF(\solvec^{\xsumIndex,\SGsumIndex},\uncertaintyInt)\\[3pt]
\WF_{(0)}(\solvec^{\xsumIndex,\SGsumIndex})
\end{matrix} \bigg)\cdot 
\bigg(\begin{matrix}
\adj^{\xsumIndex,\SGsumIndex}\\[3pt]
\adj^{\xsumIndex,\SGsumIndex}_{(0)}
\end{matrix}\bigg) \,\mathrm{d}\timevar,
\end{equation}
where $\adj^{\xsumIndex,\SGsumIndex}=\adj^{\xsumIndex,\SGsumIndex}(\timevar)\in\mathbb{R}^{\ncells\cdot\MEElements}$ and $\adj^{\xsumIndex,\SGsumIndex}_{(0)}\in\mathbb{R}^{\ncells\cdot\MEElements}$ are the Lagrange multipliers that represent the adjoint state.
We denote by $\adj$ the vectorized values $\adj^{\xsumIndex,\SGsumIndex}$ and $\adj^{\xsumIndex,\SGsumIndex}_{(0)}$ for all $\xsumIndex=0,\ldots,\xdegree$ and $\SGsumIndex=0,\ldots,\SGtruncorder$.
The chronological order of the elements in these vectors is chosen according to \rmref{rm:order}.

The first-order optimality condition to \eqref{eq:OP} is given by solving 
$$\mathrm{d}\lagrange(\solvec,\uncertaintyInt,\adj) = 0.$$
The derivative with respect to the adjoint state $\adj$ results in the state equation, while the derivative with respect to the state $\solvec$ yields the adjoint system and the optimality condition is obtained by the derivative with respect to the control $\uncertaintyInt$, c.f. \cite{Hinze2009}.

For arbitrary directions $\directionxi\in\mathbb{R}^2$, $\directionhk=\directionhk(\timevar)\in\mathbb{R}^{\ncells\cdot\MEElements}$, where $\direction$ denotes the vectorized values $\directionhk$ with the same chronological order as in $\solvec^{\xsumIndex,\SGsumIndex}$ and $ \adj^{\xsumIndex,\SGsumIndex}$,
we obtain the following Gâteaux derivatives of the cost functional
\begin{align*}
	\mathrm{d}_{\solvec}  \costfun(\solvec,\uncertaintyInt)[\direction]
	& = \sum_{\xsumIndex=0}^{\xdegree}\sum_{\SGsumIndex=0}^\SGtruncorder (\solvec^{\xsumIndex,\SGsumIndex}(\tend)-\solvec^{\xsumIndex,\SGsumIndex}_D)\cdot \directionhk(\tend),\\[4pt]
	 \mathrm{d}_{\uncertaintyInt} \costfun(\solvec,\uncertaintyInt)[\directionxi] 
	& = \xidelta\,(\uncertaintyInt-\uncertaintyInt_p)\cdot\directionxi.
\end{align*}

The directional derivatives for the second part of the Lagrangian \eqref{eq:lagrangian} read
\begin{align*}
	 &\intT\!  \mathrm{d}_{\solvec^{\xsumIndexvar, \SGsumIndexvar}} \sum_{\xsumIndex=0}^{\xdegree}\sum_{\SGsumIndex=0}^\SGtruncorder \bigg(\begin{matrix}
	\WF(\solvec^{\xsumIndex,\SGsumIndex},\uncertaintyInt)\\[3pt]
	\WF_{(0)}(\solvec^{\xsumIndex,\SGsumIndex})
	\end{matrix}\bigg) [\direction] \cdot \bigg(\begin{matrix}
	\adj^{\xsumIndex,\SGsumIndex}\\[3pt]
	\adj^{\xsumIndex,\SGsumIndex}_{(0)}
	\end{matrix}\bigg) \mathrm{d}\timevar \\[4pt]
	& \hspace*{1.7cm}= \sum_{\xsumIndex=0}^{\xdegree}\sum_{\SGsumIndex=0}^\SGtruncorder \intT\! 
	M_1^{\xsumIndex,\SGsumIndex} \partial_\timevar \directionhk \cdot \adj^{\xsumIndex,\SGsumIndex} 
	+  \Big(\mathrm{d}_{\solvec^{\xsumIndexvar, \SGsumIndexvar}}M_3^{\xsumIndex,\SGsumIndex}(\solvec) - \mathrm{d}_{\solvec^{\xsumIndexvar, \SGsumIndexvar}}M_2^{\xsumIndex,\SGsumIndex}(\solvec)\Big)\, \directionhkvar \cdot \adj^{\xsumIndexvar,\SGsumIndexvar}  \,\mathrm{d}\timevar \\
	&\hspace*{6.2cm}+  \directionhk(0)\cdot\adj^{\xsumIndex,\SGsumIndex}_{(0)}, \\[5pt]
	&\intT\! \mathrm{d}_{\uncertaintyInt} \bigg(\begin{matrix}
	\WF(\solvec^{\xsumIndex,\SGsumIndex},\uncertaintyInt)\\[3pt]
	\WF_{(0)}(\solvec^{\xsumIndex,\SGsumIndex})
	\end{matrix}\bigg)[\directionxi] \cdot 
	\bigg(\begin{matrix}
	\adj^{\xsumIndex,\SGsumIndex}\\[3pt]
	\adj^{\xsumIndex,\SGsumIndex}_{(0)}
	\end{matrix}\bigg) \mathrm{d}\timevar \\[4pt]
	& \hspace*{1.7cm} = \intT\! \bigg(\Big(
\partial_{\uncertaintyInt} M_3^{\xsumIndex,\SGsumIndex}(\solvec)
	- \partial_{\uncertaintyInt} M_2^{\xsumIndex,\SGsumIndex}(\solvec) \Big) \cdot \directionxi \bigg) \cdot\adj^{\xsumIndex,\SGsumIndex} 
	- \Big(\partial_{\uncertaintyInt}\solvec^{\xsumIndex,\SGsumIndex}_{(0)} \cdot \directionxi\Big) \cdot \adj^{\xsumIndex,\SGsumIndex}_{(0)} \,\mathrm{d}\timevar,
\end{align*}
with Lagrange multipliers $\adj^{\xsumIndex,\SGsumIndex}\in\mathbb{R}^{\ncells\cdot\MEElements}$ for $\xsumIndex=0,\ldots,\xdegree$ and $\SGsumIndex=0,\ldots,\SGtruncorder$.
Using \rmref{rem:refCell} yields the matrices $\mathrm{d}_{\solvec^{\xsumIndexvar, \SGsumIndexvar}} M_2^{\xsumIndex,\SGsumIndex}(\solvec)$ and $\mathrm{d}_{\solvec^{\xsumIndexvar, \SGsumIndexvar}} M_3^{\xsumIndex,\SGsumIndex}(\solvec)^T \in \mathbb{R}^{\ncells\cdot\MEElements\times\ncells\cdot\MEElements}$
	\begin{align*}
\mathrm{d}_{\solvec^{\xsumIndexvar, \SGsumIndexvar}} M_2^{\xsumIndex,\SGsumIndex}(\solvec) &= 
	\frac{1}{\Delta\x}\,\text{diag}\Big(\int_{\cell{\cellind}}\!	\int_{\randomElement{\MEIndex}}\! \flux_\solution( \localsolution{\cellind}{\cellindR})\,\polybasis{\xsumIndexvar}_{\cellind}\partial_\x\polybasis{\xsumIndex}_{\cellind}\,
	\localxiBasisPoly{\SGsumIndexvar}{\MEIndex}\localxiBasisPoly{\SGsumIndex}{\MEIndex}\,\xiPDFdxi\mathrm{d}\x  \Big)\,,\\
	\mathrm{d}_{\solvec^{\xsumIndexvar, \SGsumIndexvar}} M_3^{\xsumIndex,\SGsumIndex}(\solvec)^T& = \frac{1}{2\Delta\x}\,\text{tridiag}\!\begin{pmatrix*}[l]
		\hspace*{0.8cm}\int_{\randomElement{\MEIndex}}\!\big(\flux_\solution\big(\localsolution{\cellind}{\cellindR} (\x_{\cellind+\frac12})\big)-\viscosityconstant\big)\,\polybasis{\xsumIndex}_{\cellind}(\x_{\cellind+\frac12})\polybasis{\xsumIndexvar}_{\cellind}(\x_{\cellind-\frac12})\localxiBasisPoly{\SGsumIndex}{\MEIndex}\localxiBasisPoly{\SGsumIndexvar}{\MEIndex}\,\xiPDFdxi\\
	\int_{\randomElement{\MEIndex}}\!\big(\flux_\solution\big(\localsolution{\cellind}{\cellindR} (\x_{\cellind+\frac12})\big)+\viscosityconstant\big)\,\polybasis{\xsumIndex}_{\cellind}(\x_{\cellind+\frac12})\polybasis{\xsumIndexvar}_{\cellind}(\x_{\cellind+\frac12})\localxiBasisPoly{\SGsumIndex}{\MEIndex}\localxiBasisPoly{\SGsumIndexvar}{\MEIndex}\,\xiPDFdxi\\
	\hspace*{0.9cm}+\int_{\randomElement{\MEIndex}}\!\big(-\flux_\solution\big(\localsolution{\cellind}{\cellindR} (\x_{\cellind-\frac12})\big)+\viscosityconstant\big)\,\polybasis{\xsumIndex}_{\cellind}(\x_{\cellind-\frac12})\polybasis{\xsumIndexvar}_{\cellind}(\x_{\cellind-\frac12})\localxiBasisPoly{\SGsumIndex}{\MEIndex}\localxiBasisPoly{\SGsumIndexvar}{\MEIndex}\,\xiPDFdxi\\
	\hspace*{0.6cm}\int_{\randomElement{\MEIndex}}\!\big(-\flux_\solution\big(\localsolution{\cellind}{\cellindR} (\x_{\cellind-\frac12})\big)-\viscosityconstant\big)\,\polybasis{\xsumIndex}_{\cellind}(\x_{\cellind-\frac12})\polybasis{\xsumIndexvar}_{\cellind}(\x_{\cellind+\frac12})\localxiBasisPoly{\SGsumIndex}{\MEIndex}\localxiBasisPoly{\SGsumIndexvar}{\MEIndex}\,\xiPDFdxi \,\\
	\end{pmatrix*}\!.
	\end{align*}
We require the transposed matrix $\mathrm{d}_{\solvec^{\xsumIndexvar, \SGsumIndexvar}} M_3^{\xsumIndex,\SGsumIndex}(\solvec)^T$ for the adjoint problem in the following theorem which is why we do not state $\mathrm{d}_{\solvec^{\xsumIndexvar, \SGsumIndexvar}} M_3^{\xsumIndex,\SGsumIndex}(\solvec)$ explicitly.
For the derivatives with respect to the distribution parameters we obtain
	\begin{align*}
\partial_{\uncertaintyInt}\big(M_1^{\xsumIndex,\SGsumIndex}\solvec^{\xsumIndex,\SGsumIndex}\big) &=\partial_{\uncertaintyInt}\bigg( \text{diag}\Big(\int_{\refcell}\!(\polybasis{\xsumIndex})^2\,\mathrm{d}\x 
\int_{\refME}\!(\phi^{\SGsumIndex})^2\,\mathrm{d}\uncertainty \Big)\,\solvec^{\xsumIndex,\SGsumIndex}\bigg)=0,\\
\partial_{\uncertaintyR}\big(M_2^{\xsumIndex,\SGsumIndex}(\solvec)\big) &= 
\frac{1}{\Delta\x}\,\text{vec}\bigg(\int_{\refcell}\!
\int_{\refME}\!\flux_\uncertainty\Big(\localsolution{\cellind}{\cellindR},\,\big(\uncertaintyR-\uncertaintyL\big)\frac{\big(\uncertainty+\MEIndex-\frac12\big)}{\MEElements}+\uncertainty_L\Big)\Big(\frac{\uncertainty+\MEIndex-\frac12}{\MEElements}\Big)
\partial_\x\polybasis{\xsumIndex} \,\phi^{\SGsumIndex}\,\mathrm{d}\uncertainty\mathrm{d}\x\bigg),
\end{align*}
where $\partial_{\uncertaintyL}M_2^{\xsumIndex,\SGsumIndex}(\solvec)$ and $\partial_{\uncertaintyInt}M_3^{\xsumIndex,\SGsumIndex}(\solvec)$ follow analogously.  Moreover, we denote by $\flux_\solution$ and $\flux_\uncertainty$ the derivative of the flux function $\flux(\solution,\uncertainty)$ with respect to the first and second component, respectively. The term $(\uncertaintyR-\uncertaintyL)(\uncertainty+\MEIndex-\frac12)/\MEElements+\uncertainty_L$ describes the variable transformation $\uncertainty\,\Delta\uncertainty+\uncertainty_\MEIndex$ using $\uncertaintyL$ and $\uncertaintyR$. Furthermore, we have
$$\partial_{\uncertaintyR}\!\big(\solvec^{\xsumIndex,\SGsumIndex}_{(0)}\big) = \big(M_1^{\xsumIndex,\SGsumIndex}\big)\!^{-1}
\!\text{diag}\bigg(\!\int_{\refcell}\!\int_{\refME}\!\solution_{(0)}'\Big(\!\!\x\,\Delta\x+\x_i,\,\big(\uncertaintyR-\uncertaintyL\big)\!\frac{\big(\uncertainty+\MEIndex-\frac12\big)}{\MEElements}+\uncertainty_L\Big)\!\!\Big(\frac{\uncertainty+\MEIndex-\frac12}{\MEElements}\Big)\polybasis{\xsumIndex}\phi^{\SGsumIndex}\,\mathrm{d}\uncertainty\mathrm{d}\x\!\!\bigg)\!.$$
Here, $\solution_{(0)}'$ describes the derivative of $\solution_{(0)}$ with respect to the second argument.
\begin{remark}
	In order to calculate these quantities,  the flux $\flux(\solution,\uncertainty)$ has to be differentiable with respect to $\solution$ and $\uncertainty$ and $\solution_{(0)}(\x,\uncertainty)$ has to be differentiable with respect to $\uncertainty$. 
\end{remark}
\begin{theorem}
Let $\big(\uncertaintyInt_*,\solvec_*\big)$ be an optimal pair for \eqref{eq:OP}. Then the following first-order optimality condition holds
\begin{align}\label{eq:opticond}
\sum_{\xsumIndex=0}^{\xdegree}\sum_{\SGsumIndex=0}^\SGtruncorder \intT\!  \Big(\partial_{\uncertaintyInt} M_3^{\xsumIndex,\SGsumIndex} (\solvec_*)
&- \partial_{\uncertaintyInt} M_2^{\xsumIndex,\SGsumIndex} (\solvec_*) \Big) \cdot\adj^{\xsumIndex,\SGsumIndex}\,\mathrm{d}\timevar +  \partial_{\uncertaintyInt}\solvec^{\xsumIndex,\SGsumIndex}_{(0)} \cdot \big(M_1^{\xsumIndex,\SGsumIndex} \adj^{\xsumIndex,\SGsumIndex}(0)\big)
=  \xidelta\,(\uncertaintyInt_p-\uncertaintyInt_*),
\end{align}
where $\adj^{\xsumIndex,\SGsumIndex} $ satisfies the adjoint problem given by
\begin{equation}\label{eq:adjsys}
	M_1^{\xsumIndex,\SGsumIndex} \,\partial_\timevar \adj^{\xsumIndex,\SGsumIndex}  
	= \sum_{\xsumIndexvar=0}^{\xdegree}\sum_{\SGsumIndexvar=0}^\SGtruncorder \Big(\big(\mathrm{d}_{\solvec^{\xsumIndexvar, \SGsumIndexvar}}M_3^{\xsumIndex,\SGsumIndex}(\solvec_*)\big)^T - \mathrm{d}_{\solvec^{\xsumIndexvar, \SGsumIndexvar}}M_2^{\xsumIndex,\SGsumIndex}(\solvec_*)\Big)\,\adj^{\xsumIndexvar,\SGsumIndexvar},
\end{equation}
for $\xsumIndex=0,\ldots,\xdegree$ and $\SGsumIndex=0,\ldots,\SGtruncorder$ with terminal condition $M_1^{\xsumIndex,\SGsumIndex}\adj^{\xsumIndex,\SGsumIndex}(\tend) = (\solvec^{\xsumIndex,\SGsumIndex}_D-\solvec_*^{\xsumIndex,\SGsumIndex})$.
\end{theorem}
\begin{proof} For the optimal pair, we have $\mathrm{d}_\solvec \lagrange(\solvec_*,\uncertaintyInt_*,\adj) =0$, where we derive via integration by parts 
\begin{align*}
	\intT \!M_1^{\xsumIndex,\SGsumIndex}\partial_t\directionhk\cdot\adj^{\xsumIndex,\SGsumIndex}\,\mathrm{d}\timevar 
	&= \intT\! \partial_t\big(M_1^{\xsumIndex,\SGsumIndex}\directionhk\big)\cdot\adj^{\xsumIndex,\SGsumIndex}\,\mathrm{d}\timevar  \\
	&= \big[\directionhk\cdot\big(M_1^{\xsumIndex,\SGsumIndex}\big)^T\adj^{\xsumIndex,\SGsumIndex}\big]_0^T - \intT\!\directionhk\cdot\big(M_1^{\xsumIndex,\SGsumIndex}\big)^T\partial_\timevar\adj^{\xsumIndex,\SGsumIndex}\,\mathrm{d}\timevar,\\
		\intT \big(M_3^{\xsumIndex,\SGsumIndex}(\solvec_*) - M_2^{\xsumIndex,\SGsumIndex}(\solvec_*))\big)\directionhkvar\cdot\adj^{\xsumIndexvar,\SGsumIndexvar} \,\mathrm{d}\timevar
	&= \intT \directionhkvar\cdot\big(M_3^{\xsumIndex,\SGsumIndex}(\solvec_*) - M_2^{\xsumIndex,\SGsumIndex}(\solvec_*)\big)^T\adj^{\xsumIndexvar,\SGsumIndexvar} \,\mathrm{d}\timevar.
\end{align*}
Choosing smooth functions with $\directionhk(0)=\directionhk(T)=0$ and $\directionhk(t) \ne 0$ for $t\in (0,T)$ as test functions, together with the symmetry of $M_1^{\xsumIndex,\SGsumIndex}$ and $M_2^{\xsumIndex,\SGsumIndex}(\solvec_*)$, we obtain \eqref{eq:adjsys}. 
Choosing $\directionhk$ smooth with $\directionhk(0)=\directionhk(T)\ne 0$ yields the terminal condition, where $M_1^{\xsumIndex,\SGsumIndex}$ has full rank, thus it is invertible. Now, using the conditions derived so far and $\directionhk\ne0 $,  we get $\adj^{\xsumIndex,\SGsumIndex}_{(0)} = M_1^{\xsumIndex,\SGsumIndex}\adj^{\xsumIndex,\SGsumIndex}(0)$, together with $\mathrm{d}_{\uncertaintyInt} \lagrange(\solvec,\uncertaintyInt,\adj) =0$ leads to \eqref{eq:opticond}.
\end{proof}
\begin{remark}\label{rem:opt_type}
Note, that if the flux $\flux(\solution,\uncertainty)$ does not depend on $\uncertainty$, the optimality condition \eqref{eq:opticond} only contains the optimization in the uncertain initial condition through the second term $\partial_{\uncertaintyInt}\solvec^{\xsumIndex,\SGsumIndex}_{(0)} \cdot \big(M_1^{\xsumIndex,\SGsumIndex}\adj^{\xsumIndex,\SGsumIndex}(0)\big)$. On the other hand, if the initial condition is independent of $\uncertainty$, this term vanishes.
\end{remark}

\subsection{Reduced Cost Functional}
\refone{In the derivation of the first-order optimality conditions of \eqref{eq:OP} we used the state and control dependence of the cost functional $\costfun(\solution,\uncertaintyInt)$ explicitly. Now, we propose a gradient-based optimization algorithm that acts only on the control. We therefore define the control-to-state operator $$\controltostate(\uncertaintyInt)=\solvec,$$
that is deduced through application of the discontinuous stochastic Galerkin scheme \eqref{eq:DG}. The existence and uniqueness of solutions is discussed in \rmref{rem:existunique}.	Further, we introduce the corresponding reduced cost functional
\begin{equation}\label{eq:redcostfun} \redcostfun(\uncertaintyInt) := \costfun(\controltostate(\uncertaintyInt),\uncertaintyInt) = \costfun(\solvec,\uncertaintyInt). \end{equation}} 
\refone{In the following we assume the control-to-state operator to be Gâteaux differentiable,} otherwise the so-called shift derivative has to be introduced, for more information we refer to \cite{Hajian2019,Ulbrich2001}.
Additionally, we define 
$$ \Fvar(\solvec,\uncertaintyInt,\adj) := \sum_{\xsumIndex=0}^{\xdegree}\sum_{\SGsumIndex=0}^\SGtruncorder \bigg(\begin{matrix}
\WF(\solvec^{\xsumIndex,\SGsumIndex},\uncertaintyInt)\\[3pt]
\WF_{(0)}(\solvec^{\xsumIndex,\SGsumIndex})
\end{matrix}\bigg) \cdot \bigg(\begin{matrix}
\adj^{\xsumIndex,\SGsumIndex}\\[3pt]
\adj^{\xsumIndex,\SGsumIndex}_{(0)}
\end{matrix}\bigg).$$
Analogous to \cite{Hinze2009,Totzeck}, using  $\Fvar(\solvec,\uncertaintyInt,\adj)=0$, we obtain 
$$0 = \mathrm{d}_{\uncertaintyInt}\Fvar(\controltostate(\uncertaintyInt),\uncertaintyInt,\adj)
=\mathrm{d}_\solvec\Fvar(\controltostate(\uncertaintyInt),\uncertaintyInt,\adj)[\mathrm{d}\controltostate(\uncertaintyInt)] 
+ \mathrm{d}_{\uncertaintyInt}\Fvar(\controltostate(\uncertaintyInt),\uncertaintyInt,\adj).$$
The derivative of the Lagrangian with respect to the state $\solvec$ yields
$$\mathrm{d}_\solvec\Fvar(\solvec,\uncertaintyInt,\adj)^*=-\mathrm{d}_\solvec\costfun(\solvec,\uncertaintyInt).$$
With these equations, the  Gâteaux derivative of the reduced cost functional in arbitrary direction $\directionxi\in\mathbb{R}^2$ is given by
\begin{subequations}\label{eq:derredcostfun}
\begin{align}
\mathrm{d}\redcostfun(\uncertaintyInt)[\directionxi]&=\mathrm{d}_\solvec\costfun(\solvec,\uncertaintyInt)\cdot\mathrm{d}\controltostate(\uncertaintyInt)[\directionxi] + \mathrm{d}_{\uncertaintyInt}\costfun(\solvec,\uncertaintyInt)[\directionxi]\\
&=\mathrm{d}_{\uncertaintyInt}\Fvar(\controltostate(\uncertaintyInt),\uncertaintyInt,\adj) \cdot \directionxi + \mathrm{d}_{\uncertaintyInt}\costfun(\solvec,\uncertaintyInt)[\directionxi]\\
&=\mathrm{d}_{\uncertaintyInt}\lagrange(\solvec,\uncertaintyInt,\adj) \cdot  \directionxi ,
\end{align}
\end{subequations}
which corresponds to the first-order optimality condition from \eqref{eq:opticond}.
\section{Numerical Scheme}\label{sec:scheme}
For the numerical solution of our optimization problem \eqref{eq:OP} we first have to solve the discretized forward problem \eqref{eq:spde} as described in \secref{sec:DSG} through the discontinuous stochastic Galerkin method in order to obtain the state variable for predefined equally sized time slices of $[0,\,\tend]$, namely $\timevar_\timeind = \timeind\,\Delta\timevar, \timeind=0,\ldots,\nbtime$ with $\Delta\timevar=\frac{\tend}{\nbtime}$. We will use these time slices later to solve the time integral in the optimality condition \eqref{eq:opticond} via quadrature. \refone{For the quadrature rule with respect to $\timevar$, we might use the equidistant time slices with Newton-Cotes weights.} Moreover, we store the state information at all times $\timevar$ that are used in the adjoint system \eqref{eq:adjsys}.
We solve the adjoint problem backwards in time with \refone{Runge-Kutta} which yields the adjoint state at our time slices. 

Note that \cite{Hager2000} showed a reduced order of Runge Kutta in solving the adjoint system with more than two stages, this issue is also discussed for conservation laws in \cite{Hintermuller2018,Hajian2019}. Using a three stage SSP RK method will only yield a second-order approximation in the adjoint problem. Additional conditions have to be imposed on the RK scheme to reach third- and fourth-order methods in the backward solution, c.f. \cite{Hager2000}. According to \cite{Hajian2019}, our SSP TVD RK discretization of the state equation ensures that the discrete adjoint is stable in every time step. Moreover, in the case of the Lax Friedrichs numerical flux, we choose the time step size for the forward and backward scheme as half of the usual CFL number. In our numerical calculations we transformed the adjoint system to an initial value problem using a new time variable $\tend-\timevar$. 

\refone{
\begin{remark}
In general, the CFL number of the state and the adjoint equation differ. Here, the 'first optimize, then discretize' approach on the ODE level is advantageous, as we can use different time discretizations for the state and the adjoint systems. When the state equation is discretized first, there is less choice for the adjoint system.  		
\end{remark}
}

Finally, we perform a optimization method in order to update the control. 
We first consider a steepest descent method where we update
\begin{equation}\label{eq:steepdes}
\uncertaintyInt_{\text{new}} = \uncertaintyInt  - \alpha \,\mathrm{d} \redcostfun(\uncertaintyInt),
\end{equation}
with the reduced cost functional from \eqref{eq:redcostfun}.
The step size $\alpha>0$ is obtained by the Armijo step size rule \cite{Hinze2009}. 
Moreover, $\uncertaintyInt$ denotes the control in the current iteration, whereas at time $t=0$ we start with the initial guess $\uncertaintyInt_p$.
 According to \eqref{eq:derredcostfun}, the derivative of the reduced cost functional given by 
\begin{subequations}\label{eq:gradJ}
\begin{align}
\mathrm{d}\redcostfun(\uncertaintyInt) =\!\!\intT\!\!\! \Big(\sum_{\xsumIndex=0}^{\xdegree}\sum_{\SGsumIndex=0}^\SGtruncorder \partial_{\uncertaintyInt} M_3^{\xsumIndex,\SGsumIndex} (\solvec)
&- \partial_{\uncertaintyInt} M_2^{\xsumIndex,\SGsumIndex} (\solvec) \Big) \cdot\adj^{\xsumIndex,\SGsumIndex}\,\mathrm{d}\timevar \\
&+  \partial_{\uncertaintyInt}\solvec^{\xsumIndex,\SGsumIndex}_{(0)} \cdot \big(M_1^{\xsumIndex,\SGsumIndex} \adj^{\xsumIndex,\SGsumIndex}(0)\big)
+  \xidelta\,(\uncertaintyInt-\uncertaintyInt_p),
\end{align}
\end{subequations}
where the integrals over $[0, \tend]$ are solved using an appropriate quadrature rule with points $\timevar_0,\ldots,\timevar_{\nbtime}$.

Alternatively, we can perform the optimization process through a BFGS scheme, described for example in \cite{FredericBonnans2006}, where we approximate the inverse of the Hessian matrix in each iteration and use it within a quasi-Newton scheme.
Here, the update of the control \refone{is similar to the steepest descent method \eqref{eq:steepdes} but with an additionally introduced matrix $H^{-1}$ and} reads 
\begin{equation}\label{eq:BFGS}
\uncertaintyInt_{\text{new}} = \uncertaintyInt  - \alpha \,H^{-1}\mathrm{d} \redcostfun(\uncertaintyInt),
\end{equation}
where $\alpha>0$ again denotes the step size of the Armijo rule and $H^{-1}$ is the approximate inverse of the Hessian. We initially use the identity matrix for $H^{-1}$ and then update it for the next iteration \refone{by a rank-two correction of the inverse} to
$$H_{\text{new}}^{-1} = H^{-1} - \frac{s^TyH^{-1} + H^{-1}y^Ts}{y^Ts} + \bigg(1+\frac{y^TH^{-1}y}{y^Ts}\bigg)\frac{ss^T}{y^Ts}\,,$$
with $s=-\alpha H^{-1}\mathrm{d} \redcostfun(\uncertaintyInt)$ and $y=\mathrm{d} \redcostfun(\uncertaintyInt_{\text{new}}) - \mathrm{d} \redcostfun(\uncertaintyInt)$. We reset $H^{-1}$ back to the identity matrix if $-(\uncertaintyInt)^T(H^{-1}\mathrm{d} \redcostfun(\uncertaintyInt))>0$. \refone{This formula adds two symmetric rank-one matrices to $H^{-1}$ in order to update the Hessian and is therefore known as a rank-two correction method. The exact update can be obtained by an expansion of the gradient with help of the Hessian and by approximating the resulting gradient difference, for more information see \cite{FredericBonnans2006}.}

We perform the whole method until the determination of the parameter $\uncertaintyInt_{\text{new}}$ has converged, namely until the $L_2$ norm of the direction of descent $\|\mathrm{d}\redcostfun(\uncertaintyInt_{\text{new}})\|_{\refone{2}}$ is smaller than a predefined tolerance.

The whole numerical scheme is summarized in the following algorithm.
\begin{algorithm}[H] 
	\caption{Discontinuous stochastic Galerkin (DsG) scheme for parameter estimation}
	\label{algo:OCalgo}
	\begin{algorithmic}[1]
		\STATE Set $L = \|\mathrm{d}\redcostfun(\uncertaintyInt)\|_{\refone{2}}$\smallskip
		\WHILE{$L > \text{tol}$}
		\smallskip
		\STATE Solve \eqref{eq:spde} using the DsG scheme \eqref{eq:DG} to calculate $\solvec(\timevar_\timeind)$ for all $\timeind=0,\ldots,\nbtime$\smallskip\STATE Solve the adjoint system \eqref{eq:adjsys} backwards to get $\adj(\timevar_\timeind)$, $\timeind=0,\ldots,\nbtime$\smallskip
		\STATE Use the Armijo rule to compute the step size and update with \eqref{eq:steepdes} or \eqref{eq:BFGS} in order to obtain $\uncertaintyInt_{\text{new}}$ \smallskip
		\STATE Update $L = \|\mathrm{d}\redcostfun(\uncertaintyInt_{\text{new}})\|_{\refone{2}}$ \smallskip
		\ENDWHILE
	\end{algorithmic}
\end{algorithm}

\begin{remark}\label{rem:LA}
	If we only consider an uncertain linear advection equation, namely if the flux is given by $\flux(\solution,\uncertainty) =  \advec(\uncertainty)\,\solution$, the discontinuous Galerkin scheme is linear in $\solvec$, meaning that we can write $M_2^{\xsumIndex,\SGsumIndex}(\solvec)=M_2^{\xsumIndex,\SGsumIndex}\solvec$ and $M_3^{\xsumIndex,\SGsumIndex}(\solvec) = M_3^{\xsumIndex,\SGsumIndex}\solvec$ in \eqref{eq:DSGM}.
	Thus, the derivatives simplify and lead to an adjoint system of the form 
	$$M_1^{\xsumIndex,\SGsumIndex} \,\partial_\timevar \adj^{\xsumIndex,\SGsumIndex}  
	=  \big(M_3^{\xsumIndex,\SGsumIndex} - M_2^{\xsumIndex,\SGsumIndex}\big)^T\,\adj.$$
	In this case, we do not have to store $\solvec$ at every point in time of the backward model and we can derive $M_2^{\xsumIndex,\SGsumIndex}$ and $M_3^{\xsumIndex,\SGsumIndex}$ independently for every $\uncertaintyInt$ which decreases the run time of the algorithm tremendously.
\end{remark}

\section{Numerical Results}\label{sec:results}
In the following numerical experiments we analyze the discontinuous stochastic Galerkin scheme for parameter estimation described in \algoref{algo:OCalgo}. We apply the method to the stochastic linear advection and Burgers' equations, supplemented with different types of uncertainty, and study the identification of a predefined reference solution as well as the influence of the discretization parameters.

\subsection{Linear Advection}
We begin the numerical calculations by considering the linear advection equation
\begin{subequations}\label{eq:LA}
	\begin{align}
	\frac{\partial}{\partial\timevar}\solution(\timevar ,\x, \uncertainty)+ \advec(\uncertainty) \frac{\partial}{\partial\x} \solution(\timevar ,\x, \uncertainty)&= 0, \quad \quad \text{for} \ \ \x\in\Domain,\, \timevar\in [0,\,\tend],\,\uncertainty\in\SD,\\
	\solution(0,\x,\uncertainty) &= \solution_{(0)}(\x,\uncertainty), \quad \quad \text{for} \ \ \x\in\PD,\, \xi\in\SD,
	\end{align}
\end{subequations}
with uncertain wave speed $\advec=\advec(\uncertainty)$. The flux is therefore given by $\flux(\solution,\uncertainty) =  \advec(\uncertainty)\,\solution$, thus we can apply \rmref{rem:LA}.
In each of the following test cases, the BFGS scheme was able to outperform the steepest descent method in terms of the number of iterations and Armijo steps, which is why we restrict the results to BFGS.

\subsubsection{Parameter estimation under discontinuous initial conditions}\label{sec:LAshock}
In this subsection we study \algoref{algo:OCalgo} applied to the uncertain linear advection equation \eqref{eq:LA} with discontinuous initial conditions, namely
\begin{equation}\label{eq:initial_shock}\solution_{(0)}(\x) = \begin{cases}1,&\text{for }0.4<\x<0.6,\\ 0,&\text{for }0\leq\x\leq 0.4 \text{ and } 0.6\leq\x\leq 1.\end{cases}\end{equation}

The uncertain wave speed is given by
$$\advec(\uncertainty) = 2\uncertainty.$$

Moreover, the uncertainty is chosen to be uniform $\uncertainty\sim\mathcal{U}(\uncertaintyL,\uncertaintyR)$ and we calculate the data $u_D$, i.e., the reference solution, at time $\tend = 0.01$ with $\uncertainty\sim\mathcal{U}(-1,1)$. We start with an initial guess \reftwo{$\uncertaintyInt_p = [-0.8,\,1.4]$}. For our numerical calculations we use $\ncells=200$ spatial cells and $\MEElements=20$ \ME{}s, the polynomial degrees of the local approximation read $\xdegree =1$ and $\SGtruncorder = 4$, providing a two-stage SSP RK method. In the cost functional, we set the scaling factor of $\costfunxi$ as $\xidelta=10^{-2}$. We perform the optimization in \algoref{algo:OCalgo} until the tolerance $tol=10^{-2}$ is reached. The results are shown in \figref{fig:shock_xi} -- \ref{fig:shock_u-uD}. \refone{The influence of the parameters for \algoref{algo:OCalgo} will be discussed later within \tabref{tab:params}. In general, a solid resolution of the spatial and stochastic domain should be given and $\xidelta$ is chosen such that $\costfunxi$ does not dominate in $\costfun$, whereas the values $10^{-2}$ or $10^{-3}$ were a good choice in our examples. We set $tol$ such that the optimization process stopped as soon as there seemed no improvement in $\uncertaintyInt$ getting closer to its true value. }

We observe that the parameters $\uncertaintyL$ and $\uncertaintyR$ of the uniform distribution converge to the reference solution, namely to $[-1,\,1]$. In \figref{fig:shock_xi} they quickly approach these values starting at \reftwo{$[-0.8,\,1.4]$} until they slightly propagate around the true value. We need \reftwo{12} iterations in order to achieve the tolerance $10^{-2}$. 

The plots of the cost functional in \figref{fig:shock_J} illustrate \reftwo{the decrease of $\costfun(\solvec,\uncertaintyInt)$  from around 12 in the first iteration until it convergences to $10^{-4}$.} The second part $\costfunxi(\uncertaintyInt,\uncertaintyInt_p)$ of this functional in \figref{fig:shock_Jxi} approaches a value larger than the initial one since we chose our initial guess \reftwo{$[-0.8,\,1.4]$}  as $\uncertaintyInt_p$. This is why we scale this part through $\xidelta\ll1$. 

Furthermore, we display the \reftwo{$L_2$ norm} of the difference between the current solution at time $\tend=0.01$ and the data $\solution_D$ in \figref{fig:shock_u-uD}, which consequently takes a similar form to the cost functional. This value approaches $10^{-5}$ which reflects our computational error. The number of Armijo stages within our calculations was between 4 and 5 in each iteration, starting with the initial value 0.125.

\externaltikz{LAshock}{
\begin{figure}[htb]
	\centering
	\hspace*{-.4cm}
	\begin{minipage}[b]{0.45\textwidth}
 \begin{tikzpicture}
  \begin{axis}[width=1.2\figurewidth, height=1.1\figureheight, xlabel=$\#$Iteration, ylabel=$\xi$, ylabel style = {rotate=-90}, legend style={at={(0.97,0.7)}}] 
   \addplot[mark=none, solid, tuklred, line width=.7pt] file {Images/LAshock/LAshock_xiL.txt};
   \addplot[mark=none, dashed, tuklblue, line width=.8pt] file {Images/LAshock/LAshock_RefxiL.txt};
   \addplot[mark=none, solid, tuklred, line width=.7pt] file {Images/LAshock/LAshock_xiR.txt};
   \addplot[mark=none, dashed, tuklblue, line width=.8pt] file {Images/LAshock/LAshock_RefxiR.txt};
   \legend{$\uncertaintyInt$, $\uncertaintyInt_{\text{Ref}}$}
  \end{axis}
 \end{tikzpicture}
 \caption{\reftwo{Performance of $\xi$ during the optimization process. Example \eqref{eq:initial_shock}.}}
 \label{fig:shock_xi}
 	\end{minipage}
 \hspace*{.8cm}
 \begin{minipage}[b]{0.45\textwidth}
 	\begin{tikzpicture}
 \begin{axis}[width=1.2\figurewidth, height=1.1\figureheight, xlabel=$\#$Iteration, ylabel=$\costfun$, ylabel style = {rotate=-90}, ymode = log] 
 \addplot[mark=none, solid, tuklblue, line width=.7pt] file {Images/LAshock/LAshock_J.txt};
 \end{axis}
 \end{tikzpicture}
 \caption{\refone{Cost functional during the optimization process.} \reftwo{Example \eqref{eq:initial_shock}.}}	
 \label{fig:shock_J}
\end{minipage}
\end{figure}

\begin{figure}[htb]
	\centering
	\hspace*{-.4cm}
	\begin{minipage}[b]{0.45\textwidth}
	\begin{tikzpicture}
\begin{axis}[width=1.2\figurewidth, height=1.1\figureheight, xlabel=$\#$Iteration, ylabel=$\costfunxi$, ylabel style = {rotate=-90}] 
\addplot[mark=none, solid, tuklblue, line width=.7pt] file {Images/LAshock/LAshock_Jxi.txt};
\end{axis}
\end{tikzpicture}
	\caption{\reftwo{$2^{\text{nd}}$ term of cost functional during the optimization process. Example \eqref{eq:initial_shock}.}}
\label{fig:shock_Jxi}
	\end{minipage}
	\hspace*{.6cm}
\begin{minipage}[b]{0.45\textwidth}
	\begin{tikzpicture}
\begin{axis}[width=1.2\figurewidth, height=1.1\figureheight, xlabel=$\#$Iteration, ylabel=\reftwo{$\|\solution-\solution_D\|_{2}$},  y tick label style={ /pgf/number format/.cd, fixed,precision = 6,/tikz/.cd}, ymode = log] 
\addplot[mark=none, solid, tuklblue, line width=.7pt] file {Images/LAshock/LAshock_u-uD.txt};
\end{axis}
\end{tikzpicture}
\caption{\reftwo{$L_2$ norm of $\solution(\tend) - \solution_D(\tend)$ during the optimization process.} \refone{Example \eqref{eq:initial_shock}.}}
\label{fig:shock_u-uD}
	\end{minipage}
\end{figure}


}

\subsubsection[Choice of Parameters in Algorithm 1]{Choice of Parameters in \algoref{algo:OCalgo}}\label{sec:smoothtest}
We supplement the uncertain linear advection equation \eqref{eq:LA} with smooth initial conditions 
\begin{equation}\label{eq:initial_sinus}\solution_{(0)}(\x) = \sin(2\pi\x), \qquad \x\in[0,\,1]  \end{equation}
and periodic boundary conditions, whereas the advection term is given by
$$\advec(\uncertainty) = \uncertainty.$$
We perform \algoref{algo:OCalgo} with different parameter settings in order to study how they influence the convergence of $\uncertaintyInt$ to a reference solution $\uncertainty\sim\mathcal{U}(-1,1)$ at time $\tend=0.01$ with the prior $\uncertaintyInt_p=[1,\,-1]$. We initially start with $\alpha=1$ for the Armijo step size rule in each of the calculations.

The results are described in \tabref{tab:params}. They show that the choice of the parameters does not influence the number of iterations in \algoref{algo:OCalgo}. In each case, we require around 18 -- 21 iterations in order to achieve the tolerance $\|\mathrm{d}\redcostfun(\uncertaintyInt)\|_{\refone{2}}<1e-5$.

In \tabref{tab:Nx} and \tabref{tab:ME} we increase the number of cells in the physical and stochastic domain, and observe that the number of iterations stays the same while the run time of the algorithm is increased and the resulting parameter estimate $\uncertaintyInt$ gets closer to the true value, i.e., the cost functional is decreased. This indicates the mesh-independence of our algorithm. However, using $\MEElements=1$ and therefore no \ME{} ansatz, we are not able to identify the distribution parameters since the terminal condition of the adjoint system, hence the adjoint state and thus the gradient of the reduced cost functional, get too small.

Increasing the polynomial degrees $\xdegree$ and $\SGtruncorder$ of the discontinuous stochastic Galerkin approach in \tabref{tab:Kx} and \tabref{tab:KO} results in the same number of iterations and estimates $\uncertaintyInt$ in each calculation while the run time is increased. The choice of $\xidelta$ in the second part of the cost functional is demonstrated in \tabref{tab:delta}. Without introducing this parameter, i.e., setting $\xidelta=1$, the optimization method is not able to find an direction of descent since the impact of the prior guess $\uncertaintyInt_p$ is too large within the cost functional. If we neglect $\uncertaintyInt_p$ completely and use $\xidelta=0$, we require twice as much iterations compared to every other parameter setting. The remaining values for $\xidelta$ yield similar results as before. Altogether we deduce that the choice of the different parameters does not influence the number of iterations in most of the cases but the run time of the algorithm as well as the approximation quality of the parameter estimation.

 \begin{table}[!h]
	\begin{subtable}[c]{0.48\textwidth}
		\begin{tabular*}{\textwidth}{C{.7cm}C{1cm}C{1.3cm}c}
			\toprule
			\begin{tabular}{@{}c@{}} \mbox{}\\\mbox{}\end{tabular} 
			$\ncells$ & $\#$ it. & time$[$s$]$ & $\uncertaintyInt$ \\
			\midrule
			100 & 18 & 25 & $[-0.979,\,0.979]$\\
			300 & 18 & 308 & $[-0.990,\,0.990]$\\
			500 & 19 & 435 & $[-0.994,\,0.994]$\\
			700 & 20 & 990 & $[-0.996,\,0.996]$\\
			\bottomrule
		\end{tabular*}
		\subcaption{$\MEElements=20$, $\xdegree=1$, $\SGtruncorder=4$, $\xidelta=1e-2$, $tol=1e-5$
		}
		\label{tab:Nx}	
	\end{subtable}
	\hspace*{0.15cm}
	\begin{subtable}[c]{0.48\textwidth}
		\begin{tabular*}{\textwidth}{C{.7cm}C{1cm}C{1.3cm}c}
			\toprule
			\begin{tabular}{@{}c@{}} \mbox{}\\\mbox{}\end{tabular} 
			$\MEElements$ & $\#$ it. & time$[$s$]$ & $\uncertaintyInt$\\
			\midrule
			1 & \multicolumn{3}{c}{$\|\mathrm{d}\redcostfun(\uncertaintyInt)\|<tol$ in the first it.}\\
			20 & 18 & 308 & $[-0.989,\,0.989]$\\
			40 & 21 & 817 & $[-0.995,\,0.995]$\\
			80 &  21 & 6558 & $[-0.997,\,0.997]$\\
			\bottomrule
		\end{tabular*}
		\subcaption{$\ncells=300$, $\xdegree=1$, $\SGtruncorder=4$, $\xidelta=1e-2$, $tol=1e-5$
		}
		\label{tab:ME}
	\end{subtable}
	\begin{subtable}[c]{0.48\textwidth}
		\begin{tabular*}{\textwidth}{C{.7cm}C{1cm}C{1.3cm}c}
			\multicolumn{4}{c}{}\\
			\multicolumn{4}{c}{}\\
			\toprule
			\begin{tabular}{@{}c@{}} \mbox{}\\\mbox{}\end{tabular} $\xdegree$& $\#$ it. & time$[$s$]$ & $\uncertaintyInt$\\
			\midrule
			0 & 18 & 82 & $[-0.992,\,0.992]$\\
			1 & 19 & 435 & $[-0.993,\,0.993]$\\
			2 &  19 &  1840 & $[-0.994,\,0.994]$\\
			3 &  21 & 1924 & $[-0.995,\,0.995]$\\
			\bottomrule
		\end{tabular*}
		\subcaption{$\ncells=500$, $\MEElements=20$, $\SGtruncorder=4$, $\xidelta=1e-2$, $tol=1e-5$
		}
		\label{tab:Kx}
	\end{subtable}
	\hspace*{0.15cm}
	\begin{subtable}[c]{0.48\textwidth}
		\begin{tabular*}{\textwidth}{C{.7cm}C{1cm}C{1.3cm}c}
			\multicolumn{4}{c}{}\\
			\multicolumn{4}{c}{}\\
			\toprule
			\begin{tabular}{@{}c@{}} \mbox{}\\\mbox{}\end{tabular} $\SGtruncorder$& $\#$ it. & time$[$s$]$ & $\uncertaintyInt$\\
			\midrule
			1 & 18 & 124 & $[-0.992,\,0.992]$\\
			2 & 21 & 215 & $[-0.993,\,0.993]$\\
			4 & 19 & 435 & $[-0.993,\,0.993]$\\
			8 & 18 & 996 & $[-0.994,\,0.994]$\\
			\bottomrule
		\end{tabular*}
		\subcaption{$\ncells=500$, $\MEElements=20$, $\xdegree\!=1$, $\xidelta=1e\!-\!2$, $tol=1e-5$
		}
		\label{tab:KO}
	\end{subtable}
	\centering
	\begin{subtable}[c]{0.4805\textwidth}
		\begin{tabular*}{\textwidth}{C{.9cm}C{1cm}C{1.3cm}c}
			\multicolumn{4}{c}{}\\
			\multicolumn{4}{c}{}\\
			\toprule
			\begin{tabular}{@{}c@{}} \mbox{}\\\mbox{}\end{tabular} $\xidelta$ & $\#$ it. & time $[$s$]$ & $\uncertaintyInt$\\
			\midrule
			1 & \multicolumn{3}{c}{fails, no direction of descent}\\
			$1e-1$ & 18 & 156 & $[-0.903,\,0.903]$\\
			$1e-2$ & 18 & 308 & $[-0.990,\,0.990]$\\
			$1e-3$ & 18 & 470 & $[-0.999,\,0.999]$\\
			0 & 40 & 355 & $[-1.000,\,1.000]$\\
			\bottomrule
		\end{tabular*}
		\subcaption{$\ncells=300$, $\MEElements=20$, $\xdegree=1$, $\SGtruncorder=4$, $tol=1e-5$
		}
		\label{tab:delta}
	\end{subtable}
	\caption{Number of iterations ($\#$ it.), run time in seconds (time$[$s$]$) and distribution parameters in the last iteration ($\uncertaintyInt$) for different parameters of \algoref{algo:OCalgo}. Example \eqref{eq:initial_sinus}.}
	\label{tab:params}
\end{table}

\subsection{Burgers' Equation}
In this subsection, we consider the Burgers' equation with uncertain initial conditions
\begin{subequations}\label{eq:Burger}
	\begin{align}
	\frac{\partial}{\partial\timevar}\solution(\timevar ,\x, \uncertainty)+ \frac12 \frac{\partial}{\partial\x} \solution(\timevar ,\x, \uncertainty)^2&= 0, \quad \quad \text{for} \ \ \x\in\Domain,\, \timevar\in [0,\,\tend],\,\uncertainty\in\SD,\\
	\solution(0,\x,\uncertainty) &= \sin(2\pi\x) + \frac12\uncertainty, \quad \quad \text{for} \ \ \x\in\PD,\, \xi\in\SD.
	\end{align}
\end{subequations}

\refone{We consider the solution of the Burgers' equation until $\tend=0.05$ with $\Domain=[0,\,1]$. The uncertainty is uniform $\uncertainty\sim\mathcal{U}(\uncertaintyL,\uncertaintyR)$, whereas $u_D$ is obtained by calculating the solution to \eqref{eq:Burger} at $\tend=0.05$ with a uniform distribution in $[-1,\,1]$.} Initially, we start with $\uncertaintyInt_p =[-0.8,\,0.8] $. We use $\ncells=200$ spatial cells and $\MEElements=40$ \ME{}s, the polynomial degrees of the local approximation are $\xdegree =2$ and $\SGtruncorder = 4$. Additionally, we prescribe the scaling factor $\xidelta=10^{-2}$ and the tolerance $tol=10^{-2}$.

In \figref{fig:burger_xi} -- \ref{fig:burger_adjoint} we demonstrate the results. 
Again, we observe convergence of $\uncertaintyL$ and $\uncertaintyR$ to the true value $[-1,\,1]$, whereas 5 iterations are now required to reach the tolerance $10^{-2}$. 
The cost functional decreases from around 1.4 to 0.2 in the first iteration and then approaches $7\cdot10^{-3}$. 
In this example 1 and 8 Armijo stages were performed, starting initially with $\alpha=0.5$. 
\figref{fig:burger_state} shows the local polynomial approximation \eqref{localMEDGDGApproximation} of the state 
in the whole $\x-\uncertainty$ plane at time $\tend=0.05$, demonstrating the propagation of the sinus wave and its translation through the uncertainty. We derive the corresponding discontinuous stochastic Galerkin polynomial from the adjoint state 
\begin{equation}\label{eq:SGadj}p\big|_{\cell{\cellind}\times D_\MEIndex}(\timevar,\x,\uncertainty) = \sum_{\xsumIndex=0}^{\xdegree}\sum_{\SGsumIndex=0}^\SGtruncorder p^{\xsumIndex,\SGsumIndex}_{\cellind,\MEIndex}(\timevar)\polybasis{\xsumIndex}_{\cellind}(\x)\localxiBasisPoly{\SGsumIndex}{\MEIndex}(\uncertainty)\end{equation}
for all $\cellind=1,\ldots,\ncells$ and $\MEIndex=1,\ldots,\MEElements$ and illustrate it in \figref{fig:burger_adjoint}. It is interesting to study in future work how this polynomial relates to an adjoint equation on PDE level. At time $\tend=0.05$ we obtain the terminal condition of the adjoint system, namely 
$M_1^{\xsumIndex,\SGsumIndex}\adj^{\xsumIndex,\SGsumIndex}(\tend) = (\solvec^{\xsumIndex,\SGsumIndex}_D-\solvec^{\xsumIndex,\SGsumIndex})$.

\externaltikz{burger}{\begin{figure}[htb]
	\centering
	\hspace*{-.4cm}
	\begin{minipage}[b]{0.45\textwidth}
 \begin{tikzpicture}
\begin{axis}[width=1.2\figurewidth, height=1.1\figureheight, xlabel=$\#$Iteration, ylabel=$\xi$, ylabel style = {rotate=-90}, ylabel style={at={(-0.1,0.5)}}, legend style={at={(0.97,0.7)}}] 
\addplot[mark=none, solid, tuklred, line width=.7pt] file {Images/burger/burger_xiL.txt};
\addplot[mark=none, dashed, tuklblue, line width=.8pt] file {Images/burger/burger_RefxiL.txt};
\addplot[mark=none, solid, tuklred, line width=.7pt] file {Images/burger/burger_xiR.txt};
\addplot[mark=none, dashed, tuklblue, line width=.8pt] file {Images/burger/burger_RefxiR.txt};
\legend{$\uncertaintyInt$, $\uncertaintyInt_{\text{Ref}}$}
\end{axis}
\end{tikzpicture}
 \caption{Performance of $\xi$ during the optimization process. Example \eqref{eq:Burger}.}
\label{fig:burger_xi}
	\end{minipage}
	\hspace*{.8cm}
	\begin{minipage}[b]{0.45\textwidth}
	\begin{tikzpicture}
\begin{axis}[width=1.2\figurewidth, height=1.1\figureheight, xlabel=$\#$Iteration, ylabel=$\costfun$, ylabel style = {rotate=-90}, ymode = log] 
\addplot[mark=none, solid, tuklblue, line width=.7pt] file {Images/burger/burger_J.txt};
\end{axis}
\end{tikzpicture}
	\caption{\refone{Cost functional during the optimization process. Example \eqref{eq:Burger}.}}
\label{fig:burger_xishock_J}
	\end{minipage}
\end{figure}


\begin{figure}[htb]
	\centering
	\hspace*{-.1cm}
	\begin{minipage}[b]{0.45\textwidth}		
		\begin{tikzpicture}
				\begin{groupplot}[
		group style={group size=1 by 1, horizontal sep = 1.5cm,  vertical sep = 2.5cm},
		axis on top,
		scale only axis,
		enlargelimits=false,
		width = 1\figurewidth,
		height = \figureheight,
		colormap name = {parula},
		colorbar horizontal,
		title style = {yshift = 0.6cm} ,
		scaled x ticks=false,
		colorbar style={,
			at={(0,1.02)},
			anchor=south west,
			height=0.03\textwidth,
			width=1\figurewidth,
			xticklabel style={font=\footnotesize,anchor=south, /pgf/number format/.cd, fixed,precision = 4,/tikz/.cd},
			xticklabel shift = -7pt,
		},
		]
		\nextgroupplot[
		xmin=0.000000,
		xmax=1.000000,
		ymin=-0.800000,
		ymax=0.800000,
		xlabel={\x},
		ylabel={\uncertainty},
		ylabel style = {rotate=-90},
		ylabel style={at={(-0.1,0.5)}},
		point meta min = -1.401656,
		point meta max = 1.401656,
		]
		\addplot graphics[xmin=0.000000,xmax=1.000000,ymin=-1.000000,ymax=1.000000] {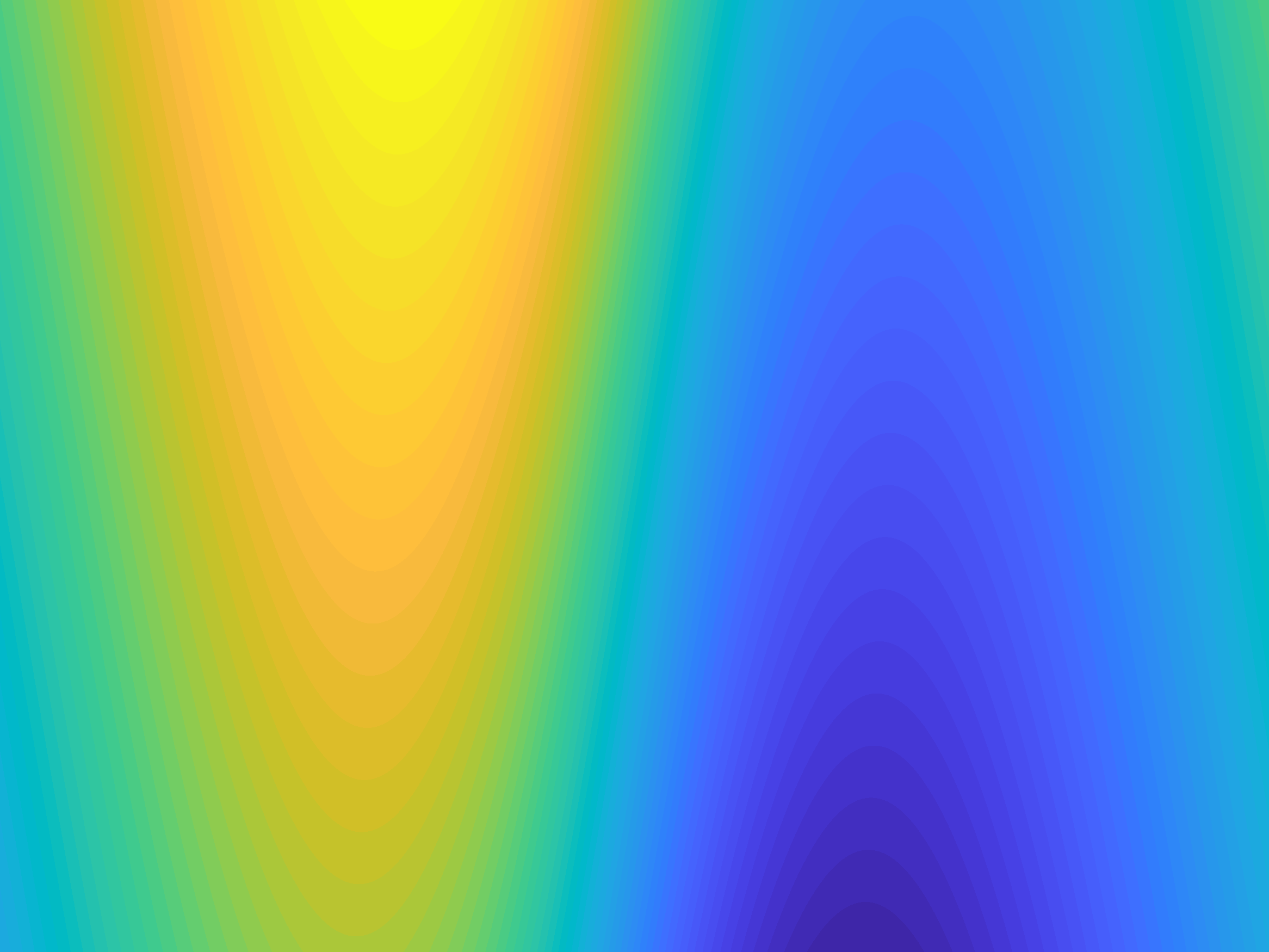};
		\end{groupplot}
		\end{tikzpicture}	
		\caption{\refone{Numerical solution $u$ to the conservation law within the $\x-\uncertainty$ plane at time $\tend=0.05$ in the first iteration, i.e., $\uncertaintyInt=\uncertaintyInt_p$. Example \eqref{eq:Burger}.}}	
		\label{fig:burger_state}
	\end{minipage}
\hspace*{1.0cm}
\begin{minipage}[b]{0.47\textwidth}		
	\begin{tikzpicture}
	\begin{groupplot}[
	group style={group size=1 by 1, horizontal sep = 1.5cm,  vertical sep = 2.5cm},
	axis on top,
	scale only axis,
	enlargelimits=false,
	width = 1\figurewidth,
	height = \figureheight,
	colormap name = {parula},
	colorbar horizontal,
	title style = {yshift = 0.6cm} ,
	scaled x ticks=false,
	colorbar style={,
		at={(0,1.02)},
		anchor=south west,
		height=0.03\textwidth,
		width=1\figurewidth,
		xticklabel style={font=\footnotesize,anchor=south, /pgf/number format/.cd, fixed,precision = 4,/tikz/.cd},
		xticklabel shift = -7pt,
	},
	]
	\nextgroupplot[
	xmin=0.000000,
	xmax=1.000000,
	ymin=-0.800000,
	ymax=0.800000,
	xlabel={\x},
	ylabel={\uncertainty},
	ylabel style = {rotate=-90},
	ylabel style={at={(-0.1,0.5)}},
	point meta min = -0.142241,
	point meta max = 0.142241,
	]
	\addplot graphics[xmin=0.000000,xmax=1.000000,ymin=-1.000000,ymax=1.000000] {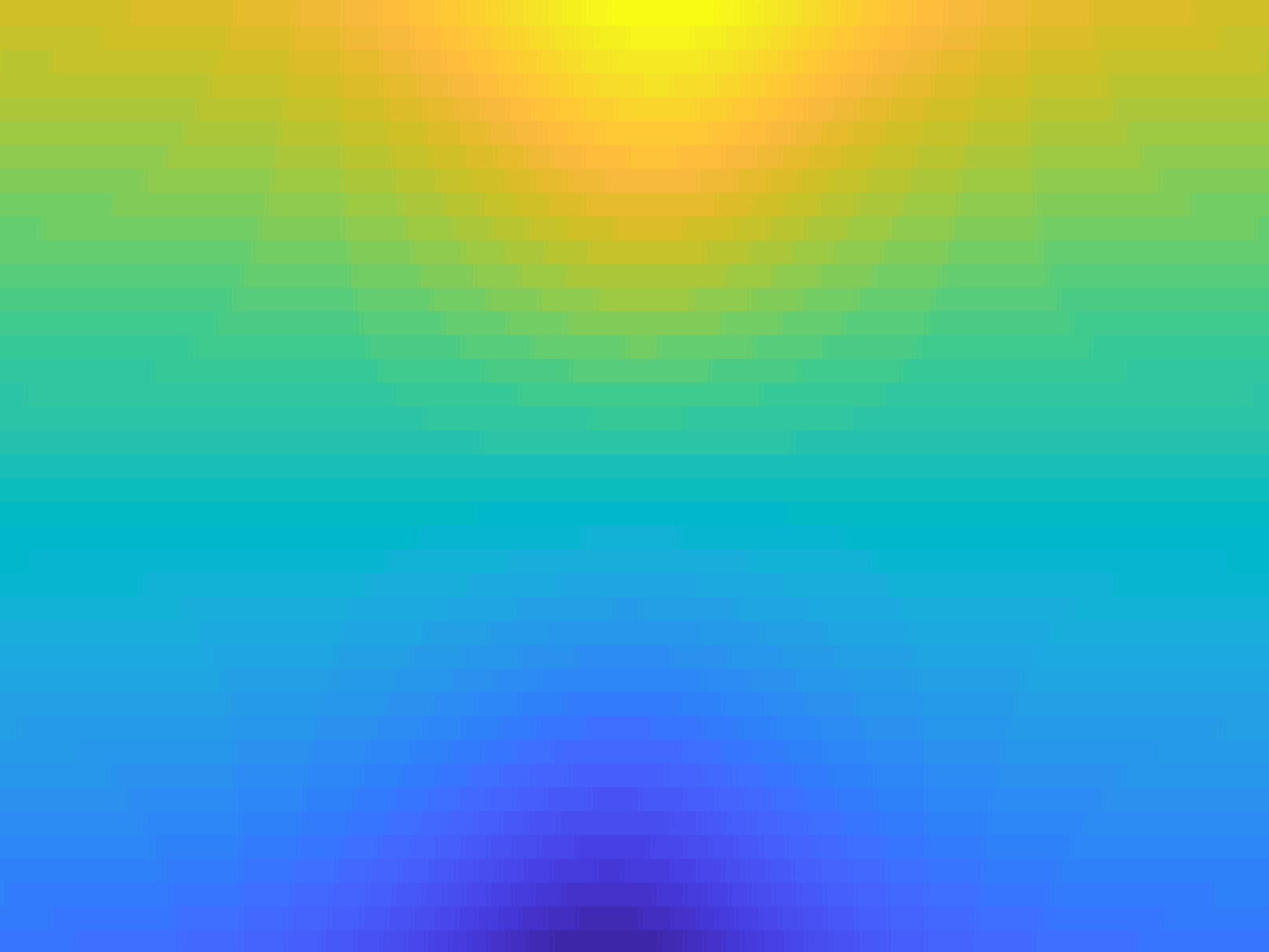};
	\end{groupplot}
	\end{tikzpicture}
	\caption{\refone{Discontinuous SG polynomial of the adjoint $p$, defined in \eqref{eq:SGadj},
		within the $\x-\uncertainty$ plane at time $\tend=0.05$ in the first iteration. Example \eqref{eq:Burger}}.}
	\label{fig:burger_adjoint}
\end{minipage}
\end{figure}

}

\section{Conclusions and Outlook}
We studied the estimation of parameters for the underlying random distribution in uncertain scalar conservation laws as an optimization problem which we solved on the time-continuous level. For the discretization in the physical and stochastic domain we employed the discontinuous stochastic Galerkin scheme, using a \ME{} stochastic Galerkin ansatz in the stochastic and a discontinuous Galerkin approach in the spatial variable. We then computed the first order optimality conditions of this semi-discrete system. This yield our algorithm for the parameter identification which we applied on different numerical test cases for the uncertain linear advection and Burgers' equation. We observed the convergence of the distribution parameters to a predefined reference solution and the independence of the number of iterations from the discretization parameters. Our results show that it is possible to use intrusive UQ schemes such as stochastic Galerkin type methods for the theoretical framework of parameter estimations in uncertain conservation laws. Moreover, we deduced an algorithm that uses a high-order discretization in space, time and the stochasticity in order to solve the optimization problem.

Future work should incorporate the estimation of parameters for multi-dimensional random variables in the flux and the initial state. In this context, other (non-)intrusive UQ methods might be relevant to adopt into the methodology since they usually outperform sG for high-dimensional stochastic domains. 

Moreover, it might be of interest to analyze how to handle for example uncertain initial states that are not differentiable with respect to the uncertainty. \refone{This article is intended to provide a first approach on how to address parameter estimation problems in the context of high-order intrusive UQ methods for hyperbolic conservation laws. Further studies might include an error analysis of the parameter choices within the algorithm in the context of \tabref{tab:params}.}
As soon as an adjoint framework for hyperbolic stochastic PDEs, i.e., in the spirit of \cite{schafe2}, is established, we can also study the relation between the discontinuous stochastic Galerkin polynomial of the adjoint coefficients and a possible adjoint equation on PDE level.

\reftwo{A comparison to the results obtained by the parameter identification problem that is constrained by the system of stochastic differential equations resulting from the space-discretization of the uncertain PDE is another interesting project for future work.}

\section*{Acknowledgements}
Funding by Deutsche Forschungsgemeinschaft (DFG) within the RTG GrK 1932 ``Stochastic Models for Innovations in the Engineering Science'' is gratefully acknowledged.

\section*{References}
\bibliographystyle{siam}
\bibliography{library,biblio,bibliography}

\end{document}